\documentclass[11pt]{article}
\usepackage{graphicx,wrapfig}
\usepackage{flafter}
\usepackage{leftidx}
\usepackage{mathrsfs}
\usepackage{amssymb,amsmath}%\eqref'
\usepackage{bbding}
\usepackage{fancyhdr}%µ÷ÓÃÅÅÓ¡Ò³Ã¼ÓëÒ³½ÅµÄºê°ü
\usepackage{mathrsfs}
\usepackage{array}
\usepackage{color}
\usepackage{stmaryrd}
\usepackage{amsfonts}
\usepackage{latexsym}
\usepackage{psfrag}
\usepackage{subfigure}
\usepackage{fancyhdr,graphicx}
\usepackage{multicol}
\usepackage{dsfont}
\usepackage{bbm}
\usepackage{booktabs}
\usepackage[center]{caption2}
\usepackage{cite}
\usepackage{epstopdf}
\usepackage{booktabs}%´´½¨Ã»ÓÐÊúÏß·Ö¸îµÄ±í¸ñ
\usepackage [latin1]{inputenc}% ·ÀÖ¹Ìá½»Ê±ÔËÐÐ³ö´í
\usepackage{multirow}%ºÏ²¢¶àÐÐµ¥Ôª¸ñ
\usepackage{enumerate}%ÁÐ±íºê°ü
\usepackage{enumitem}
\usepackage{algpseudocode,algorithm,algorithmicx}% Éú³ÉÎ±´úÂëËã·¨Í¼
\pdfoutput=1

\allowdisplaybreaks
\renewcommand{\leq}{\leqslant}

\date{}
\newtheorem{theorem}{Theorem}[section]
\newtheorem{lemma}{Lemma}[section]

\newtheorem{example}{Example}[section]

\numberwithin{equation}{section}% °´ÕÕÕÂ½Ú±àºÅ
\newcommand{\zd}{\,\mathrm{d}}
\newcommand{\abs}[1]{\left|#1\right|}

\newcommand{\bra}[1]{\left(#1\right)}
\newcommand{\brab}[1]{\big(#1\big)}

\newcommand{\kbra}[1]{\left[#1\right]}
\newcommand{\kbrab}[1]{\big[#1\big]}

\begin{document}
\title{Adaptive linear second-order energy stable schemes for time-fractional Allen-Cahn equation with volume constraint}
\author{Bingquan Ji\thanks{Department of Mathematics, Nanjing University of Aeronautics and Astronautics,
211101, P. R. China. Bingquan Ji (jibingquanm@163.com).}
\quad Hong-lin Liao\thanks{Corresponding author. ORCID 0000-0003-0777-6832; Department of Mathematics,
Nanjing University of Aeronautics and Astronautics,
Nanjing 211106, P. R. China. Hong-lin Liao (liaohl@csrc.ac.cn and liaohl@nuaa.edu.cn)
is supported by a grant 1008-56SYAH18037
from NUAA Scientific Research Starting Fund of Introduced Talent.}
\quad Yuezheng Gong \thanks{Department of Mathematics, Nanjing University of Aeronautics and Astronautics,
Nanjing 210016, P. R. China;
Yuezheng Gong (gongyuezheng@nuaa.edu.cn) is partially supported by the NSFC grant No. 11801269, and the NSF grant No. BK20180413 of Jiangsu Province.}
\quad Luming Zhang\thanks{Department of Mathematics, Nanjing University of Aeronautics and Astronautics,
211101, P. R. China. Luming Zhang (zhanglm@nuaa.edu.cn)
is supported by the NSFC grant No. 11571181.}}
%%%%%%%%%%%%%%%%%%%%%%%%%%%%%%%%%%%%%%%%%%%%%%%%%%%%%%%%%%%%%%%%%%%%%%%%%%%%%%%%%%%%%%%%%%
\date{}
\maketitle
\normalsize

\begin{abstract}
	A time-fractional Allen-Cahn equation with volume constraint is first proposed by introducing a nonlocal time-dependent Lagrange multiplier. Adaptive linear second-order energy stable schemes are developed for the proposed model by combining invariant energy quadratization and scalar auxiliary variable approaches with the recent L1$^{+}$ formula. The new developed methods are proved to be volume-preserving and unconditionally energy stable on arbitrary nonuniform time meshes. The accelerated algorithm and adaptive time strategy are employed in numerical implement. Numerical results show that the proposed algorithms are computationally efficient in multi-scale simulations,  and appropriate for accurately resolving the intrinsically initial singularity of solution and for efficiently capturing the fast dynamics away initial time.\\
  \noindent{\emph{Keywords}:}\;\; Time-fractional Allen-Cahn equation with volume constraint; invariant energy quadratization; scalar auxiliary variable; L1$^{+}$ formula; unconditional energy stable

  \noindent{\bf AMS subject classiffications.}\;\; 35Q99, 65M06, 65M12, 74A50
\end{abstract}
\section{Introduction}

%%%%%%%%%%%%%%%%%%%%%%%%%%%%%%%%%%%%%%%%%%%%%%%%%%%%%%%%%%%%%%%%%%%%%
The gradient flow models are frequently used to describe relaxation dynamics that obey the second law of thermodynamics, ranging from materials science, fluid dynamics and engineering \cite{Allen1979A,Cahn1958Free,Gong2018lLnear}. One of well-known models is the Allen-Cahn equation, which was originally introduced to model the anti-phase domain coarsening in a binary alloy \cite{Allen1979A}. Also, in the past decades, the Allen-Cahn equation and its various variants have been applied for a wide range of phenomena due to its advantages for microstructure numerical simulations, for instance, grain growth \cite{Krill2002Computer} and crystal growth \cite{Li2012Phase-field}. However, considering the phase variable represents the volume fraction of material component, the classical Allen-Cahn equation does not conserve the initial volume. To fix this drawback, the first work was given by Rubinstein and Sternaberg, who added a time-dependent Lagrange multiplier to the original equation arising from an enforcement of conservation of volume \cite{Rbinstein1992Nonlocal}. Brassel and Bretin introduced another remedy to preserve the total volume-conservative property, i.e., they imposed local and nonlocal effects on the primitive model \cite{Brassel2011A}. Recently, the time, space and time-space fractional Allen-Cahn equations were suggested to accurately describe anomalous diffusion problems \cite{Hou2017Numerical,Zheng2017A,Liu2018Time}. However, they don't preserve the volume conservation. In this paper, we are going to develop a new time-fractional Allen-Cahn equation by enforcing a nonlocal Lagrange multiplier to cancel out the variation of volume, while without influencing the primitive energy dissipative property.

An alternative model for the gradient flow system is the Cahn-Hilliard equation, which naturally possesses the volume-preserving property \cite{Cahn1958Free}.  The Allen-Cahn model with a volume constraint has been studied and compared with the Cahn-Hilliard model in \cite{Rbinstein1992Nonlocal}. The authors suggested that the Allen-Cahn model with a volume constraint is more appropriate for simulating the interfacial dynamics of immiscible multi-component material systems. And the order of Allen-Cahn equation is substantially lower than that of the Cahn-Hilliard equation, which implies that it may be relatively easier to simulate numerically. Some interesting insight may be offered by carrying out comparison investigations of volume conservative phased field models \cite{Lee2016High-order, Lee2016Comparison}.

There have been a great amount of works to develop energy stable schemes for the gradient flow model. The early well-known numerical approaches include the convex-splitting technique and the stabilizing method. Readers are referred to \cite{Shen2012Second,Xu2006Stability} for more details. Recently, Yang et al. proposed a new numerical idea of recasting the free energy into a quadratic functional to design linear, second-order, unconditionally energy stable schemes, which called the invariant energy quadratization (IEQ) method \cite{Yang2017numerical}. Subsequently, Shen et al. developed the scalar auxiliary variable (SAV) approach, which was shown to be more effective than the IEQ approach \cite{Jie2018The}. In fact, the common goal of IEQ and SAV strategies is to first transform the original PDE system into a new equivalent system with a quadratic energy functional and the corresponding modified energy dissipation law. Specifically, applying the energy stable algorithms derived by the two energy quadratization strategies, the volume-preserving Allen-Cahn model was compared with the classical Allen-Cahn model as well as the Cahn-Hilliard model  \cite{Jing2019Second}. For more discussions and the applications of IEQ and SAV strategies, we refer to \cite{Yang2017numerical,Jie2018The,Gong2019Energy,Jing2019Second} and  the references therein.

Along the numerical front with respect to the fractional phase field models, there are a lot of works devoted to the investigation on the solutions of the nonlocal models. Precisely, Hou et al. \cite{Hou2017Numerical} showed that the space-fractional Allen-Cahn equation could be
viewed a $L^{2}$ gradient flow for the fractional analogue version of Ginzburg-Landau free energy function. Meanwhile, the authors proved that the proposed numerical scheme preserves the energy decay property and the maximum principle in the discrete level. Li et al. \cite{Zheng2017A} investigated a space-time fractional Allen-Cahn phase field model that describes the transport of the fluid mixture of two immiscible fluid phases. They concluded that the alternative model could provide more accurate description of anomalous diffusion processes and sharper interfaces than the classical model. The first theoretical contribution regarding the energy dissipation property of the time-fractional phase models was done by Tang et al. \cite{Tang2018On}. They proved that the time-fractional phase field models indeed admit an energy dissipation law of an integral type. In addition, they applied the uniform L1 formula to construct a class of finite difference schemes,
which can preserve the theoretical energy dissipation property. Very recently, Du et al. \cite{Du2019Time} studied the time-fractional Allen-Cahn equation, where the well-posedness, solution regularity, and maximum principle were proved rigorously. In addition, several unconditionally solvable and stable time-stepping schemes were developed. Also, the related convergence of those numerical approaches were established without any extra regularity assumption on the exact solution. Zhao et al. \cite{Liu2018Time,Zhao2019On} studied a series of the time-fractional phase field models numerically, including the time-fractional Cahn-Hilliard equation with different types of variable mobilities and  time-fractional molecular beam epitaxy model. The considerable numerical evidences indicate that the effective free energy or roughness of the time-fractional phase field models during coarsening obeys a similar power scaling law as the integer ones, where the power is linearly proportional to the fractional index $\alpha$. In other words, the main difference between the time-fractional phase field models and integer ones lies in the time-scales of coarsening.

In this paper, we first apply the IEQ/SAV approaches to reformulate the time-fractional phase field models into an equivalent system. Then the nonuniform L1$^{+}$ formula proposed in \cite{Ji2019Adaptive} is applied for the equivalent time-fractional model to develop linear, second-order energy stable numerical schemes, which are proved to preserve the volume conservation law and unconditionally energy stability on
arbitrary nonuniform time meshes. Since
the solution lacks the smoothness near the initial time although it would be smooth away from  $t=0$ \cite{Jin2016An,JinLiZhou:2017}, the predicted second-order time accuracy of L1$^{+}$ formula is always restrictive. Actually, in any numerical methods for solving time-fractional diffusion equations, a basic consideration is the initial singularity of solution, see the recent works \cite{Liao2018Sharp,Liao2018Unconditional,Liao2018second}. Based upon the realistic assumptions on the exact solution, we utilize the L1$^{+}$ formula on nonuniform time steps
to compensate the intrinsically weak singularity of time-fractional models near initial time. We will show that the graded mesh can recover the optimal time accuracy when the solution is non-smooth near $t=0$ numerically. In addition, in order to overcome the global dependence of historical solutions of time-fractional Caputo derivative, a fast variant of  L1$^{+}$ formula is used to significantly reduce the computational complexity and the storage requirements. Since the evolution of time-fractional phase field models involves multiple time scales, adaptive time step strategy based on the evolution of total energy are reported to efficiently resolve widely varying time scales.

The outline of the article is arranged as follows. The time-fractional Allen-Cahn equation and its volume-conserving version as well as the  time-fractional Cahn-Hilliard model are reported
in Section 2. We then present the corresponding energy stable numerical schemes in Section 3.
In Section 4, several numerical examples are performed to confirm the theoretical findings, covering the volume conservation and energy dissipation properties, and provide new insights on the volume-conservative time-fractional Allen-Cahn equation compared with the non-volume-preserving one and Cahn-Hilliard equation. %In closing, we give some concluding remarks.

\section{Time-fractional phase field models}
Introduce a phase variable $\phi$,
for the effective free energy  of the phase model $E[\phi]$,
\begin{align}
E[\phi]=\int_{\Omega}\bra{\frac{\varepsilon^{2}}{2}|\nabla\phi|^{2}
+F(\phi)}\zd{\mathbf{x}},\quad
F(\phi)=\frac{1}{4}(1-\phi^{2})^{2},
\end{align}
in which $\varepsilon$ is a parameter describing the width of the interface,
the time-fractional Allen-Cahn equation then reads
\begin{align}\label{Problem-1}
\partial_{t}^{\alpha}\phi
=-\lambda\frac{\delta{E}}{\delta\phi},
\end{align}
where positive constant $\lambda$ is the mobility parameter,
$\frac{\delta{E}}{\delta\phi}$ is the functional derivative of $E$ with respect to phase variable $\phi$.
Here, the notation $\partial_{t}^{\alpha}:={}_{0}^{C}\!D_{t}^{\alpha}$ in {\eqref{Problem-1}} denotes
the Caputo's fractional  derivative of order $\alpha$ with respect to $t$, i.e.,
\begin{align}\label{CaputoDef}
(\partial_{t}^{\alpha}v)(t)
:=(\mathcal{I}_{t}^{1-\alpha}v')(t)
=\int_{0}^{t}\omega_{1-\alpha}(t-s)v'(s)\zd{s},\quad 0<\alpha<1,
\end{align}
involving the fractional  Riemann-Liouville integral $\mathcal{I}_{t}^{\beta}$ of order $\beta>0$, that is,
\begin{align}
(\mathcal{I}_{t}^{\beta}v)(t)
:=\int_{0}^{t}\omega_{\beta}(t-s)v(s)\zd{s},\quad\text{where}\quad  \omega_{\beta}(t):=t^{\beta-1}/\Gamma(\beta).
\end{align}
It is remarkable that,
in comparison with the energy dissipation law of the local Allen-Cahn model,
Tang et al. \cite{Tang2018On} proved that the energy stable property of the nonlocal one
is given by,
\begin{align}\label{Frac-Energy-Decay-Law}
E\kbra{\phi(T)}
-E\kbra{\phi(0)}
=-\frac{1}{\lambda}
\int_{\Omega}
\mathcal{I}_{t}^{1}
(\partial_{t}\phi\,\mathcal{I}_{t}^{1-\alpha}\partial_{t}\phi)(T)
\zd{\mathbf{x}}
\leq{0}.
\end{align}
The non-positive of the right part of above relation is determined by
\cite[Lemma 2.1]{Tang2018On}.

Evidently, acting the Riemann-Liouville fractional derivative ${}_{0}^{RL}\!D_{t}^{1-\alpha}$ on both sides of equation \eqref{Problem-1}, and using the identity
${}_{0}^{RL}\!D_{t}^{1-\alpha}{}_{0}^{C}\!D_{t}^{\alpha}v(s)=v^{\prime}(s)$,
one has the following relation
\begin{align}\label{Non-Volume}
\frac{\zd}{\zd{t}}\int_{\Omega}\phi\zd{\mathbf{x}}
=-\lambda{}_{0}^{RL}\!D_{t}^{1-\alpha}
\int_{\Omega}\frac{\delta{E}}{\delta\phi}\zd{\mathbf{x}}
\neq{0},
\end{align}
which means the time-fractional Allen-Cahn equation does not preserve the initial volume
that is consistent with the integer order one.
In order to impose the conservation of volume
\begin{align}\label{Voluma-Conservation}
\int_{\Omega}\phi(\mathbf{x},t)\zd{\mathbf{x}}
=\int_{\Omega}\phi(\mathbf{x},0)\zd{\mathbf{x}},
\end{align}
and without influencing  the original energy stable property \eqref{Frac-Energy-Decay-Law},
inspired by the volume conservation integer order Allen-Cahn equation performed in \cite{Rbinstein1992Nonlocal},
the equation {\eqref{Problem-1}} is modified by adding a nonlocal time-dependent Lagrange multiplier $\eta(t)$ as follows
\begin{align}\label{Problem-2}
\partial_{t}^{\alpha}\phi
=-\lambda\bra{\frac{\delta{E}}{\delta\phi}-\eta(t)},
\end{align}
where the expression of the new term is given by $\eta(t)=\frac{1}{\abs{\Omega}}\int_{\Omega}\frac{\delta{E}}{\delta\phi}\zd{\mathbf{x}}$
for the necessary condition to guarantee the invariant volume, i.e.,
\begin{align}\label{Invariant-Volume-Condition}
\frac{\zd}{\zd{t}}\int_{\Omega}\phi\zd{\mathbf{x}}
=-\lambda{}_{0}^{RL}\!D_{t}^{1-\alpha}
\int_{\Omega}\bra{\frac{\delta{E}}{\delta\phi}-\eta(t)}\zd{\mathbf{x}}
=0.
\end{align}

As well-known,
based upon the free energy $E[\phi]$,
another model that maintains the initial volume is
the time-fractional Cahn-Hilliard equation \cite{Tang2018On},
\begin{align}\label{Problem-3}
\partial_{t}^{\alpha}\phi
=\lambda\Delta\frac{\delta{E}}{\delta\phi}.
\end{align}
After a small calculation analogous to the derivation of \eqref{Non-Volume},
we see that
\begin{align}
\frac{\zd}{\zd{t}}\int_{\Omega}\phi\zd{\mathbf{x}}
=\lambda{}_{0}^{RL}\!D_{t}^{1-\alpha}
\int_{\Omega}\Delta\frac{\delta{E}}{\delta\phi}\zd{\mathbf{x}}
%=\lambda{}_{0}^{RL}\!D_{t}^{1-\alpha}
%\int_{\partial\Omega}\nabla\frac{\delta{E}}{\delta\phi}
%\cdot\mathbf{n}\zd{S}
=0,
\end{align}
where the periodic boundary condition is chosen to ensure that the boundary integrals vanish.
Meanwhile, the model \eqref{Problem-3} preserves the energy stable property,
\begin{align}
E\kbra{\phi(T)}
-E\kbra{\phi(0)}
=-\frac{1}{\lambda}
\int_{\Omega}
\mathcal{I}_{t}^{1}
(\nabla\psi\,\mathcal{I}_{t}^{1-\alpha}\nabla\psi)(T)
\zd{\mathbf{x}}
\leq{0},
\end{align}
in which $\psi=-\Delta^{-1}\partial_{t}\phi$ is the solution of the equation
$-\Delta\psi=\partial_{t}\phi$ with periodic boundary condition,
see \cite[Lemma 2.3]{Tang2018On} for more details.

To our knowledge, for the above time-fractional phase field models,
there are limited results in the literature on numerical approaches
preserving the discrete volume conservation
as well as energy dissipation law,
especially on nonuniform time grids.
Therefore the first objective
of this paper is to build  nonuniform time-stepping methods
for the continuous systems to inherit the corresponding
invariant or dissipative properties enjoyed by the original systems.

We consider the nonuniform time levels
$0=t_{0}<t_{1}<\cdots<t_{k-1}<t_{k}<\cdots<t_{N}=T$
with the time-step sizes $\tau_{k}:=t_{k}-t_{k-1}$ for $1\leq{k}\leq{N}$
and the maximum time-step size $\tau:=\max_{1\leq{k}\leq{N}}\tau_{k}$.
Also, let the local time-step ratio $\rho_k:=\tau_k/\tau_{k+1}$ and the maximum step ratio $\rho:=\max_{k\geq 1}\rho_k$. Given a grid function $\{v^{k}\}$,
put $\triangledown_{\tau}v^{k}:=v^{k}-v^{k-1}$, $\partial_{\tau}v^{k-\frac12}:=\triangledown_{\tau}v^{k}/\tau_k$
and $v^{k-\frac{1}{2}}:=(v^{k}+v^{k-1})/2$ for $k\geq{1}$.
Always,
let $(\Pi_{1,k}v)(t)$ denote the linear interpolant of a function $v(t)$ at two nodes $t_{k-1}$ and $t_{k}$,
and define a piecewise linear approximation
\begin{align}\label{linear interpolation}
\Pi_{1}v:=\Pi_{1,k}v\quad
\text{so that}\quad
(\Pi_{1}v)'(t)=\partial_{\tau}v^{k-\frac12}\quad
\text{for $t_{k-1}<{t}\leq t_{k}$ and $k\geq1$}.
\end{align}
%%%%%%%%%%%%%%%%%%%%%%%%%%%%%%%%%%%%%%%%%%%%%%%%%%%%%%%%%%%%%%%%%%%%%%%%%%%%
%%%%%%%%%%%%%%%%%%%%%%%%%%%%%%%%%%%%%%%%%%%%%%%%%%%%%%%%%%%%%%%%%%%%%%%%%%%%
\begin{figure}[htb!]
\centering
\includegraphics[width=3.0in,height=2.0in]{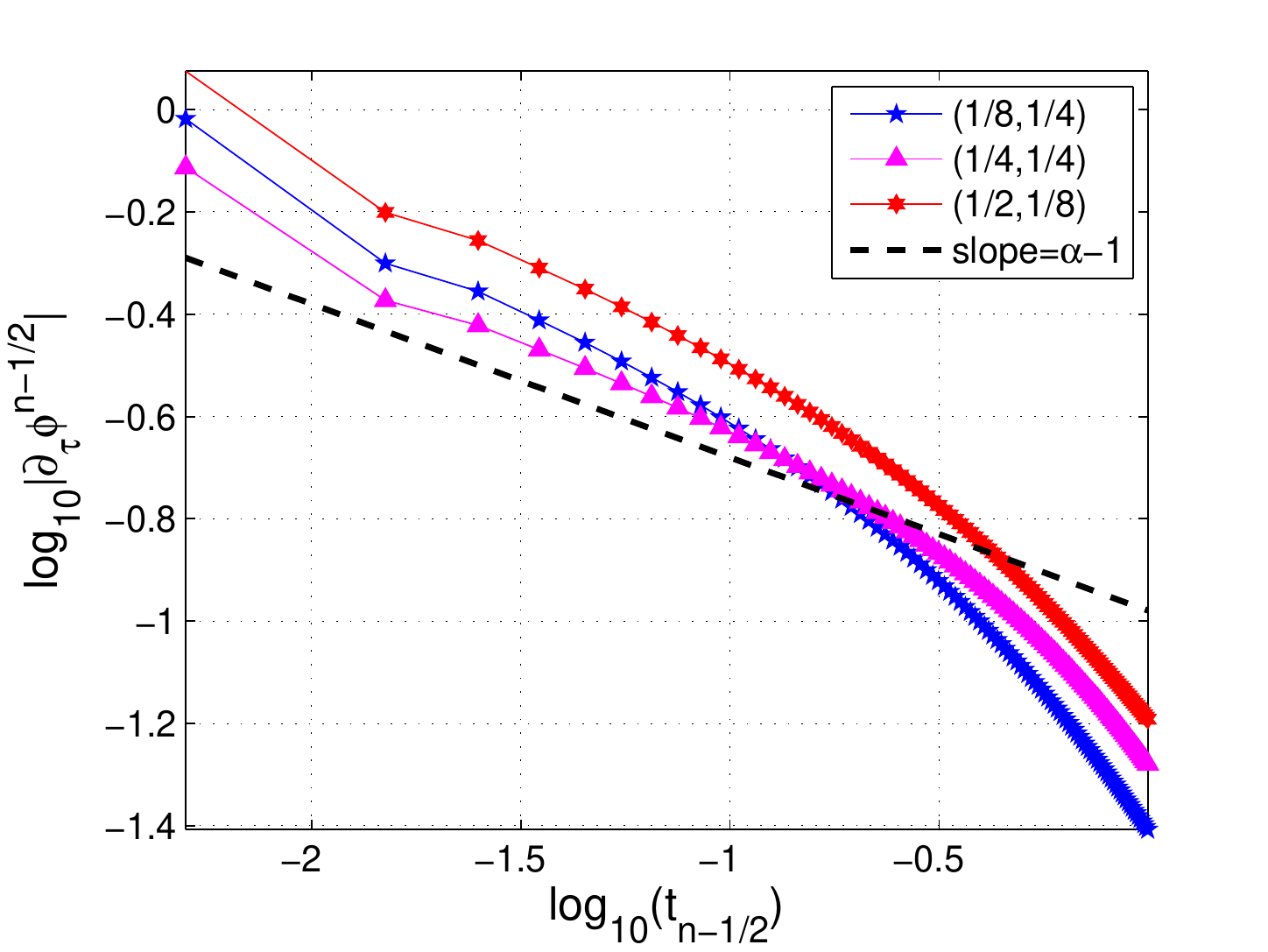}
\includegraphics[width=3.0in,height=2.0in]{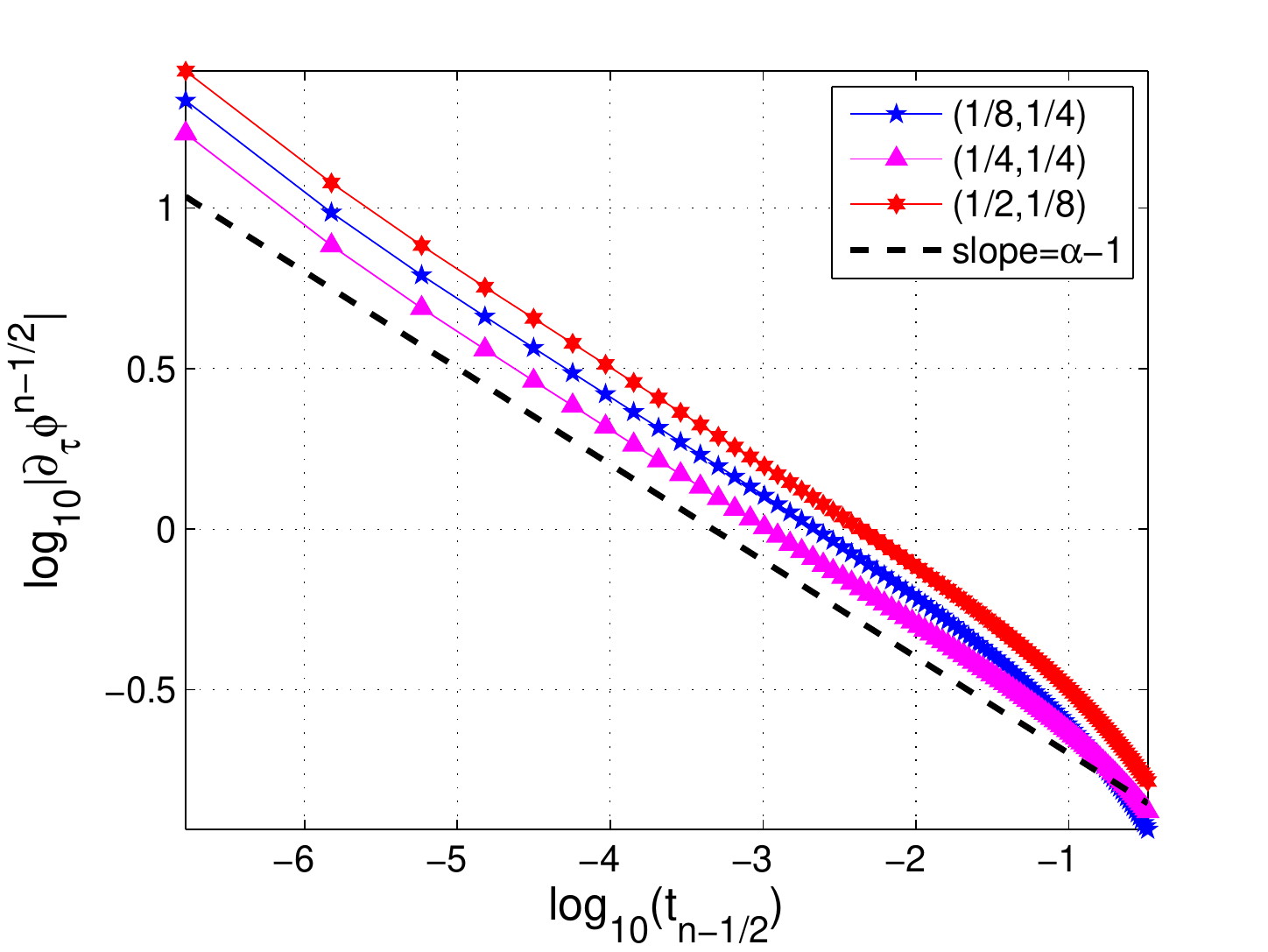}
\caption{The log-log plot of the difference quotient $\partial_{\tau}\phi^{k-\frac12}$ versus time
for problem {\eqref{Problem-2}} with fractional order $\alpha=0.7$ and $\gamma=1,\,3$ (from left to right), respectively.}
\label{Nonlocal-AC-Singularity-alpha-07}
\end{figure}
%%%%%%%%%%%%%%%%%%%%%%%%%%%%%%%%%%%%%%%%%%%%%%%%%%%%%%%%%%%%%%%%%%%%%%%%%%%%
%%%%%%%%%%%%%%%%%%%%%%%%%%%%%%%%%%%%%%%%%%%%%%%%%%%%%%%%%%%%%%%%%%%%%%%%%%%%

To reveal the initial singularity of solution of
the time-fractional phase filed models,
we apply the L1$^{+}$ formula
that we see shortly
to the time-fractional problem {\eqref{Problem-2}},
and more details can be found in subsection 3.1 and Example
\ref{Simulating-Four-Drops}.
The drawings in Figure
\ref{Nonlocal-AC-Singularity-alpha-07}
depict the discrete time derivative $\partial_{\tau}\phi^{k-\frac12}$
near $t=0$ on the graded mesh $t_k=(k/N)^{\gamma}$ when fractional order $\alpha=0.7$.
The numerical results suggest that
\[
\log|\phi_{t}(\mathbf{x},t)|
\approx(\alpha-1)\log(t)+C(\mathbf{x})\quad\text{as $t\rightarrow0$,}
\]
and tell us
that the solution
is weakly singularity like $\phi_{t}=O(t^{\alpha-1})$ near initial time,
which could be alleviated by using the graded mesh.
Hence, the second objective of present work is to resolve the essentially weak singularity
in the time-fractional phase field
by refining time mesh.

\section{Energy stable numerical approaches}

To achieve the above assertions,
our starting point is to apply the L1$^{+}$
formula to approximate the Caputo derivative,
which naturally possesses the energy dissipation property on nonuniform time levels
when it is applied to the time-fractional phased filed models.

\subsection{The L1$^{+}$ formula of Caputo derivative}

The L1$^{+}$ formula for the Caputo derivative \eqref{CaputoDef} is defined at time $t=t_{n-\frac{1}{2}}$ as follows
\begin{align}\label{New-L1-Formula}
(\partial_{\tau}^{\alpha}v)^{n-\frac{1}{2}}
:=\frac{1}{\tau_{n}}\int_{t_{n-1}}^{t_{n}}
\int_{0}^{t}\omega_{1-\alpha}(t-s)(\Pi_{1}v)'(s)\zd{s}\zd{t}
=\sum_{k=1}^{n}a_{n-k}^{(n)}\triangledown_{\tau}v^{k}\quad \text{for $n\ge1$,}
\end{align}
in which the discrete convolution kernels $a_{n-k}^{(n)}$ are given by
\begin{align}\label{New-L1-Coeff}
a_{n-k}^{(n)}
:=\frac{1}{\tau_{n}\tau_{k}}\int_{t_{n-1}}^{t_{n}}
\int_{t_{k-1}}^{\min\{t,t_{k}\}}\omega_{1-\alpha}(t-s)\zd{s}\zd{t}\quad \text{for $1\leq{k}\leq{n}$.}
\end{align}
Following the discussions given in \cite[Lemma 3.1]{Ji2019Adaptive},
we have the following remarkable property. It says that the L1$^{+}$
formula is positive semi-definite on arbitrary nonuniform meshes.
%%%%%%%%%%%%%%%%%%%%%%%%%%%%%%%%%%%%%%%%%%%%%%%%%%%%%%%%%%%%%%%%%%%%%%%%%%%%
%%%%%%%%%%%%%%%%%%%%%%%%%%%%%%%%%%%%%%%%%%%%%%%%%%%%%%%%%%%%%%%%%%%%%%%%%%%%
\begin{lemma}\label{Quadratic-Form-New-L1}
The discrete convolution kernels $a_{n-k}^{(n)}$
in \eqref{New-L1-Coeff} are positive and positive semi-definite. For any real sequence $\{w_k\}_{k=1}^n$ with $n$ entries, it holds that
\begin{align*}
\sum_{k=1}^nw_k\sum_{j=1}^ka_{k-j}^{(k)}w_j\geq{0}\quad \text{for $n\ge1$.}
\end{align*}
\end{lemma}

The definition \eqref{New-L1-Coeff} of discrete kernels $a_{j}^{(n)}$ and the integral mean-value theorem
yield the following result.
%%%%%%%%%%%%%%%%%%%%%%%%%%%%%%%%%%%%%%%%%%%%%%%%%%%%%%%%%%%%%%%%%%%%%%%%%%%%%%%%%%%%%%%%%%%
\begin{lemma}\label{New-L1-Coeff-Relation}
The positive discrete kernels $a_{n-k}^{(n)}$ in \eqref{New-L1-Coeff} fulfill
\[
a_{0}^{(n)}=\frac{1}{\Gamma(3-\alpha)\tau_{n}^{\alpha}},\quad
a_{1}^{(n)}>a_{2}^{(n)}>\cdots>a_{n-1}^{(n)}>0\quad \text{for $n\ge2$}.
\]
\end{lemma}

Simple manipulations of the first two discrete kernels reveal that
\begin{align*}
a_{0}^{(n)}-a_{1}^{(n)}
=\frac{1}{\Gamma(3-\alpha)\tau_{n}^{\alpha}\rho_{n-1}}
\brab{1+\rho_{n-1}+\rho_{n-1}^{2-\alpha}-(1+\rho_{n-1})^{2-\alpha}}.
\end{align*}
It is easily seen that $a_{0}^{(n)}<a_{1}^{(n)}$ as $\alpha\rightarrow{0}$ and $a_{0}^{(n)}>a_{1}^{(n)}$ as $\alpha\rightarrow{1}$,
that is, the value of $a_{0}^{(n)}-a_{1}^{(n)}$ may change the sign when the fractional order $\alpha$ varies over $(0,1)$.

It is to mention that, the nonuniform L1$^{+}$ formula is quite different from some nonuniform formulas
approximating the Caputo time derivative,
including the L1 formula \cite{Liao2018Sharp,Liao2018Unconditional},
 L1-2$_{\sigma}$ (Alikhanov) formula \cite{Liao2016JSC,Liao2018second},
and Caputo's BDF2-type formula \cite{Liao2016Stability}.
We compare them in Table \ref{Comparison-Formula}, in which the discrete kernels are referred to
the sequence $\big\{A_{n-k}^{(n)}\big\}$ in the form $\sum_{k=1}^nA_{n-k}^{(n)}\triangledown_{\tau}v^n$
or $\sum_{k=1}^nA_{k}^{(n)}\triangledown_{\tau}v^{n-k}$.
As seen, the  L1$^{+}$  formula has some advantages: it is second-order accuracy,
the convergence order is independent of the fractional order $\alpha$,
and is positive semi-definite in the sense of Lemma \ref{Quadratic-Form-New-L1}.
Actually, these properties make it useful in designing
linear, second-order energy stable schemes to the time-fractional phase filed models introduced in Section 2.

 %%%%%table2%%%%%%%%%%%%%%%%%%%%%%%%%%%%%%%%%%%%%%%%%%%%%%%%%%%%%%%%%%%%%%%%%%%	
\begin{table}[htb!]
\begin{center}
\caption{Numerical Caputo derivatives on nonuniform meshes.}\label{Comparison-Formula} \vspace*{0.3pt}
\def\temptablewidth{0.87\textwidth}
{\rule{\temptablewidth}{0.5pt}}
\begin{tabular*}{\temptablewidth}{@{\extracolsep{\fill}}c|c|c|c|c}
  Numerical Formula
  &L1\cite{Liao2018Sharp,Liao2018Unconditional}
  &Alikhanov\cite{Liao2016JSC,Liao2018second}
  &CBDF2\cite{Liao2016Stability}
  &L1$^{+}$ \cite{Ji2019Adaptive}\\
  \midrule
  Formal accuracy          &$2-\alpha$ &$3-\alpha$    &$3-\alpha$           &2\\

  Positive kernels         &Yes        &Yes           &$A_1^{(n)}\ngeq 0$   &Yes\\

  Monotonous kernels       &Yes        &Yes           &$A_1^{(n)}\ngeq A_2^{(n)}$        &$A_0^{(n)}\ngeq A_1^{(n)}$\\

  Positive semi-definite   &Unknown    &Unknown       &Unknown              &Yes\\
\end{tabular*}
{\rule{\temptablewidth}{0.5pt}}
\end{center}
\end{table}	

\subsection{Numerical approach using IEQ}

For the volume-conserving time-fractional Allen-Cahn model {\eqref{Problem-2}},
we introduce an auxiliary function $u(\phi)$ in term of original variable $\phi$ given by
\begin{align}
u(\phi)
=\phi^{2}-1-\beta,
\end{align}
where the artificial parameter $\beta$ is utilized to regularize the numerical approaches.
As a consequence, the free energy of the original problem is transformed into a quadratic form
\begin{align}\label{IEQ-Equivalent-Energy}
E\kbra{\phi,u}
=\int_{\Omega}\bra{
\frac{\varepsilon^{2}}{2}\abs{\nabla\phi}^{2}
+\frac{\beta}{2}\abs{\phi}^{2}
+\frac{1}{4}u^{2}}
\zd{\mathbf{x}}
-\brab{\frac{\beta}{2}+\frac{\beta^{2}}{4}}\abs{\Omega}.
\end{align}
Correspondingly, the problem {\eqref{Problem-2}} could be reformulated to the following equivalent form
\begin{align}
\partial_{t}^{\alpha}\phi
&=-\lambda
\bra{
-\varepsilon^{2}\Delta\phi
+\beta\phi
+u\phi-\eta},\label{Nonlocal-AC-IEQ-1}\\
\eta
&=\frac{1}{\abs{\Omega}}\int_{\Omega}\bra{
-\varepsilon^{2}\Delta\phi
+\beta\phi
+u\phi}\zd{\mathbf{x}},\label{Nonlocal-AC-IEQ-2}\\
\partial_{t}u
&=2\phi\partial_{t}\phi.\label{Nonlocal-AC-IEQ-3}
\end{align}
The new system is subjected to the initial conditions
\begin{align}
\phi\bra{\mathbf{x},0}
=\phi_{0}\bra{\mathbf{x}}\quad
\text{and}\quad
u\bra{0}
=u\bra{\phi_{0}\bra{\mathbf{x}}},
\end{align}
and the same boundary conditions of the primitive model.
Define the usual $L^{2}$ inner product
$\bra{f,g}=\int_{\Omega}fg\zd{\mathbf{x}}$
for all $f,g\in{L}^{2}(\Omega)$.
We see clearly that the equivalent system
preserves  the volume-preserving property \eqref{Voluma-Conservation}
by making the $L^{2}$ inner product of \eqref{Nonlocal-AC-IEQ-1} with a constant
and finding that $(\partial_{t}^{\alpha}\phi,1)=0$.
Also, taking the inner product of \eqref{Nonlocal-AC-IEQ-1} and \eqref{Nonlocal-AC-IEQ-3} with $\partial_{t}\phi$ and $u$ respectively,
summing up the resulting equalities,
and integrating the time $t$ from $t=0$ to $T$, we obtain the energy decay law
\begin{align}
E\kbra{\phi(T),u(T)}
-E\kbra{\phi(0),u(0)}
=-\frac{1}{\lambda}
\int_{\Omega}
\mathcal{I}_{t}^{1}
\brab{\partial_{t}\phi\,\mathcal{I}_{t}^{1-\alpha}\partial_{t}\phi}(T)
\zd{\mathbf{x}}
\leq{0},
\end{align}
where we use the fact $(\eta,\partial_{t}\phi)=0$ due to the condition \eqref{Invariant-Volume-Condition}.

By virtue of the equivalent PDE system
\eqref{Nonlocal-AC-IEQ-1}-\eqref{Nonlocal-AC-IEQ-3},
we construct new numerical schemes
that concern only with the time discretization,
while the spatial approximations can be diverse,
examples as finite difference,
finite element or spectral methods.
Integrating the equations \eqref{Nonlocal-AC-IEQ-1}-\eqref{Nonlocal-AC-IEQ-3} from $t=t_{n-1}$ to $t_{n}$,
respectively, results in the following equations
\begin{align}
\frac{1}{\tau_{n}}\int_{t_{n-1}}^{t_{n}}
\partial_{t}^{\alpha}\phi
\zd{t}
&=-\frac{\lambda}{\tau_{n}}\int_{t_{n-1}}^{t_{n}}
\bra{
-\varepsilon^{2}\Delta\phi
+\beta\phi
+u\phi-\eta}\zd{t},\\
\frac{1}{\tau_{n}}\int_{t_{n-1}}^{t_{n}}\eta\zd{t}
&=\frac{1}{\abs{\Omega}\tau_{n}}
\int_{t_{n-1}}^{t_{n}}
\int_{\Omega}\bra{
-\varepsilon^{2}\Delta\phi
+\beta\phi
+u\phi}\zd{\mathbf{x}}\zd{t},\\
\frac{1}{\tau_{n}}\int_{t_{n-1}}^{t_{n}}\partial_{t}u\zd{t}
&=\frac{2}{\tau_{n}}\int_{t_{n-1}}^{t_{n}}
\phi\partial_{t}\phi\zd{t}.
\end{align}
By means of the L1$^{+}$ formula \eqref{New-L1-Formula},
the trapezoidal formula,
we have the following Crank-Nicolson IEQ (CN-IEQ) time-stepping  scheme
\begin{align}
\bra{\partial_{\tau}^{\alpha}\phi}^{n-\frac{1}{2}}
&=-\lambda
\bra{
-\varepsilon^{2}\Delta\phi^{n-\frac{1}{2}}
+\beta\phi^{n-\frac{1}{2}}
+u^{n-\frac{1}{2}}\hat{\phi}^{n-\frac{1}{2}}
-\eta^{n-\frac{1}{2}}
},\label{Nonlocal-AC-CN-IEQ-1}\\
\eta^{n-\frac{1}{2}}
&=\frac{1}{\abs{\Omega}}
\int_{\Omega}
\bra{
-\varepsilon^{2}\Delta\phi^{n-\frac{1}{2}}
+\beta\phi^{n-\frac{1}{2}}
+u^{n-\frac{1}{2}}\hat{\phi}^{n-\frac{1}{2}}
}\zd{\mathbf{x}},\label{Nonlocal-AC-CN-IEQ-2}\\
\partial_{\tau}u^{n-\frac{1}{2}}
&=2\hat{\phi}^{n-\frac{1}{2}}
\partial_{\tau}\phi^{n-\frac{1}{2}},\label{Nonlocal-AC-CN-IEQ-3}
\end{align}
where $\hat{\phi}^{n-\frac{1}{2}}:=\phi^{n-1}+\triangledown_{\tau}\phi^{n-1}/(2\rho_{n-1})$ is the local extrapolation.

%%%%%%%%%%%%%%%%%%%%%%%%%%%%%%%%%%%%%%%%%%%%%%%%%%%%%%%%%%%%%%%%%%%%%%%%%%%%%%%%%%%%
\begin{theorem}\label{Nonlocal-AC-CN-IEQ-Volume}
The CN-IEQ scheme \eqref{Nonlocal-AC-CN-IEQ-1}-\eqref{Nonlocal-AC-CN-IEQ-3} conserves the volume,
\begin{align}\label{Volime-Conservation}
\int_{\Omega}\phi^{n}\zd{\mathbf{x}}
=\int_{\Omega}\phi^{n-1}\zd{\mathbf{x}},
\quad \text{for}\quad 1\leq{n}\leq{N}.
\end{align}
\end{theorem}
%%%%%%%%%%%%%%%%%%%%%%%%%%%%%%%%%%%%%%%%%%%%%%%%%%%%%%%%%%%%%%%%%%%%%%%%%%%%%%%%%%%%
%%%%%%%%%%%%%%%%%%%%%%%%%%%%%%%%%%%%%%%%%%%%%%%%%%%%%%%%%%%%%%%%%%%%%%%%%%%%%%%%%%%%
\begin{proof}
We prove the discrete volume-conserving by induction.
It is easy to check that the volume conservation holds when $n=1$.
In what follows, we assume that the relation \eqref{Volime-Conservation}
is valid for the numerical scheme
\eqref{Nonlocal-AC-CN-IEQ-1}-\eqref{Nonlocal-AC-CN-IEQ-3} with no more than $(N-1)$ indices,
where $N\geq{2}$.
It is sufficient to verify the desired assertion still holds for $n=N$.
Actually, we have
\begin{align}
\bra{a_{0}^{(n)}\triangledown_{\tau}\phi^{n},1}
&=\bra{
\bra{\partial_{\tau}^{\alpha}\phi}^{n-\frac{1}{2}},1}\nonumber\\
&=-\lambda
\bra{
-\varepsilon^{2}\Delta\phi^{n-\frac{1}{2}}
+\beta\phi^{n-\frac{1}{2}}
+u^{n-\frac{1}{2}}\hat{\phi}^{n-\frac{1}{2}}
-\eta^{n-\frac{1}{2}},1}=0,
\end{align}
where the induction assumption $\bra{\phi^{n},1}=\bra{\phi^{n-1},1},\,1\leq{n}\leq{N-1}$
has been used in the derivation of the above identity.
We then have $\bra{\phi^{N},1}=\bra{\phi^{N-1},1}$ that shows the desired result still holds for $n=N$.
Consequently, the relationship \eqref{Volime-Conservation} is valid by the induction.
\end{proof}
%%%%%%%%%%%%%%%%%%%%%%%%%%%%%%%%%%%%%%%%%%%%%%%%%%%%%%%%%%%%%%%%%%%%%%%%%%%%%%%%%%%%%%%%%%

Note that,
the remarkable property of L1$^{+}$ formula in Lemma \ref{Quadratic-Form-New-L1}
implies that the above CN-IEQ scheme \eqref{Nonlocal-AC-CN-IEQ-1}-\eqref{Nonlocal-AC-CN-IEQ-3} is naturally suitable for a general class of nonuniform time meshes.
Precisely, the following result shows that it is unconditionally energy stable.
%%%%%%%%%%%%%%%%%%%%%%%%%%%%%%%%%%%%%%%%%%%%%%%%%%%%%%%%%%%%%%%%%%%%%%%%%%%%%%%%%%%%
\begin{theorem}\label{Nonlocal-AC-CN-IEQ-Decay-Law}
The CN-IEQ scheme \eqref{Nonlocal-AC-CN-IEQ-1}-\eqref{Nonlocal-AC-CN-IEQ-3} preserves the energy dissipation law,
\begin{align}
E\kbra{\phi^{n},u^{n}}
-E\kbra{\phi^{0},u^{0}}
\leq{0},\quad \text{for}\quad 1\leq{n}\leq{N},
\end{align}
such that it is unconditionally stable, where discrete energy is given by
\begin{align*}
E\kbra{\phi^{n},u^{n}}
=\int_{\Omega}
\bra{\frac{\varepsilon^{2}}{2}\abs{\nabla\phi^{n}}^{2}
+\frac{\beta}{2}\abs{\phi^{n}}^{2}
+\frac{1}{4}(u^{n})^{2}}\zd{\mathbf{x}}
-\brab{\frac{\beta}{2}+\frac{\beta^{2}}{4}}\abs{\Omega}.
\end{align*}
\end{theorem}
%%%%%%%%%%%%%%%%%%%%%%%%%%%%%%%%%%%%%%%%%%%%%%%%%%%%%%%%%%%%%%%%%%%%%%%%%%%%%%%%%%%%
\begin{proof}
Taking the inner product of \eqref{Nonlocal-AC-CN-IEQ-1} and \eqref{Nonlocal-AC-CN-IEQ-3} with $\triangledown_{\tau}\phi^{n}$ and $2\tau_{n}u^{n-\frac{1}{2}}$,
respectively,
and adding the resulting two equalities,
we have the following equation
\begin{align}
-\frac{1}{\lambda}
\bra{\bra{\partial_{\tau}^{\alpha}\phi}^{n-\frac{1}{2}},\triangledown_{\tau}\phi^{n}}
&=
\bra{
-\varepsilon^{2}\Delta\phi^{n-\frac{1}{2}}
+\beta\phi^{n-\frac{1}{2}},\triangledown_{\tau}\phi^{n}}
+\frac{1}{4}\bra{\bra{u^{n}}^{2}-\bra{u^{n-1}}^{2},1},
\end{align}
in which the volume conservation \eqref{Volime-Conservation} has been used to show the fact
$(\eta^{n-\frac{1}{2}},\triangledown_{\tau}\phi^{n})=0$.
As a result, we get the following identity
\begin{align}
E\kbrab{\phi^{k},u^{k}}
-E\kbrab{\phi^{k-1},u^{k-1}}
=-\frac{1}{\lambda}
\bra{(\partial_{\tau}^{\alpha}\phi)^{k-\frac{1}{2}},\triangledown_{\tau}\phi^{k}}\quad \text{for $1\leq{k}\leq{n}$.}
\end{align}
By summing the superscript $k$ from $1$ to $n$, we obtain the following inequality
\begin{align*}
E\kbrab{\phi^{n},u^{n}}
-E\kbrab{\phi^{0},u^{0}}
=-\frac{1}{\lambda}\int_{\Omega}
\sum_{k=1}^n\triangledown_{\tau}\phi^{k}
\sum_{j=1}^ka_{k-j}^{(k)}\triangledown_{\tau}\phi^{k}\zd{\mathbf{x}}\leq{0}\quad \text{for $1\leq{n}\leq{N}.$}
\end{align*}
where Lemma \ref{Quadratic-Form-New-L1}
has been used in the last inequality. It completes the proof.
\end{proof}
%%%%%%%%%%%%%%%%%%%%%%%%%%%%%%%%%%%%%%%%%%%%%%%%%%%%%%%%%%%%%%%%%%%%%%%%%%%%%%%%%%%%%%%%%%
\subsection{Numerical approach using SAV}
For the time-fractional Allen-Cahn model with volume constraint  {\eqref{Problem-2}},
we here introduce a scalar auxiliary function $v(t)$ in term of original variable $\phi$ as follows
\begin{align}
v(t)
=\sqrt{\int_{\Omega}
\frac{1}{4}\brab{\phi^{2}-1-\beta}^{2}\zd{\mathbf{x}}+C_{0}},
\end{align}
where the positive constant $C_{0}$ is chosen to ensure the radicand positive and $\beta$ is the regularized parameter.
Therefore, the free energy of the primitive problem could be rewritten into
\begin{align}\label{SAV-Equivalent-Energy}
E\kbra{\phi,v}
=\int_{\Omega}\bra{
\frac{\varepsilon^{2}}{2}\abs{\nabla\phi}^{2}
+\frac{\beta}{2}\abs{\phi}^{2}}
\zd{\mathbf{x}}
+v^{2}
-C_{0}
-\brab{\frac{\beta}{2}+\frac{\beta^{2}}{4}}\abs{\Omega}.
\end{align}
We then could reformulate the problem \eqref{Problem-2} as an equivalent form
\begin{align}
\partial_{t}^{\alpha}\phi
&=-\lambda
\bra{
-\varepsilon^{2}\Delta\phi
+\beta\phi
+V(\phi)v-\eta},\label{Nonlocal-AC-SAV-1}\\
\eta
&=\frac{1}{\abs{\Omega}}\int_{\Omega}
\bra{
-\varepsilon^{2}\Delta\phi
+\beta\phi
+V(\phi)v}\zd{\mathbf{x}},\label{Nonlocal-AC-SAV-2}\\
v_{t}
&=\frac{1}{2}\int_{\Omega}V(\phi)\partial_{t}\phi\zd{\mathbf{x}},\label{Nonlocal-AC-SAV-3}
\end{align}
in which the expression of the notation $V(\phi)$ is given by
\begin{align}
V(\phi)=\frac{(\phi^{2}-1-\beta)\phi}
{\sqrt{\int_{\Omega}
\frac{1}{4}\brab{\phi^{2}-1-\beta}^{2}\zd{\mathbf{x}}+C_{0}}},
\end{align}
with the following initial conditions
\begin{align}
\phi\bra{\mathbf{x},0}
=\phi_{0}\bra{\mathbf{x}}\quad
\text{and}\quad
v(0)
=v\bra{\phi_{0}\bra{\mathbf{x}}}.
\end{align}
It is easy to check that the new system admits the volume-conserving property \eqref{Voluma-Conservation}
and the following energy dissipation law
\begin{align}
E\kbra{\phi(T),v(T)}
-E\kbra{\phi(0),v(0)}
=-\frac{1}{\lambda}
\int_{\Omega}
\mathcal{I}_{t}^{1}
(\partial_{t}\phi\,\mathcal{I}_{t}^{1-\alpha}\partial_{t}\phi)(T)
\zd{\mathbf{x}}
\leq{0}.
\end{align}
%%%%%%%%%%%%%%%%%%%%%%%%%%%%%%%%%%%%%%%%%%%%%%%%%%%%%%%%%%%%%%%%%%%%%%%%%%%%%%%%%%%%%%%%%%

As done in the above subsection, for the equivalent system \eqref{Nonlocal-AC-SAV-1}-\eqref{Nonlocal-AC-SAV-3},
we have the following Crank-Nicolson SAV (CN-SAV) scheme
\begin{align}
\bra{\partial_{\tau}^{\alpha}\phi}^{n-\frac{1}{2}}
&=-\lambda
\bra{
-\varepsilon^{2}\Delta\phi^{n-\frac{1}{2}}
+\beta\phi^{n-\frac{1}{2}}
+V(\hat{\phi}^{n-\frac{1}{2}})v^{n-\frac{1}{2}}
-\eta^{n-\frac{1}{2}}
},\label{Nonlocal-AC-CN-SAV-1}\\
\eta^{n-\frac{1}{2}}
&=\frac{1}{\abs{\Omega}}
\int_{\Omega}
\bra{
-\varepsilon^{2}\Delta\phi^{n-\frac{1}{2}}
+\beta\phi^{n-\frac{1}{2}}
+V(\hat{\phi}^{n-\frac{1}{2}})v^{n-\frac{1}{2}}
}\zd{\mathbf{x}},\label{Nonlocal-AC-CN-SAV-2}\\
\partial_{\tau}v^{n-\frac{1}{2}}
&=\frac{1}{2}\int_{\Omega}V(\hat{\phi}^{n-\frac{1}{2}})
\partial_{\tau}\phi^{n-\frac{1}{2}}\zd{\mathbf{x}},\label{Nonlocal-AC-CN-SAV-3}
\end{align}
Also, we have the following theorems on the volume conservation and energy dissipation by
following the proofs of Theorems \ref{Nonlocal-AC-CN-IEQ-Volume} and \ref{Nonlocal-AC-CN-IEQ-Decay-Law}, respectively.

%%%%%%%%%%%%%%%%%%%%%%%%%%%%%%%%%%%%%%%%%%%%%%%%%%%%%%%%%%%%%%%%%%%%%%%%%%%%%%%%%%%%
\begin{theorem}\label{Nonlocal-AC-CN-SAV-Volume}
The CN-SAV scheme \eqref{Nonlocal-AC-CN-SAV-1}-\eqref{Nonlocal-AC-CN-SAV-3} inherits the volume conservation,
\begin{align}
\int_{\Omega}\phi^{n}\zd{\mathbf{x}}
=\int_{\Omega}\phi^{n-1}\zd{\mathbf{x}},
\quad \text{for}\quad 1\leq{n}\leq{N}.
\end{align}
\end{theorem}
%%%%%%%%%%%%%%%%%%%%%%%%%%%%%%%%%%%%%%%%%%%%%%%%%%%%%%%%%%%%%%%%%%%%%%%%%%%%%%%%%%%%

\begin{theorem}\label{Nonlocal-AC-CN-SAV-Decay-Law}
The CN-IEQ scheme \eqref{Nonlocal-AC-CN-SAV-1}-\eqref{Nonlocal-AC-CN-SAV-3}
preserves the energy dissipation law,
\begin{align}
E\kbra{\phi^{n},v^{n}}
-E\kbra{\phi^{0},v^{0}}
\leq{0},\quad \text{for}\quad 1\leq{n}\leq{N},
\end{align}
such that it is unconditionally stable, in which
\begin{align*}
E\kbra{\phi^{n},v^{n}}
=\int_{\Omega}
\bra{\frac{\varepsilon^{2}}{2}\abs{\nabla\phi^{n}}^{2}
+\frac{\beta}{2}\abs{\phi^{n}}^{2}}\zd{\mathbf{x}}
+\bra{v^{n}}^{2}
-C_{0}
-\brab{\frac{\beta}{2}+\frac{\beta^{2}}{4}}\abs{\Omega}.
\end{align*}
\end{theorem}

%%%%%%%%%%%%%%%%%%%%%%%%%%%%%%%%%%%%%%%%%%%%%%%%%%%%%%%%%%%%%%%%%%%%%%%%%%%%%%%%%%%%
\subsection{Numerical approaches for the time-fractional Cahn-Hilliard model}
In order to make a comparison study between the two volume-preserving models \eqref{Problem-2} and
\eqref{Problem-3},
we include the  CN-IEQ scheme for time-fractional Cahn-Hilliard equation \eqref{Problem-3}
\begin{align*}
(\partial_{\tau}^{\alpha}\phi)^{n-\frac{1}{2}}
&=\lambda\Delta\brab{
-\varepsilon^{2}\Delta\phi^{n-\frac{1}{2}}
+\beta\phi^{n-\frac{1}{2}}
+q^{n-\frac{1}{2}}\hat{\phi}^{n-\frac{1}{2}}
},\\
\partial_{\tau}q^{n-\frac{1}{2}}
&=2\hat{\phi}^{n-\frac{1}{2}}
\partial_{\tau}\phi^{n-\frac{1}{2}},
\end{align*}
in which $q(\phi)=\phi^{2}-1-\beta$, and the following CN-SAV scheme
\begin{align*}
(\partial_{\tau}^{\alpha}\phi)^{n-\frac{1}{2}}
&=\lambda\Delta\brab{-\varepsilon^{2}\Delta\phi^{n-\frac{1}{2}}
+\beta\phi^{n-\frac{1}{2}}
+R(\hat{\phi}^{n-\frac{1}{2}})r^{n-\frac{1}{2}}},\\
\partial_{\tau}r^{n-\frac{1}{2}}
&=\frac{1}{2}\int_{\Omega}R(\hat{\phi}^{n-\frac{1}{2}})
\partial_{\tau}\phi^{n-\frac{1}{2}}\zd{\mathbf{x}},
\end{align*}
where
\[
r(t)
=\sqrt{\int_{\Omega}\frac{1}{4}\bra{\phi^{2}-1-\beta}^{2}+C_{0}},\quad
R\bra{\phi}
=\frac{
\brab{\phi^{2}-1-\beta}\phi
}
{\sqrt{\int_{\Omega}
\frac{1}{4}\bra{\phi^{2}-1-\beta}^{2}\zd{\mathbf{x}}+C_{0}}}.
\]
It is not difficult to show that the two computationally efficient approaches
both are volume-conserving and unconditionally energy stable
by following the proofs of Theorems \ref{Nonlocal-AC-CN-IEQ-Volume} or \ref{Nonlocal-AC-CN-IEQ-Decay-Law},
but we here omit the details for brevity.

%%%%%%%%%%%%%%%%%%%%%%%%%%%%%%%%%%%%%%%%%%%%%%%%%%%%%%%%%%%%%%%%%%%%%%%%%%%%%%%%%%%

\section{Adaptive time-stepping and examples}

The CN-IEQ scheme \eqref{Nonlocal-AC-CN-IEQ-1}-\eqref{Nonlocal-AC-CN-IEQ-3}
and CN-SAV scheme \eqref{Nonlocal-AC-CN-SAV-1}-\eqref{Nonlocal-AC-CN-SAV-3}
are run for the conservative time-fractional Allen-Cahn model \eqref{Problem-2} in this section.
Always, we adopt the fast algorithm to speed up the evaluation of the L1$^{+}$ formula
by setting an absolute tolerance error $\epsilon=10^{-12}$ for the underlaying SOE approximation, see \cite{Ji2019Adaptive}.
The spatial domain $\Omega$ is divided uniformly
using an equispaced mesh in each direction and the Fourier pseudo-spectral method is employed.

Also, to compensate the lack of smoothness of the solution,
the time interval $[0,T]$ is always divided into two parts $[0, T_{0}]$ and $[T_{0}, T]$ with total $N$ subintervals.
Take the graded parameter $\gamma\ge1$ and
apply the graded mesh $t_{k}=T_{0}(k/N_0)^{\gamma}$
in $[0,T_{0}]$ to resolve the initial singularity. Some different time-stepping approaches
are examined in the remainder interval $[T_{0},T]$. In the following context, the \emph{Graded Step} strategy
uses the graded mesh in the starting cell $[0,T_0]$ with the uniform mesh in the remainder interval $(T_0,T]$;
while the \emph{Adaptive Step} strategy employs the graded mesh in $[0,T_0]$
and use certain adaptive time-stepping approach described below in the remainder interval $(T_0,T]$.

%\subsection{A fast version of L1$^{+}$ formula}
%It is evident that the approximation \eqref{New-L1-Formula} is prohibitively expensive for long time simulations due to the long-time memory.
%Therefore, to reduce the computational complesity and storage requirements,
%the sum-of-exponentials (SOE) technique is employed to speed up the evaluation of the L1$^{+}$ formula.
%The crucial point is to approximate the kernel function $\omega_{1-\alpha}$ efficiently on the subinterval $[\Delta{t},\,T]$,
%see \cite[Theorem 2.5]{Jiang2017Fast}.
%\begin{lemma}\label{SOE}
%For the given $\alpha\in(0,\,1)$, an absolute tolerance error $\epsilon\ll{1}$, a cut-off time $\Delta{t}>0$ and a finial time $T$, there exists a positive integer $N_{q}$, positive quadrature nodes $\theta^{\ell}$ and corresponding positive weights $\varpi^{\ell}\,(1\leq{\ell}\leq{N_{q}})$ such that
%\begin{align}
%\bigg|
%\omega_{1-\alpha}(t)
%-\sum_{\ell=1}^{N_{q}}\varpi^{\ell}e^{-\theta^{\ell}t}\bigg|\leq\epsilon,
%\quad
%\forall\,{t}\in[\Delta{t},\,T].\nonumber
%\end{align}
%\end{lemma}
%Note that the detailed discussions  of the fast algorithm of L1$^{+}$ formula
%can be found in subsection 4.1 of \cite{Ji2019Adaptive},
%we here omit it to keep the notation as simple as possible.
%%%%%%%%%%%%%%%%%%%%%%%%%%%%%%%%%%%%%%%%%%%%%%%%%%%%%%%%%%%%%%%%%%%%%%%%%%%%%%%%%%%%
\subsection{Adaptive time-stepping strategy}
%%%%%%%%%%%%%%%%%%%%%%%%%%%%%%%%%%%%%%%%%%%%%%%%%%%%%%%%%%%%%%%%%%%%%%%%%%%%%%%%%%%%

To resolve the time evolutions accurately,
small time steps are always necessary to capture the fast dynamics; but the computations would become quite costly
for the coarsening process, see the initial random perturbation problems in Example \ref{Simulating-Coarsening-Dynamics}.
Fortunately, the proposed numerical schemes are proven in Section 3  to be unconditionally energy stable
and allow large time steps to reduce the computation cost for the coarsening process.
Thus some adaptive time-stepping strategy is useful to resolve the widely varying time scales
and significantly reduce the computational cost.
In current computations, we adjust the size of time step using the formula \cite{Qiao2011An},
\begin{align}
\tau_{ada}
=\max\Bigg\{\tau_{\min},
\frac{\tau_{\max}}{\sqrt{1+\kappa\abs{E^{\prime}(t)}^{2}}}\Bigg\}.
\end{align}
Here the parameters $\tau_{\max},\tau_{\min}$ refer to the predetermined maximum and minimum time steps,
as well as $\kappa$ is chosen by the user to adjust the level of adaptivity.
%%%%%%%%%%%%%%%%%%%%%%%%%%%%%%%%%%%%%%%%%%%%%%%%%%%%%%%%%%%%%%%%%%%%%%%%%%%%%%%%%%%%

\subsection{Accuracy verification}
%%%%%%%%%%%%%%%%%%%%%%%%%%%%%%%%%%%%%%%%%%%%%%%%%%%%%%%%%%%%%%%%%%%%%%%%%

The numerical accuracy in time of our CN-IEQ and CN-SAV schemes is examined by taking $T_0=\min\{1/\gamma,T\}$ and
using the \emph{random} mesh in the remainder time interval $(T_{0},T]$,
that is,  $\tau_{N_{0}+k}:=(T-T_{0})\epsilon_{k}/S_1$ for $1\leq k\leq N_1:=N-N_0$,
%$$\tau_{N_{0}+k}=\frac{(T-T_{0})\epsilon_{k}}{\sum_{k=1}^{N_1}\epsilon_{k}}\quad\text{for $1\leq k\leq N_1:=N-N_0$}$$
where $S_1=\sum_{k=1}^{N_1}\epsilon_{k}$ and $\epsilon_{k}\in(0,1)$ are the random numbers. The maximum norm error $e(N):=\max_{1\leq{n}\leq{N}}\|U^{n}-u^{n}\|_{\infty}$ is recorded in each run  and
 the experimental  order of convergence is computed by
$$\text{Order}\approx\frac{\log\bra{e(N)/e(2N)}}{\log\bra{\tau(N)/\tau(2N)}}$$
in which $\tau(N)$ denotes the maximum time-step size for total $N$ subintervals.

%%%%%%%%%%%%%%%%%%%%%%%%%%%%%%%%%%%%%%%%%%%%%%%%%%%%%%%%%%%%%%%%%%%%%%%%%%%%%%%%%
%%%%%%%%%%%%%%%%%%%%%%%%%%%%%%%%%%%%%%%%%%%%%%%%%%%%%%%%%%%%%%%%%%%%%%%%%%%%%%%%%

%%%%%%%%%%%%%%%%%%%%%%%%%%%%%%%%%%%%%%%%%%%%%%%%%%%%%%%%%%%%%%%%%%%%%%%%%%%%%%%%%
\begin{table}[htb!]
\begin{center}
\caption{Numerical accuracy of CN-IEQ scheme \eqref{Nonlocal-AC-CN-IEQ-1}-\eqref{Nonlocal-AC-CN-IEQ-3}  with $\alpha=0.8,\,\sigma=0.4$}\label{Nonlocal-AC-CN-IEQ-Error} \vspace*{0.3pt}
\def\temptablewidth{1.0\textwidth}
{\rule{\temptablewidth}{0.5pt}}
\begin{tabular*}{\temptablewidth}{@{\extracolsep{\fill}}cccccccccc}
\multirow{2}{*}{$N$} &\multirow{2}{*}{$\tau$} &\multicolumn{2}{c}{$\gamma=2$} &\multirow{2}{*}{$\tau$} &\multicolumn{2}{c}{$\gamma=5$} &\multirow{2}{*}{$\tau$}&\multicolumn{2}{c}{$\gamma=6$} \\
             \cline{3-4}          \cline{6-7}         \cline{9-10}
             &     &$e(N)$   &Order   &      &$e(N)$   &Order    &     &$e(N)$    &Order\\
\midrule
  10    &1.72e-01  &2.41e-02 &$-$  &2.69e-01 &2.44e-02 &$-$  &2.69e-01 &3.08e-02  &$-$\\
  20    &1.00e-01  &1.38e-02 &1.02 &1.24e-01 &4.59e-03 &2.17 &1.31e-01 &6.42e-03  &2.18\\
  40    &5.39e-02  &8.50e-03 &0.79 &1.07e-02 &3.46e-05 &2.00 &6.75e-02 &1.45e-03  &2.24\\
  80    &3.07e-02  &5.22e-03 &0.87 &3.49e-02 &3.16e-04 &1.90 &3.33e-02 &3.13e-04  &2.18\\
\end{tabular*}
{\rule{\temptablewidth}{0.5pt}}
\end{center}
\end{table}
%%%%%%%%%%%%%%%%%%%%%%%%%%%%%%%%%%%%%%%%%%%%%%%%%%%%%%%%%%%%%%%%%%%%%%%%%%%%%%%%%%%%%%	 	

%%%%%%%%%%%%%%%%%%%%%%%%%%%%%%%%%%%%%%%%%%%%%%%%%%%%%%%%%%%%%%%%%%%%%%%%%%%%%%%%%%%%%%
\begin{table}[htb!]
\begin{center}
\caption{Numerical accuracy of CN-SAV scheme \eqref{Nonlocal-AC-CN-SAV-1}-\eqref{Nonlocal-AC-CN-SAV-3}  with $\alpha=0.8,\,\sigma=0.4$}\label{Nonlocal-AC-CN-SAV-Error} \vspace*{0.3pt}
\def\temptablewidth{1.0\textwidth}
{\rule{\temptablewidth}{0.5pt}}
\begin{tabular*}{\temptablewidth}{@{\extracolsep{\fill}}cccccccccc}
\multirow{2}{*}{$N$} &\multirow{2}{*}{$\tau$} &\multicolumn{2}{c}{$\gamma=2$} &\multirow{2}{*}{$\tau$} &\multicolumn{2}{c}{$\gamma=5$} &\multirow{2}{*}{$\tau$}&\multicolumn{2}{c}{$\gamma=6$} \\
             \cline{3-4}          \cline{6-7}         \cline{9-10}
             &    &$e(N)$   &Order   &       &$e(N)$    &Order    &      &$e(N)$    &Order\\
\midrule
  10    &1.86e-01 &2.40e-02 &$-$  &3.23e-01	 &1.86e-02  &$-$  &2.88e-01	 &2.33e-02  &$-$\\
  20    &1.02e-01 &1.38e-02	&0.92 &1.34e-01	 &2.99e-03  &2.08 &1.27e-01	 &6.78e-03	&1.51\\
  40    &5.31e-02 &8.51e-03 &0.74 &6.95e-02	 &1.06e-03  &1.57 &6.08e-02	 &9.51e-04	&2.66\\
  80    &2.76e-02 &5.22e-03	&0.75 &3.83e-02	 &3.16e-04  &2.03 &3.36e-02	 &2.55e-04	&2.21\\
\end{tabular*}
{\rule{\temptablewidth}{0.5pt}}
\end{center}
\end{table}		 	
%%%%%%%%%%%%%%%%%%%%%%%%%%%%%%%%%%%%%%%%%%%%%%%%%%%%%%%%%%%%%%%%%%%%%%%%%%%%%%%%%%%%%%

\begin{example}\label{Accuracy-Test-CN-IEQ-SAV}
Consider the model
$\partial_{t}^{\alpha}\phi
=-\frac{\delta{E}}{\delta\phi}+\eta(t)+g(\mathbf{x},t)$ with $\varepsilon^{2}=0.5$
for $\mathbf{x}\in(0,2\pi)^{2}$ and $0<t<1$
such that it has an exact solution $\phi=\omega_{1+\sigma}(t)\sin(x)\sin(y)$.
\end{example}

The spatial domain is discretized by using $128 \times 128$ meshes.
We chose the fractional order $\alpha=0.8$, the regular parameter $\sigma=0.4$ and
the artificial parameters $\beta=1$ and $C_{0}=1$.
Tables \ref{Nonlocal-AC-CN-IEQ-Error}  and  \ref{Nonlocal-AC-CN-SAV-Error} list
the numerical results of CN-IEQ
and CN-SAV approaches with different graded parameters $\gamma$.
It is seen that the time accuracy is of order
$O(\tau^{\gamma\sigma})$ when $\gamma<2/\sigma$,
and the second-order accuracy is achieved when $\gamma\geq\gamma_{\mathrm{opt}}=2/\sigma$.
They suggest that the time accuracy is about of $O(\tau^{\min\{\gamma\sigma,2\}})$ in time
although no theoretical proof is available up to now.

\subsection{Numerical comparisons}

%%%%%%%%%%%%%%%%%%%%%%%%%%%%%%%%%%%%%%%%%%%%%%%%%%%%%%%%%%%%%%%%%%%%%%%%%%%%%%%%%%%%%%
%%%%%%%%%%%%%%%%%%%%%%%%%%%%%%%%%%%%%%%%%%%%%%%%%%%%%%%%%%%%%%%%%%%%%%%%%%%%%%%%%%%%%%
\begin{figure}[htb!]
\centering
\includegraphics[width=3.0in,height=2.0in]{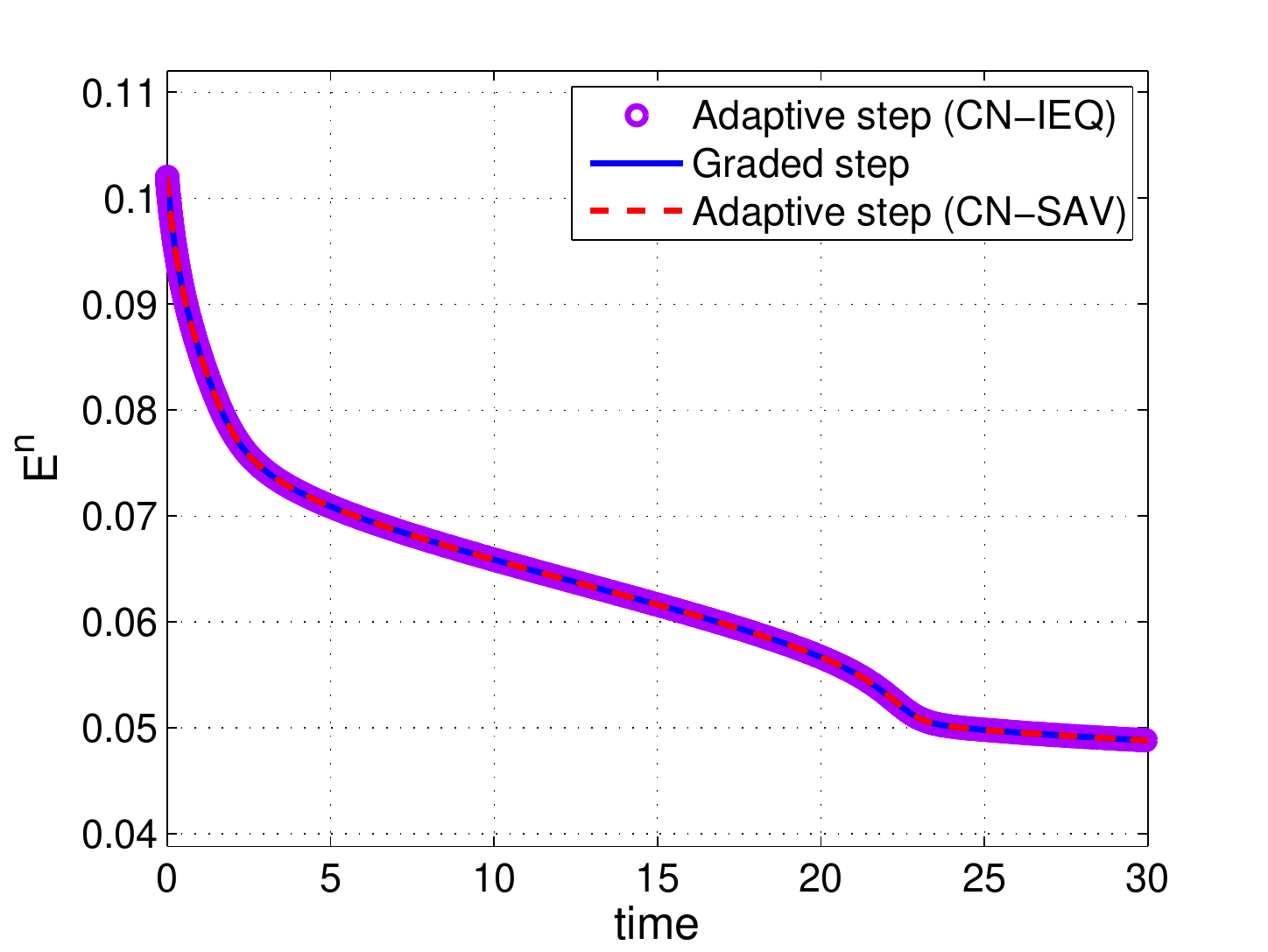}
\includegraphics[width=3.0in,height=2.0in]{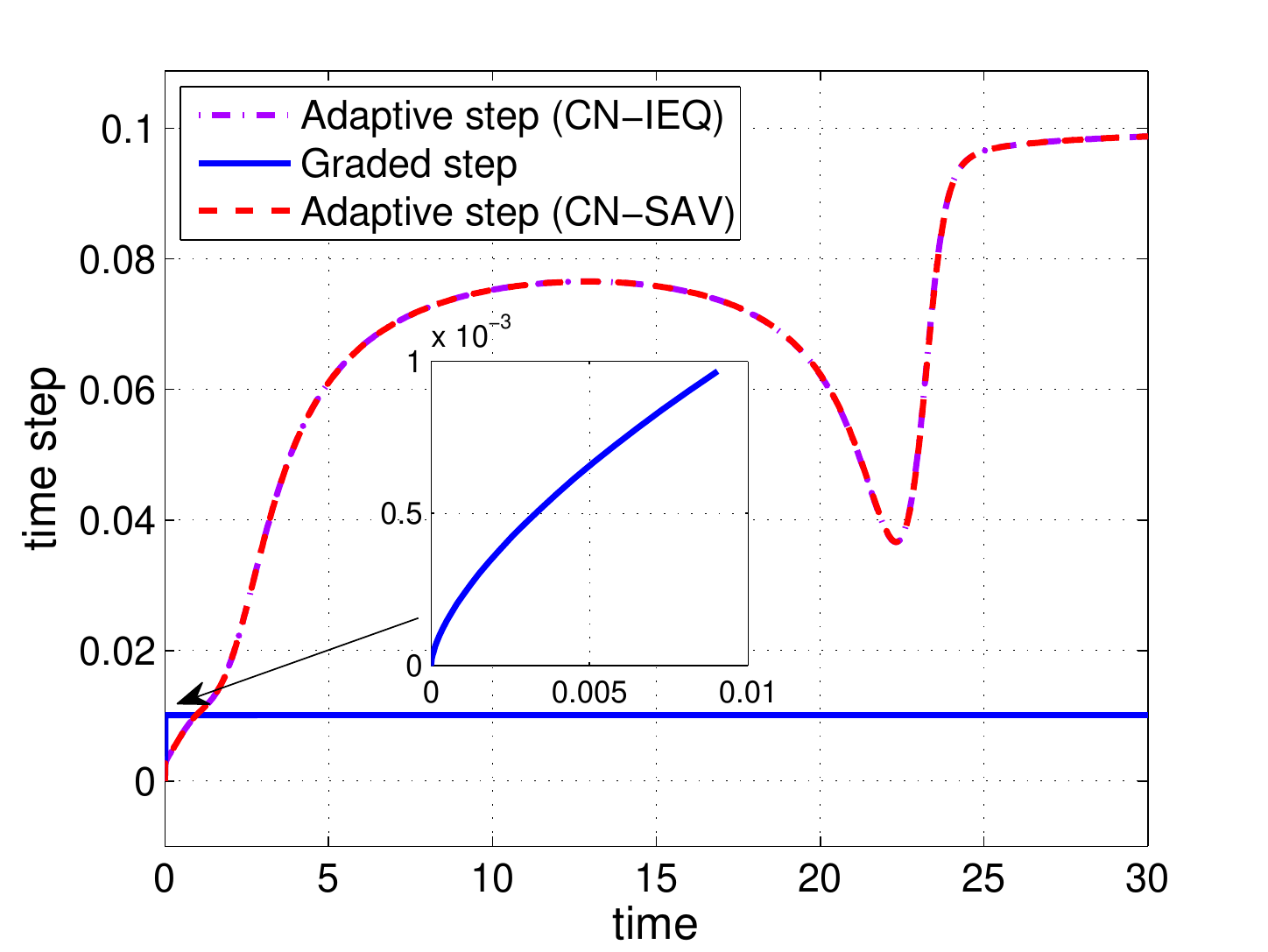}
\caption{Evolutions of energy (left) and time steps (right) of
   the conservative time-fractional Allen-Cahn equation
    using different time strategies until final time $T=30$.}
\label{Comparison-Adaptive-Energy-Curves}
\end{figure}
%%%%%%%%%%%%%%%%%%%%%%%%%%%%%%%%%%%%%%%%%%%%%%%%%%%%%%%%%%%%%%%%%%%%%%%%%%%%%%%%%%%%%%
%%%%%%%%%%%%%%%%%%%%%%%%%%%%%%%%%%%%%%%%%%%%%%%%%%%%%%%%%%%%%%%%%%%%%%%%%%%%%%%%%%%%%%
\begin{figure}[htb!]
\centering
\includegraphics[width=1.47in]{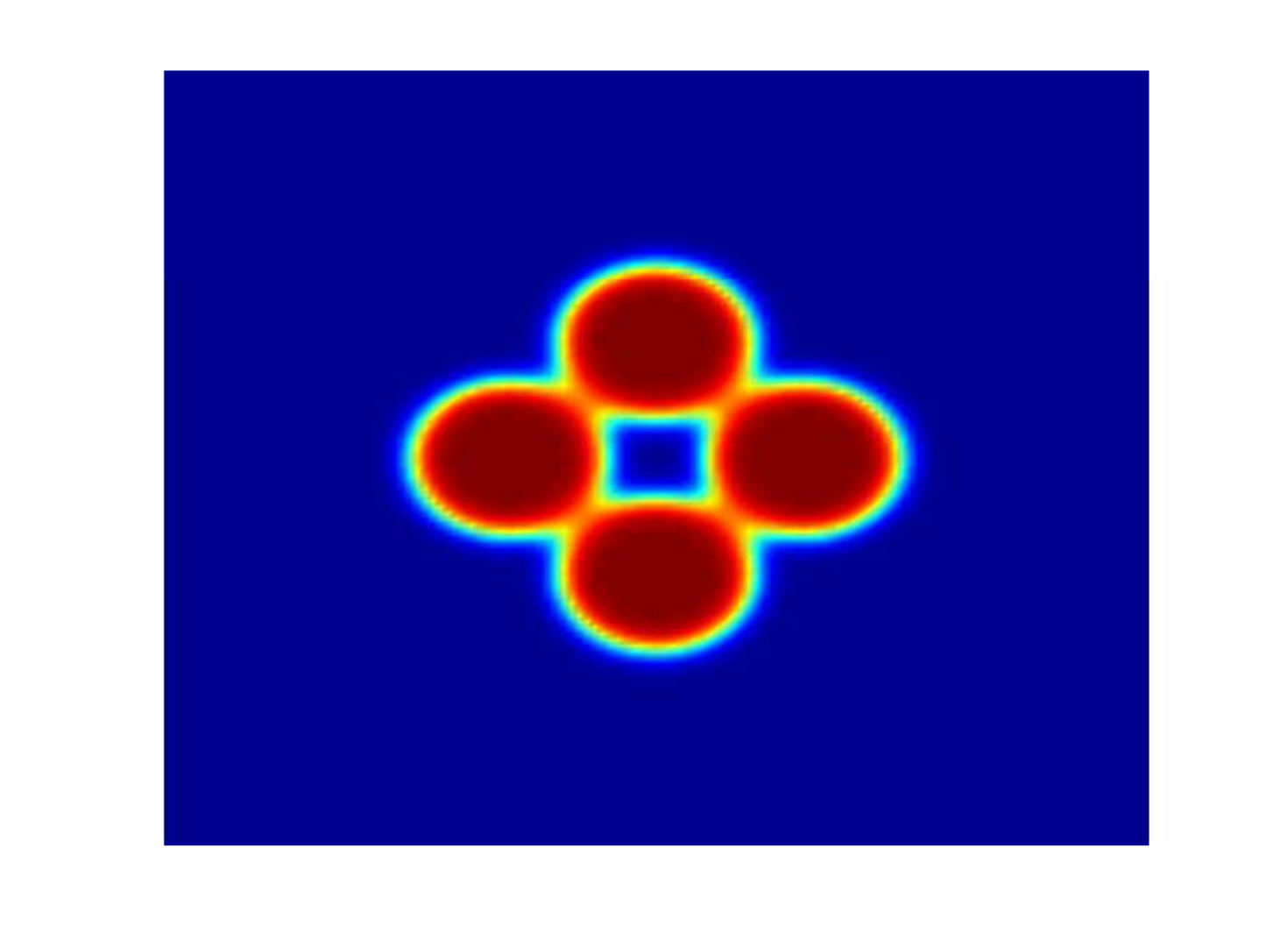}
\includegraphics[width=1.47in]{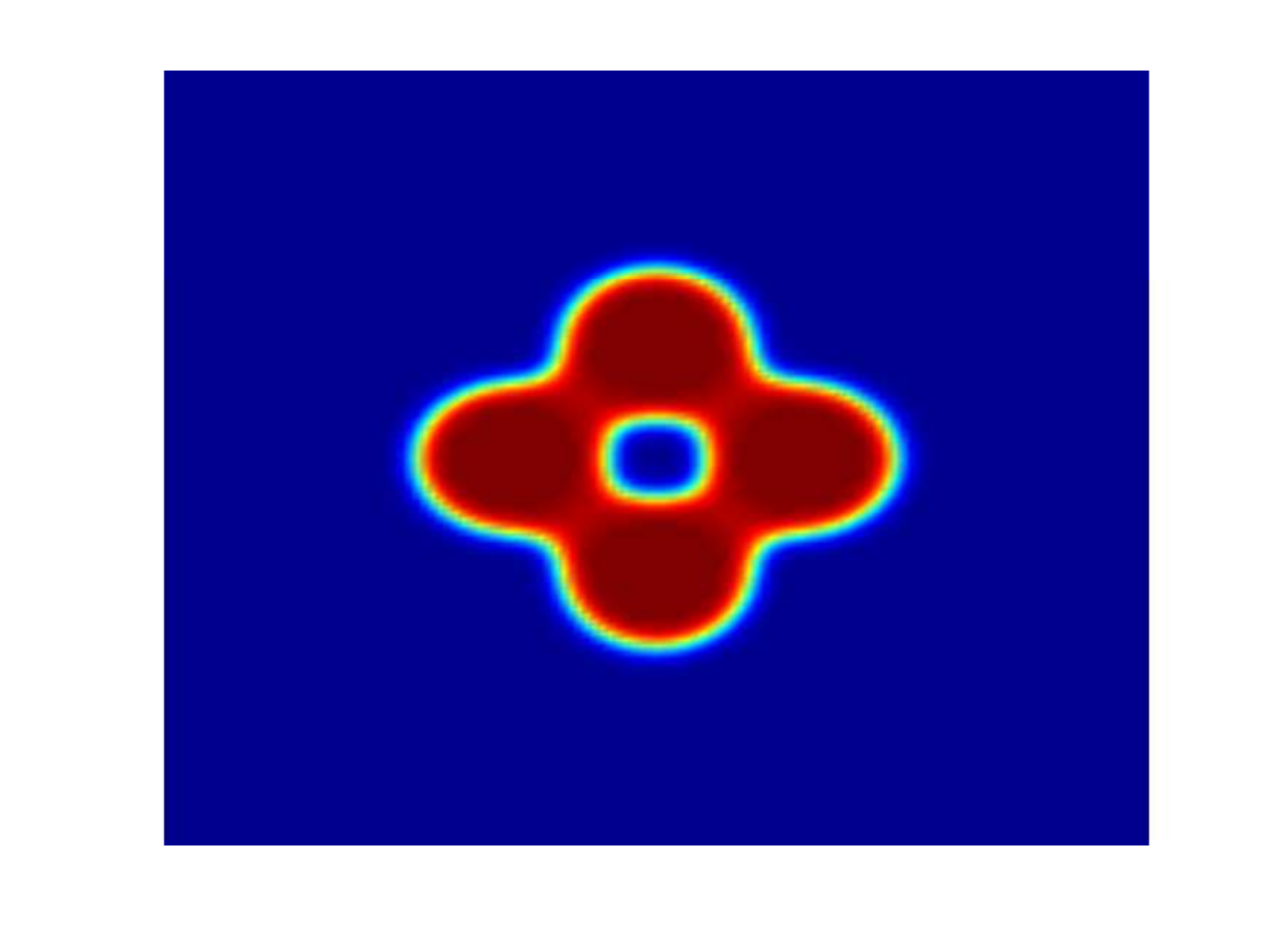}
\includegraphics[width=1.47in]{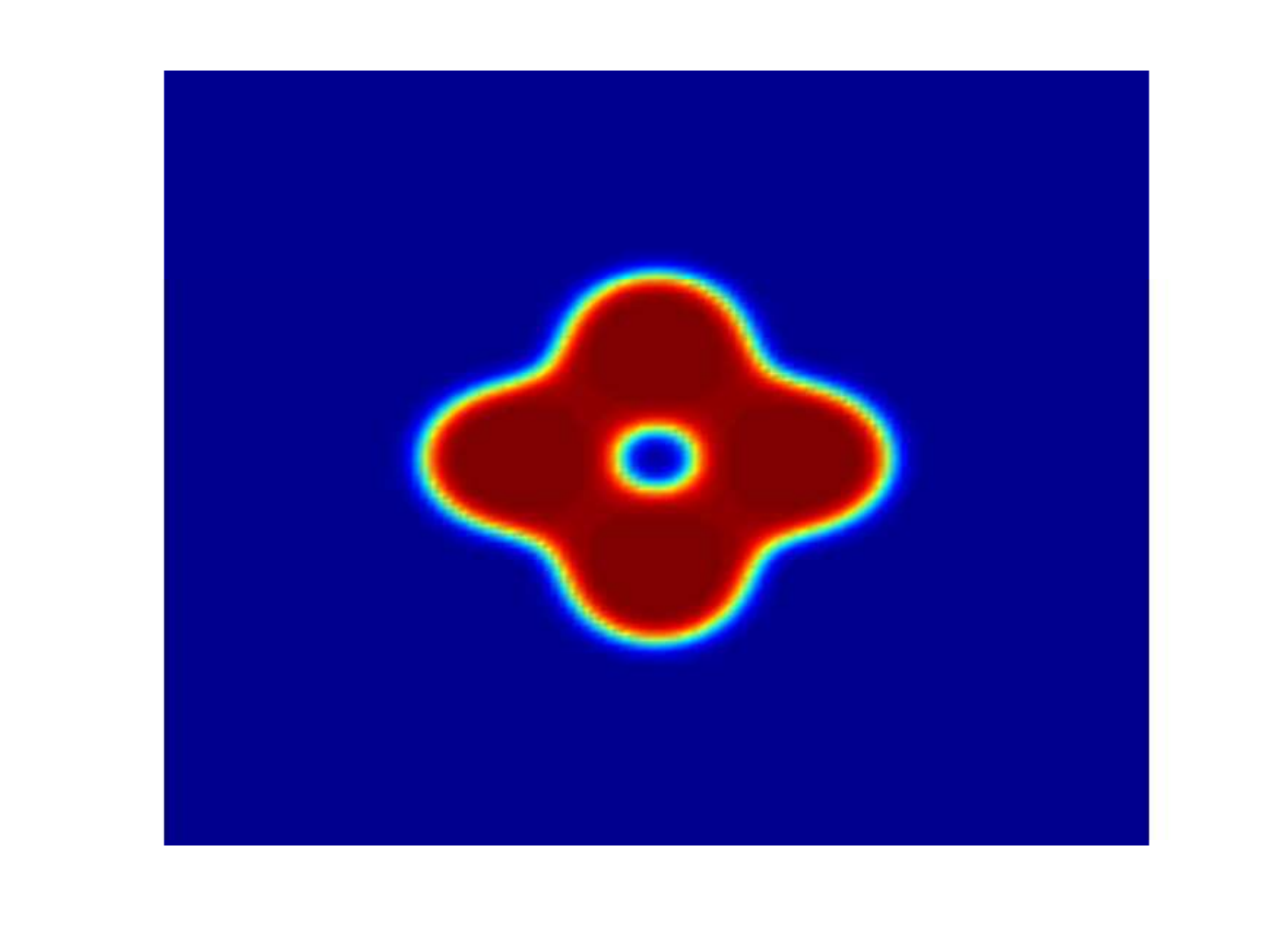}
\includegraphics[width=1.47in]{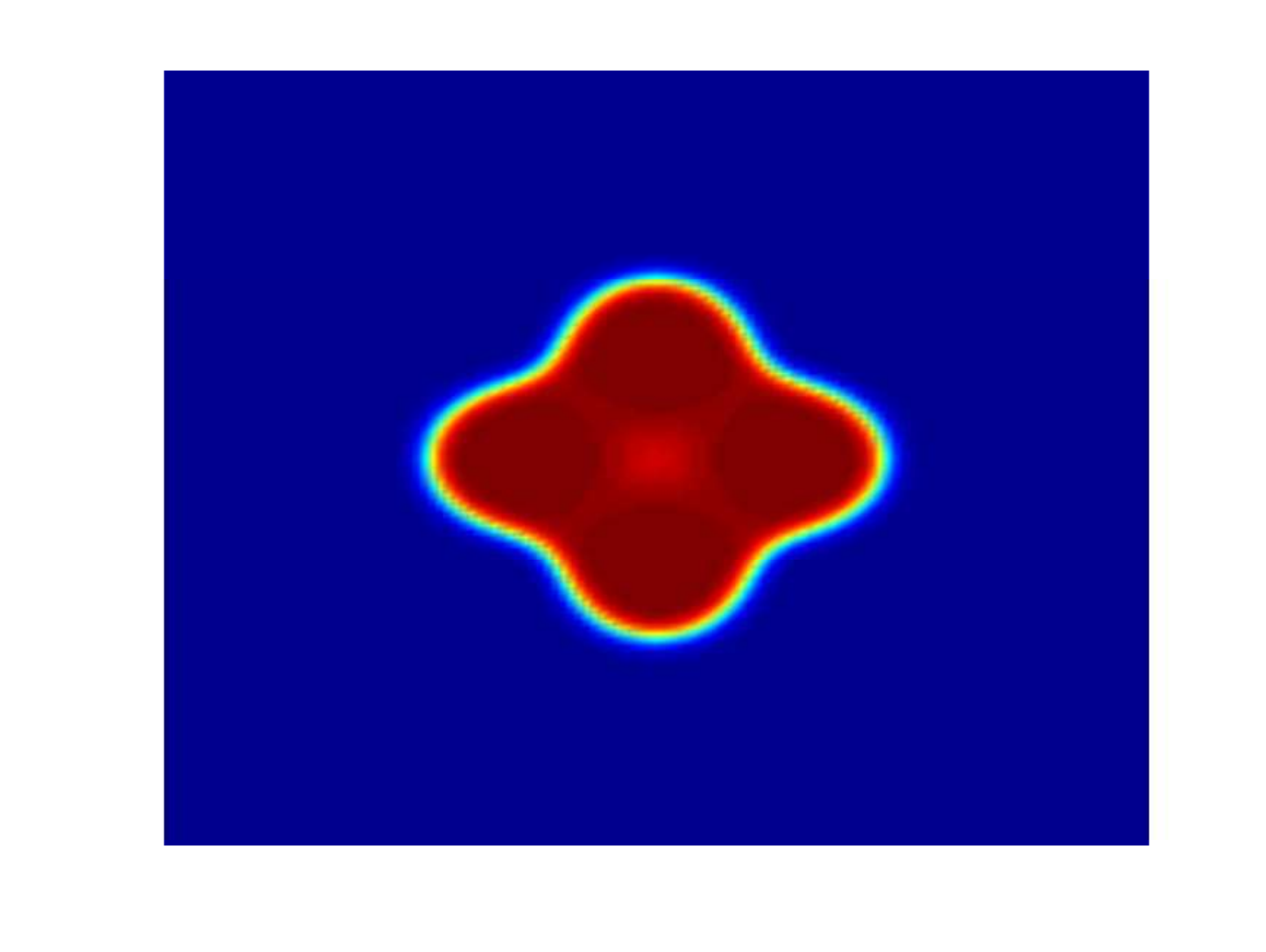}\\
\includegraphics[width=1.47in]{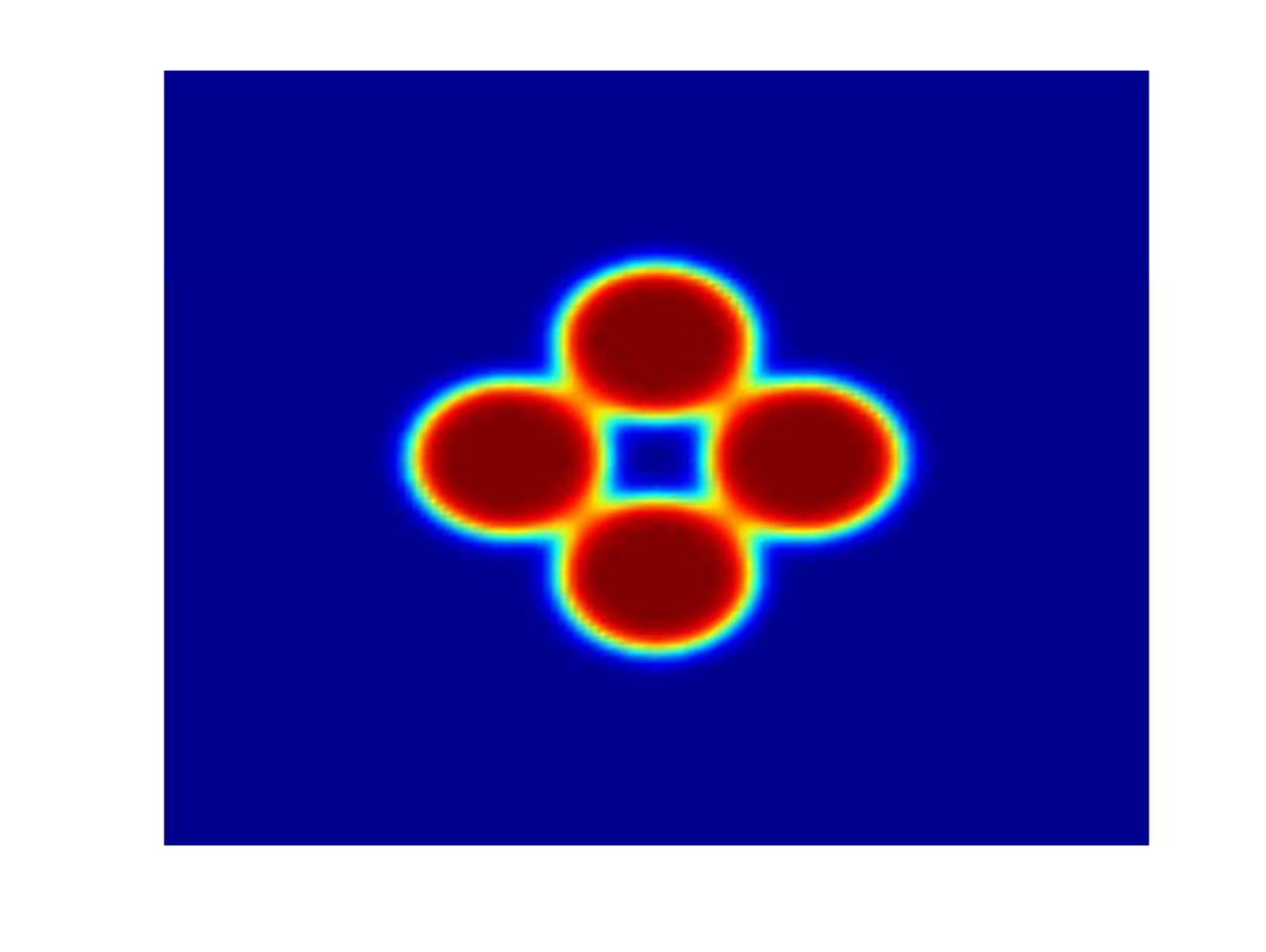}
\includegraphics[width=1.47in]{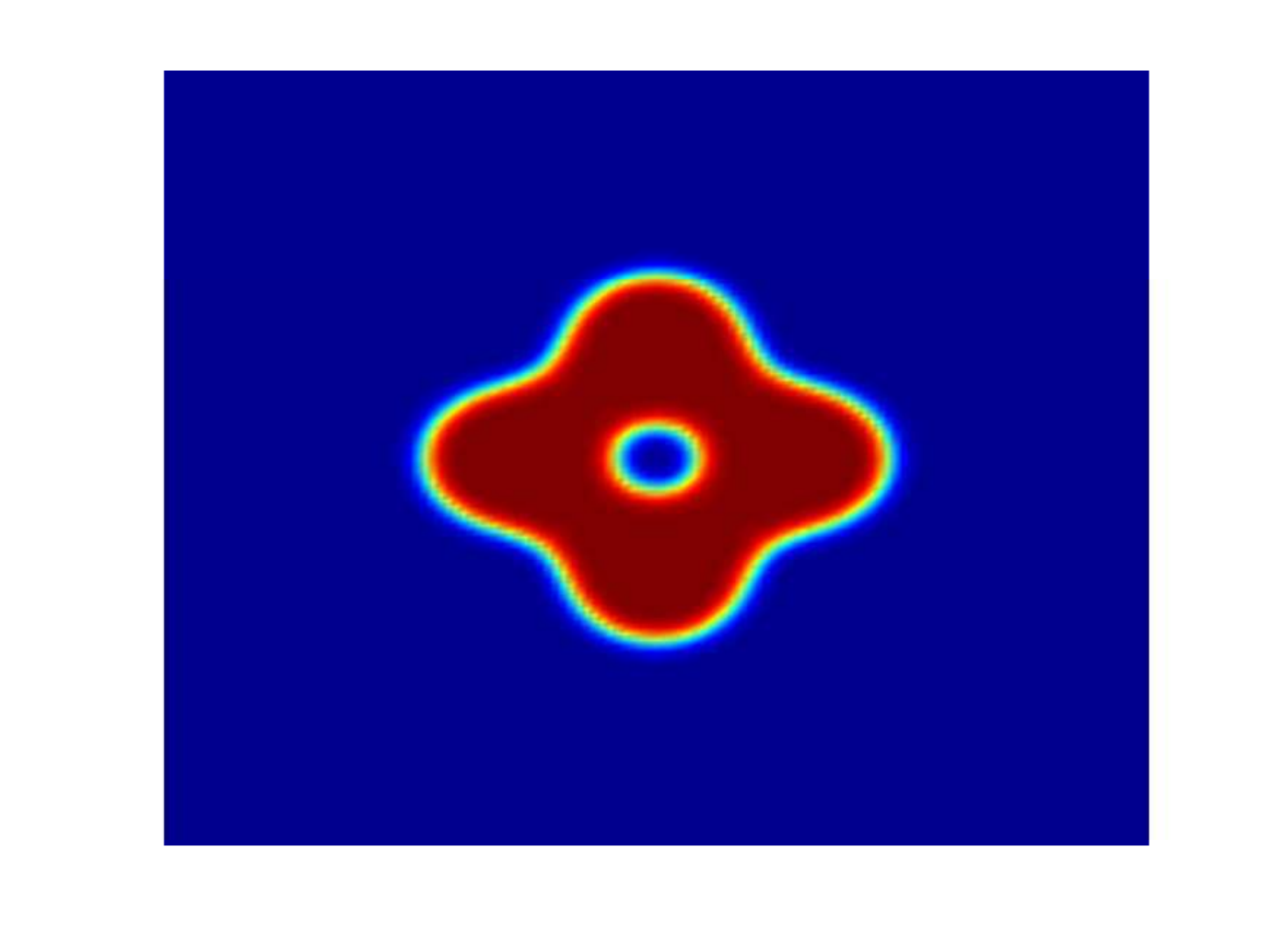}
\includegraphics[width=1.47in]{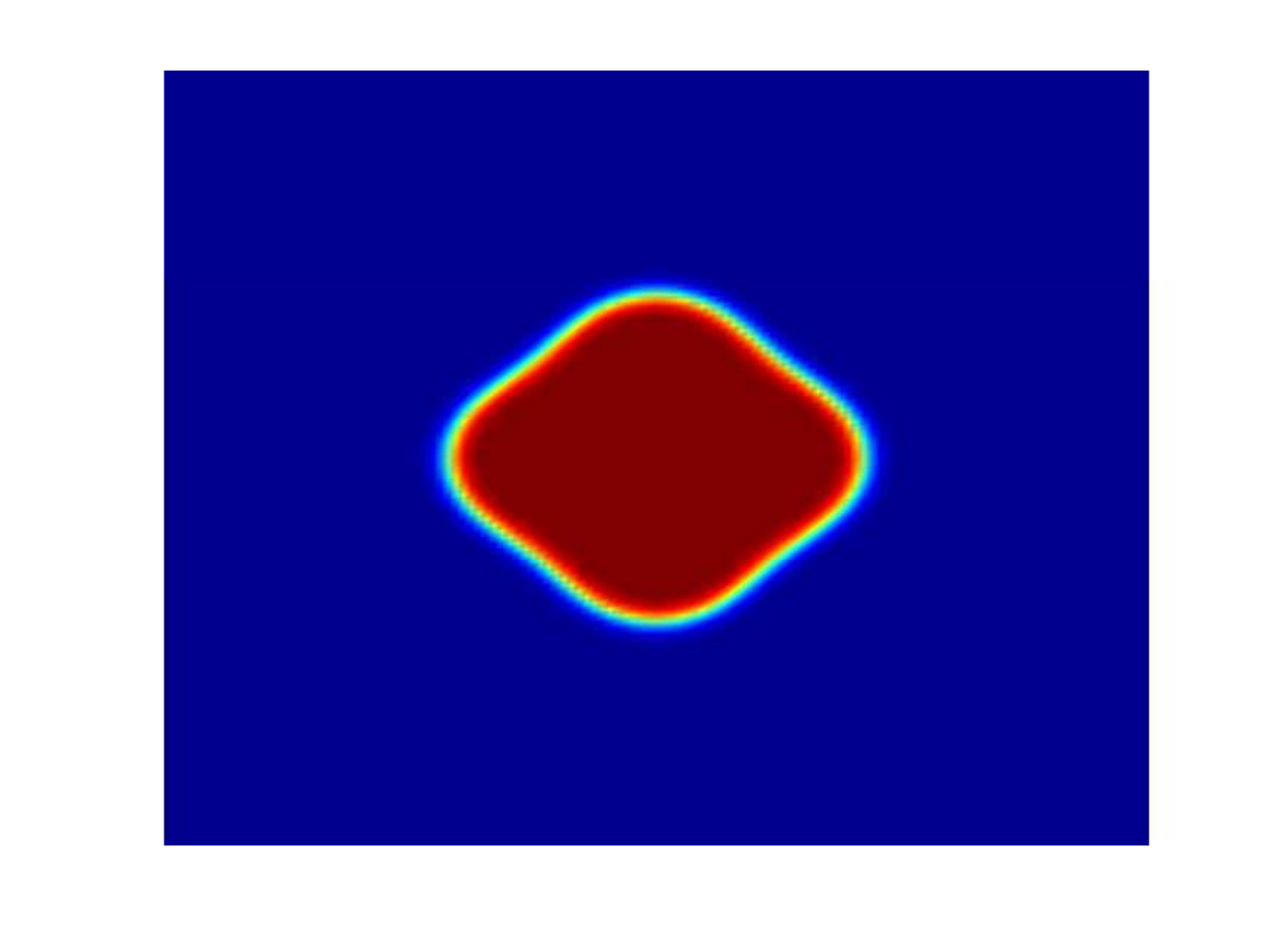}
\includegraphics[width=1.47in]{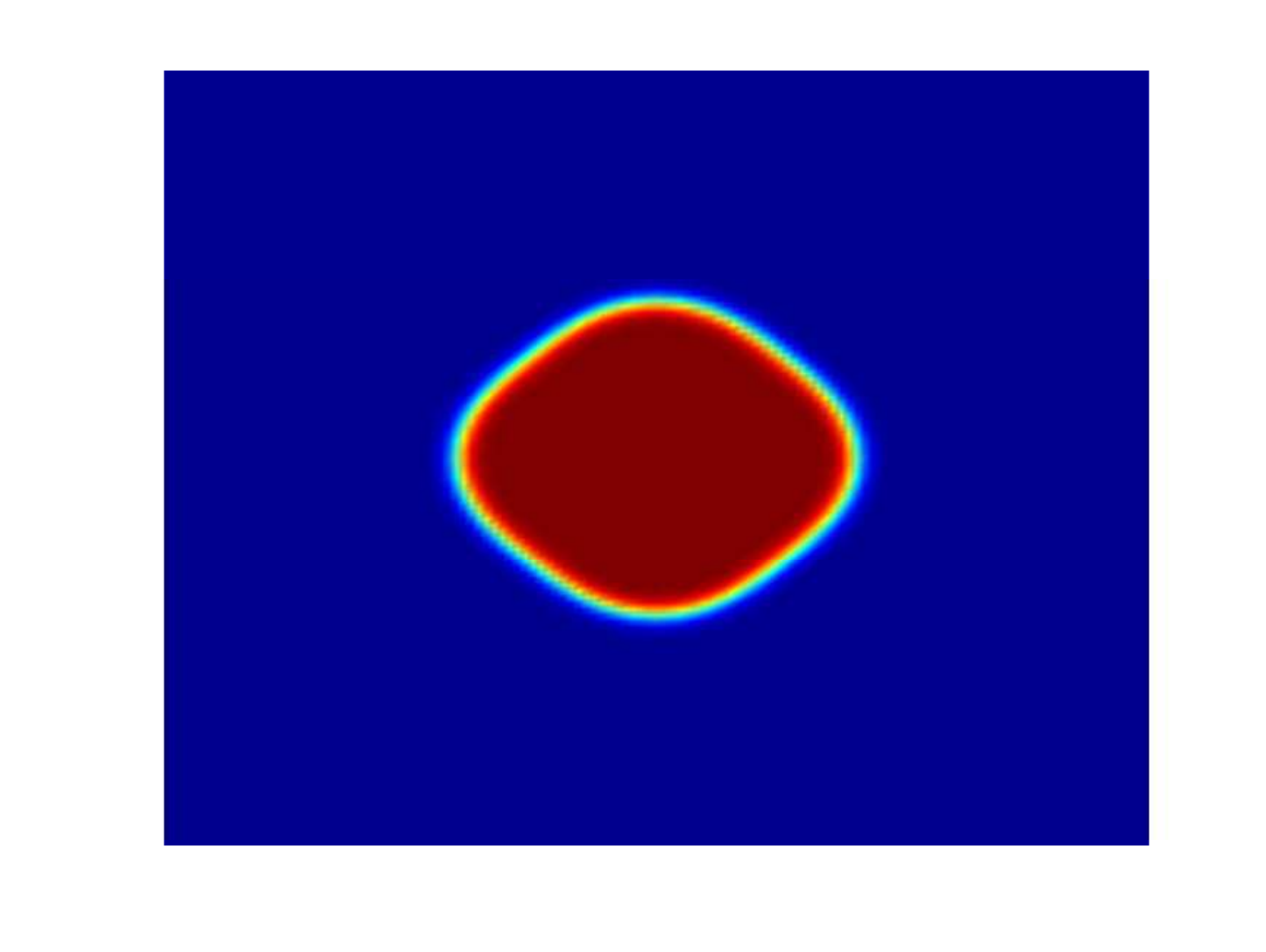}\\
\includegraphics[width=1.47in]{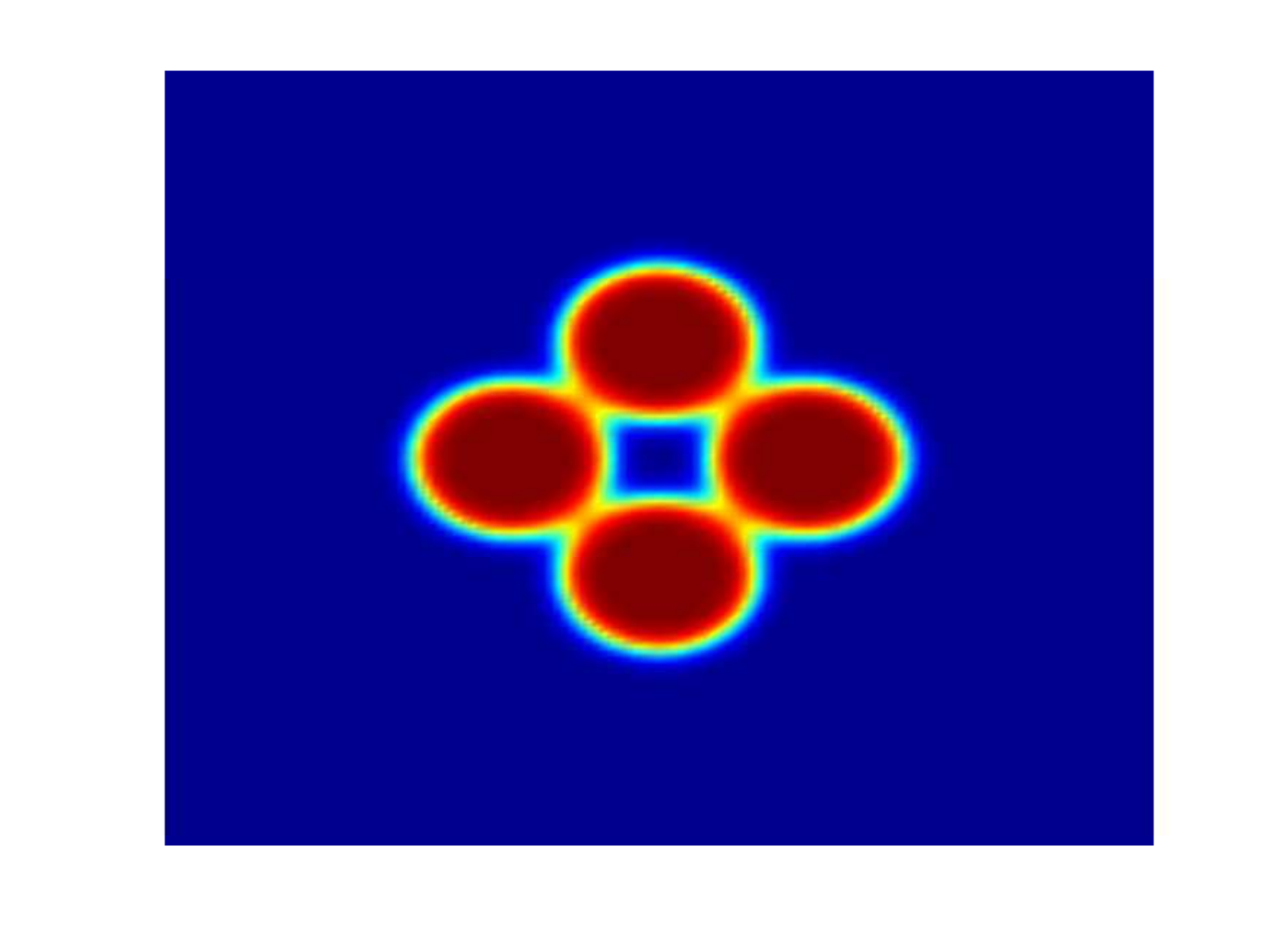}
\includegraphics[width=1.47in]{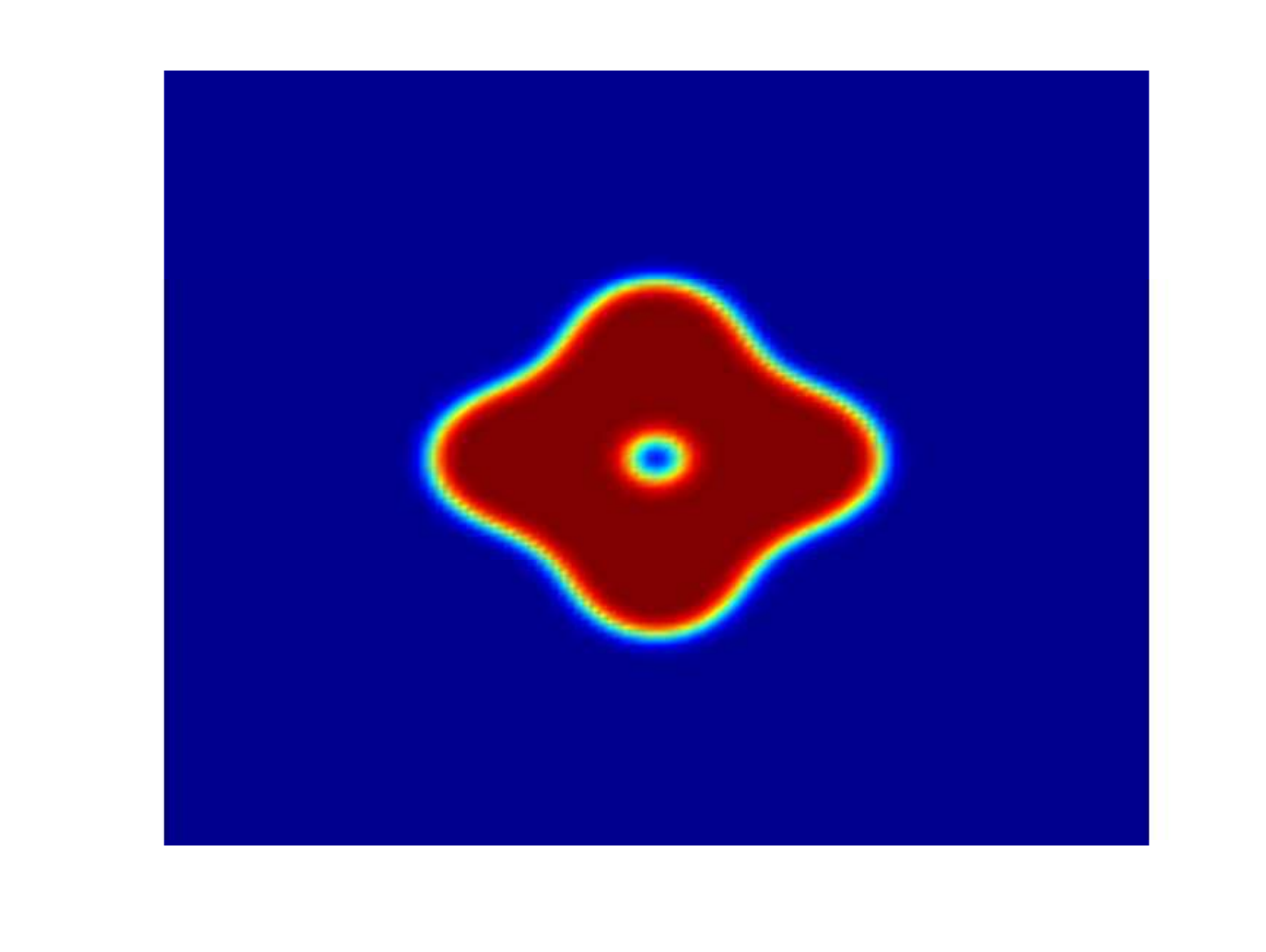}
\includegraphics[width=1.47in]{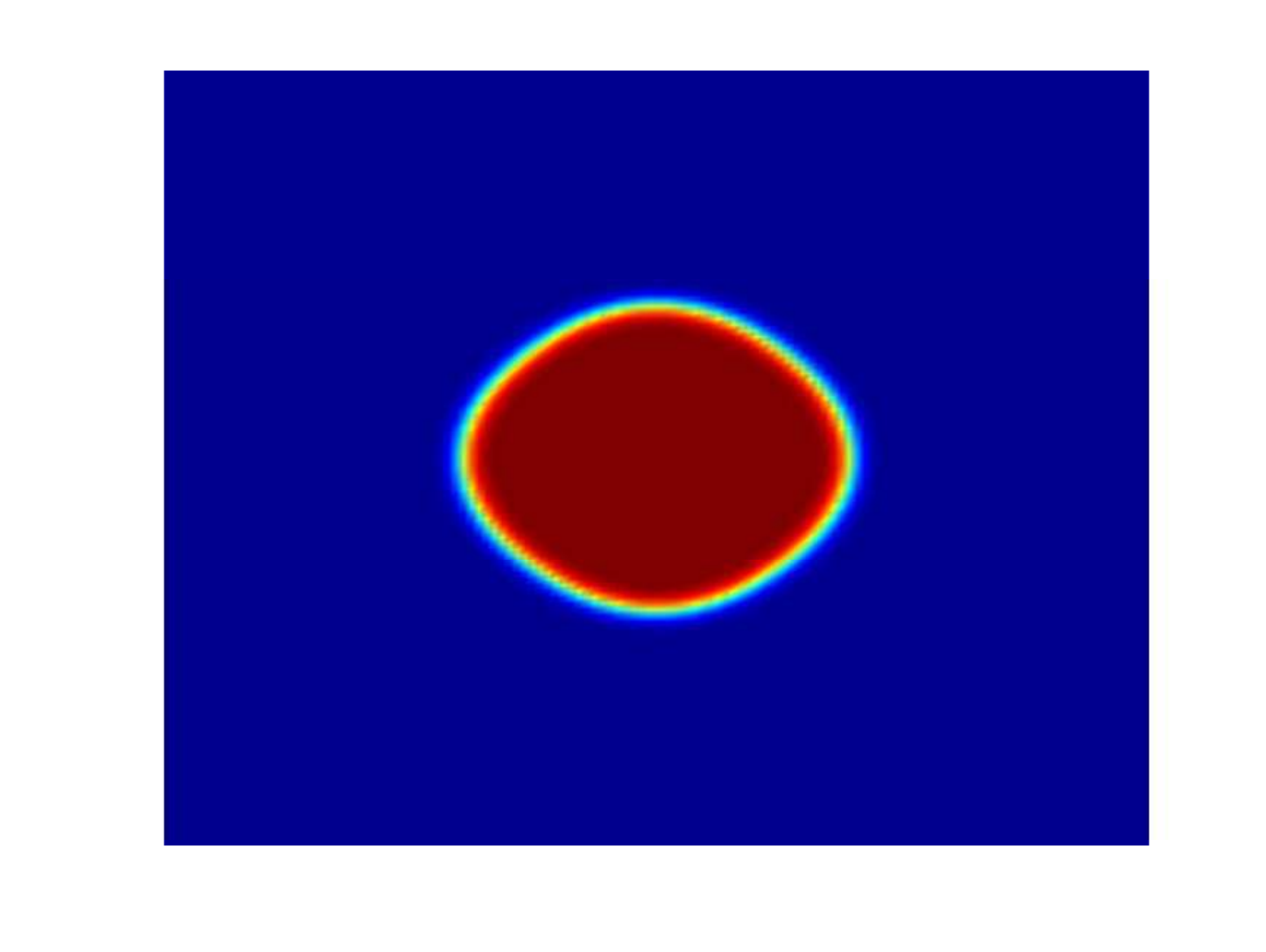}
\includegraphics[width=1.47in]{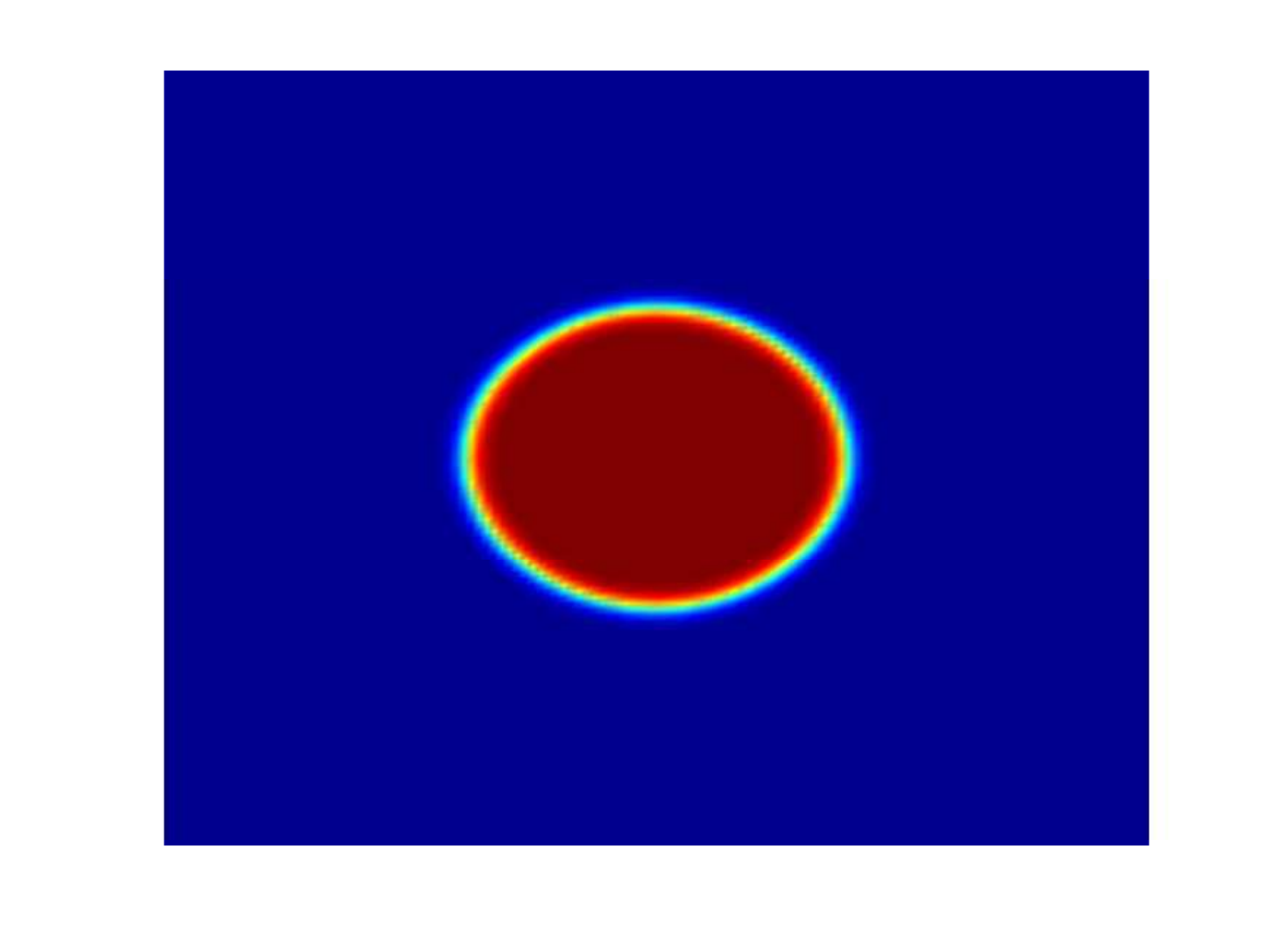}\\
\caption{Solution snapshots at $t=1, 30, 100, 200$ (from left to right)
  for three fractional orders $\alpha=0.4,\,0.7$ and $0.9$ (from top to bottom), respectively.}
\label{Nonlocal-AC-Drops}
\end{figure}
%%%%%%%%%%%%%%%%%%%%%%%%%%%%%%%%%%%%%%%%%%%%%%%%%%%%%%%%%%%%%%%%%%%%%%%%%%%%%%%%%%%%%
%%%%%%%%%%%%%%%%%%%%%%%%%%%%%%%%%%%%%%%%%%%%%%%%%%%%%%%%%%%%%%%%%%%%%%%%%%%%%%%%%%%%%
\begin{figure}[htb!]
\centering
\includegraphics[width=3.0in,height=2.0in]{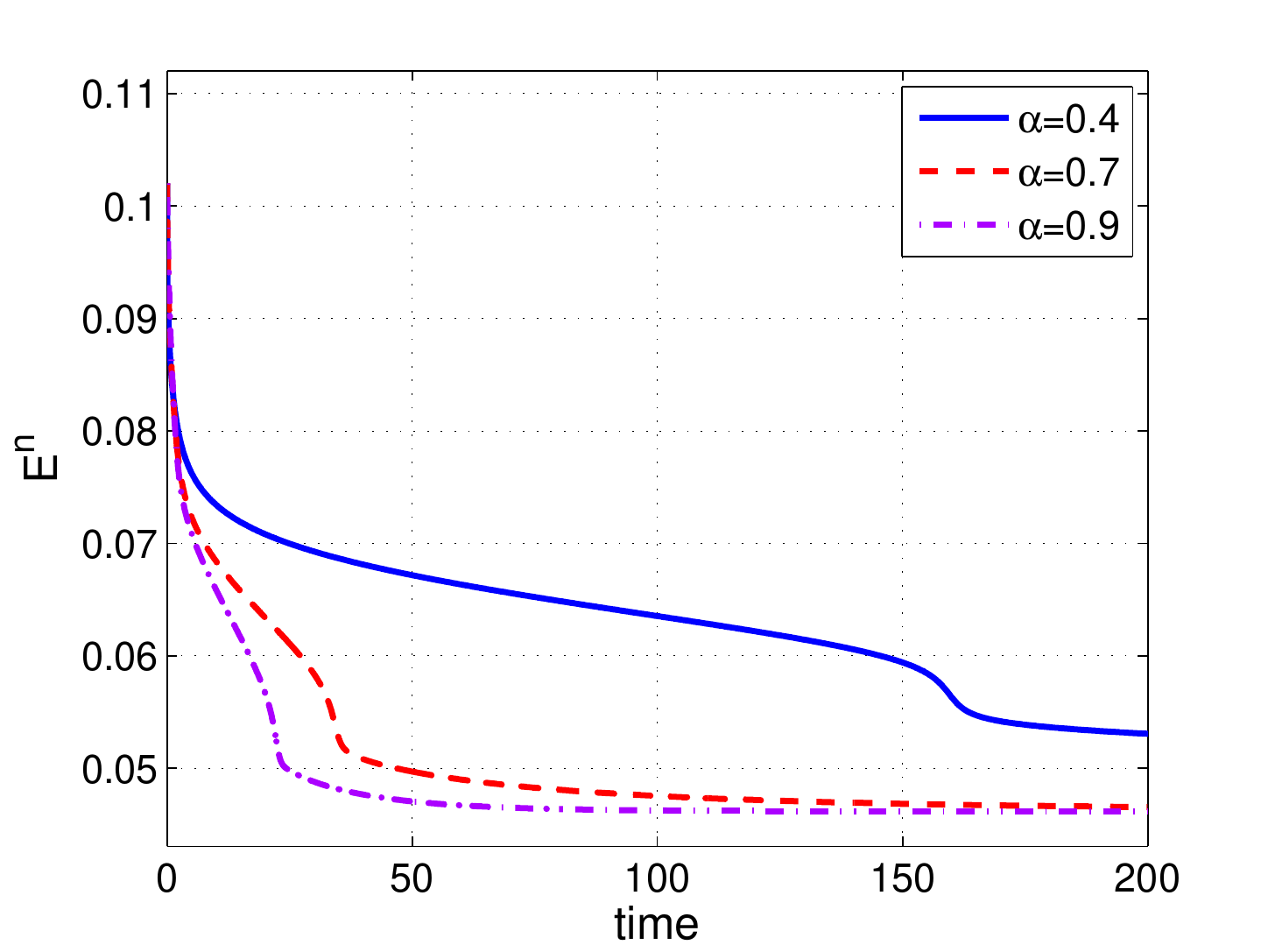}
\includegraphics[width=3.0in,height=2.0in]{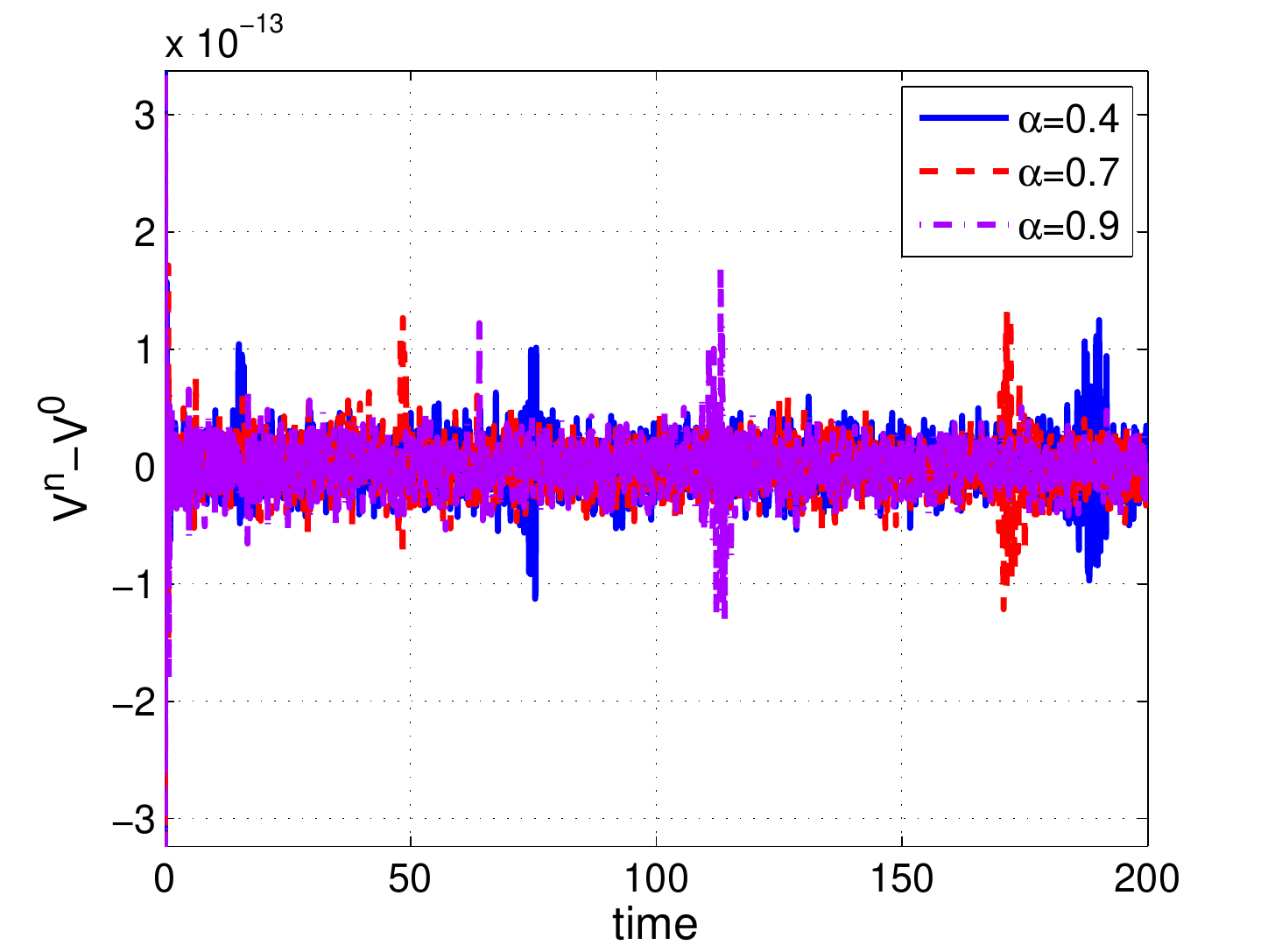}
\caption{Evolutions of energy and volume curves (from left to right) for
  the conservative time-fractional Allen-Cahn equation with fractional order $\alpha=0.4,\,0.7$ and $0.9$, respectively.}
\label{Comparison-Nonlocal-AC-Energy-Mass}
\end{figure}
%%%%%%%%%%%%%%%%%%%%%%%%%%%%%%%%%%%%%%%%%%%%%%%%%%%%%%%%%%%%%%%%%%%%%%%%%%%%%%%%%%%%%

\begin{example}\label{Simulating-Four-Drops}
Consider three different phase field models,
covering time-fractional Allen-Cahn, the conservative version
and time-fractional  Cahn-Hilliard equations, with the coefficients $\lambda=1$ and $\varepsilon=0.02$.
The CN-IEQ and CN-SAV methods with the parameters $\beta=4$ and $C_{0}=1$
are applied to simulate the merging of four drops
with an initial condition
\begin{align}
\phi_{0}\bra{\mathbf{x}}
=&-\tanh\bra{\bra{(x-0.3)^{2}+y^{2}-0.2^2}/\varepsilon}
\tanh\bra{\bra{(x+0.3)^{2}+y^{2}-0.2^2}/\varepsilon}\nonumber\\
&\times\tanh\bra{\bra{x^{2}+(y-0.3)^{2}-0.2^2}/\varepsilon}
\tanh\bra{\bra{x^{2}+(y+0.3)^{2}-0.2^2}/\varepsilon}.
\end{align}
The computational domain $\Omega=(-1,1)^{2}$ is divided uniformly into 128 parts in each direction.
\end{example}
%%%%%%%%%%%%%%%%%%%%%%%%%%%%%%%%%%%%%%%%%%%%%%%%%%%%%%%%%%%%%%%%%%%%%%%%%%%%%%%%%%%%%%

%%%%%%%%%%%%%%%%%%%%%%%%%%%%%%%%%%%%%%%%%%%%%%%%%%%%%%%%%%%%%%%%%%%%%%%%%%%%%%%%%%%%%%%%%%
\begin{figure}[htb!]
\centering
\includegraphics[width=1.47in]{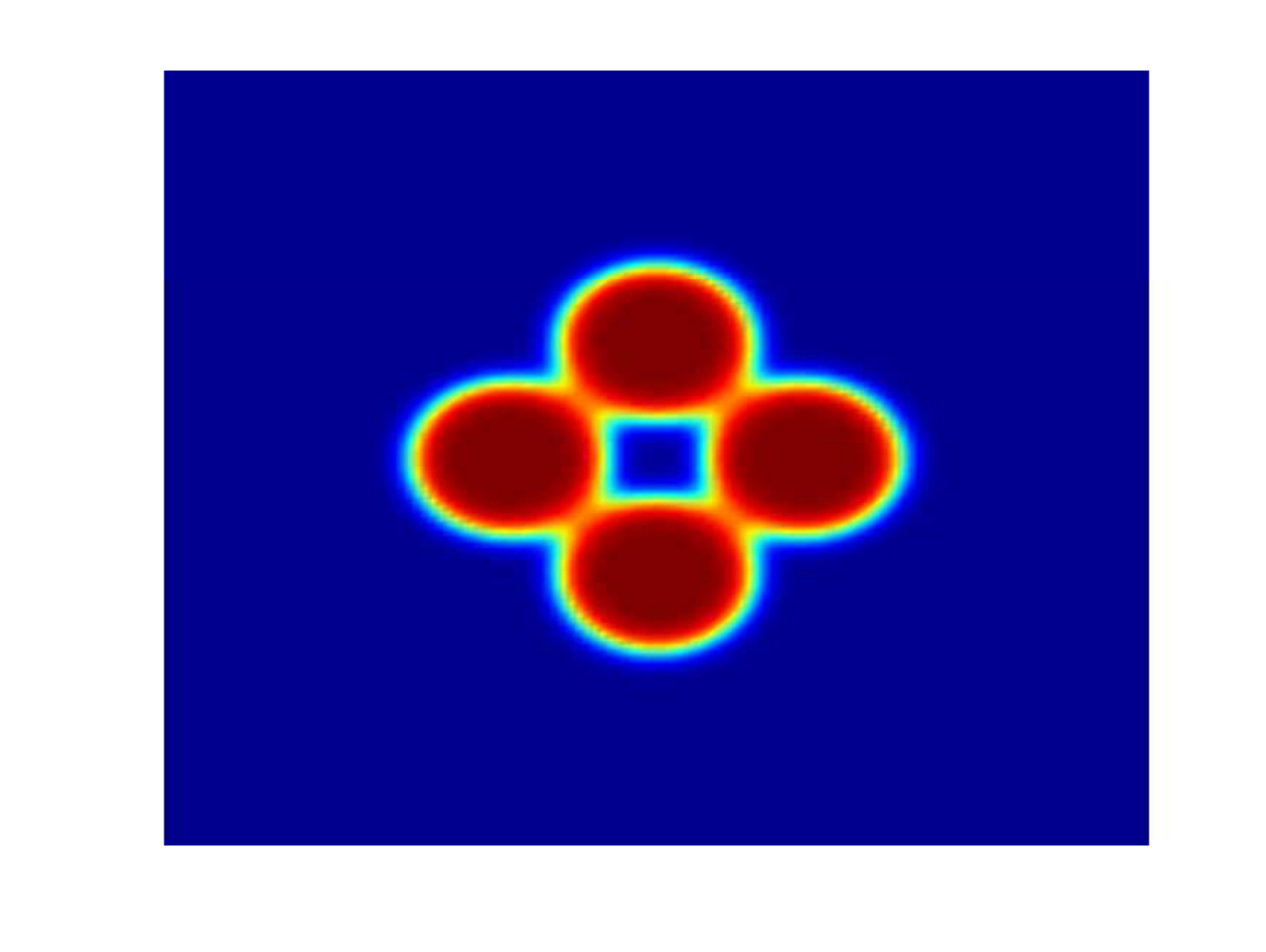}
\includegraphics[width=1.47in]{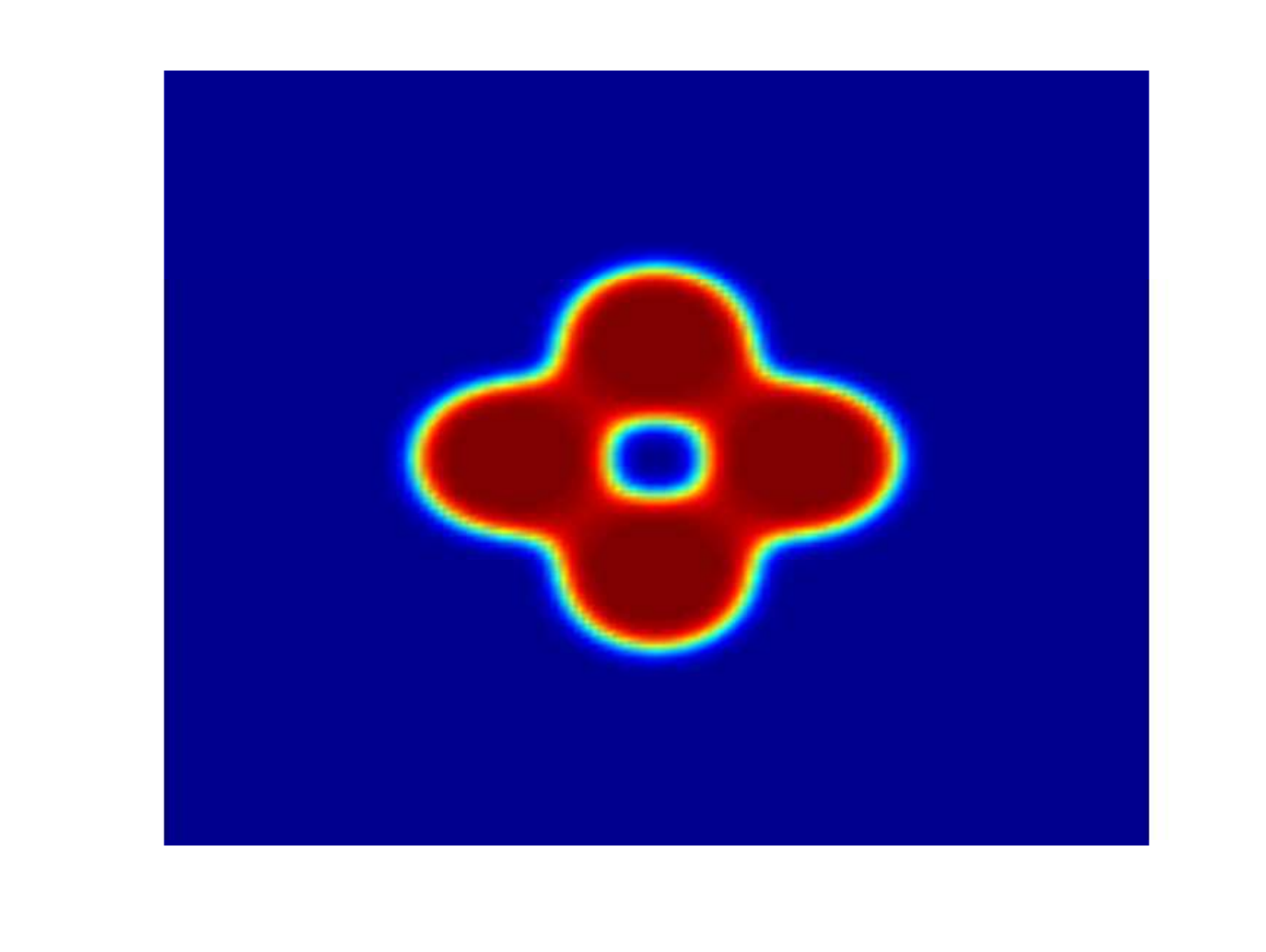}
\includegraphics[width=1.47in]{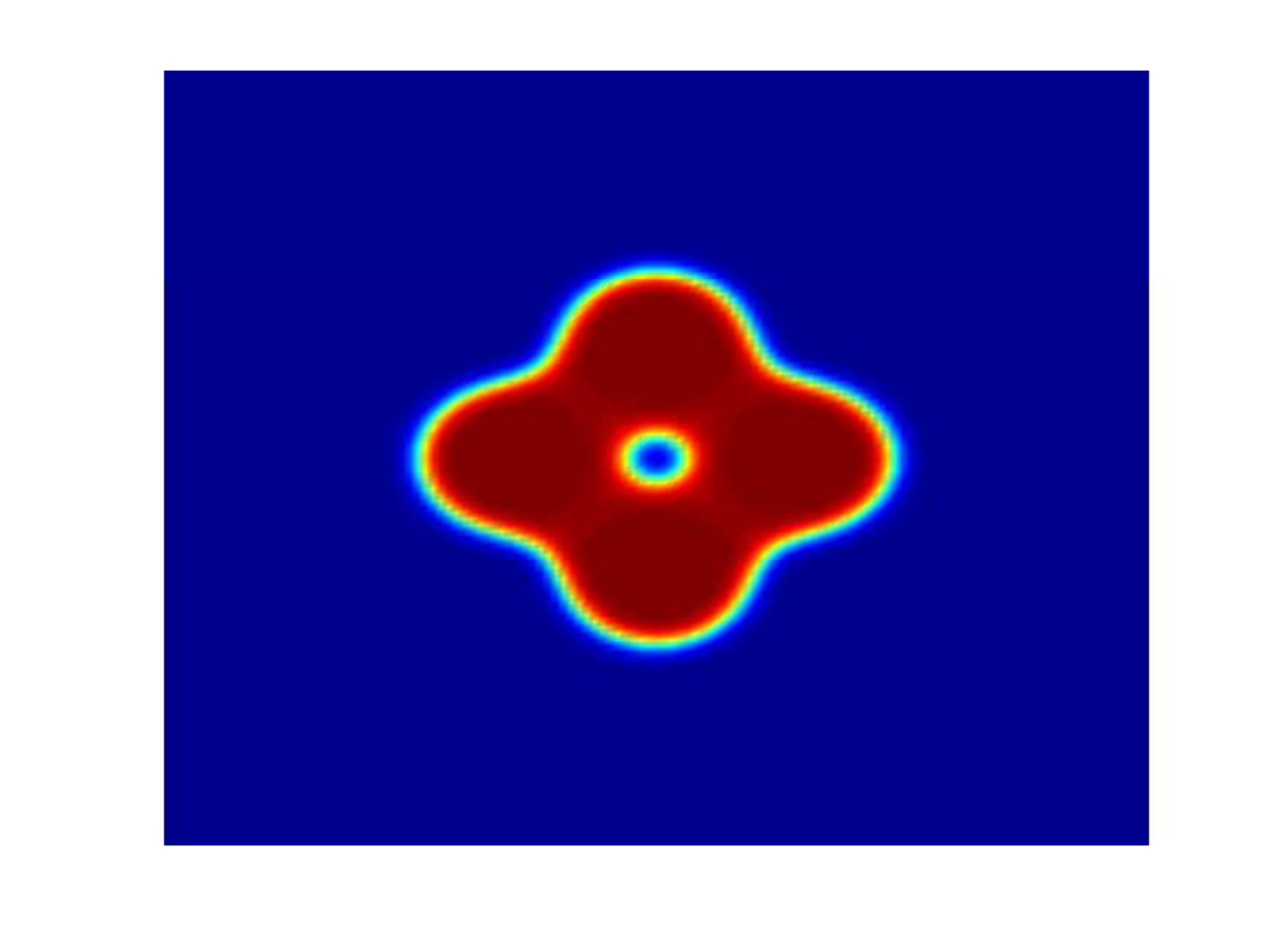}
\includegraphics[width=1.47in]{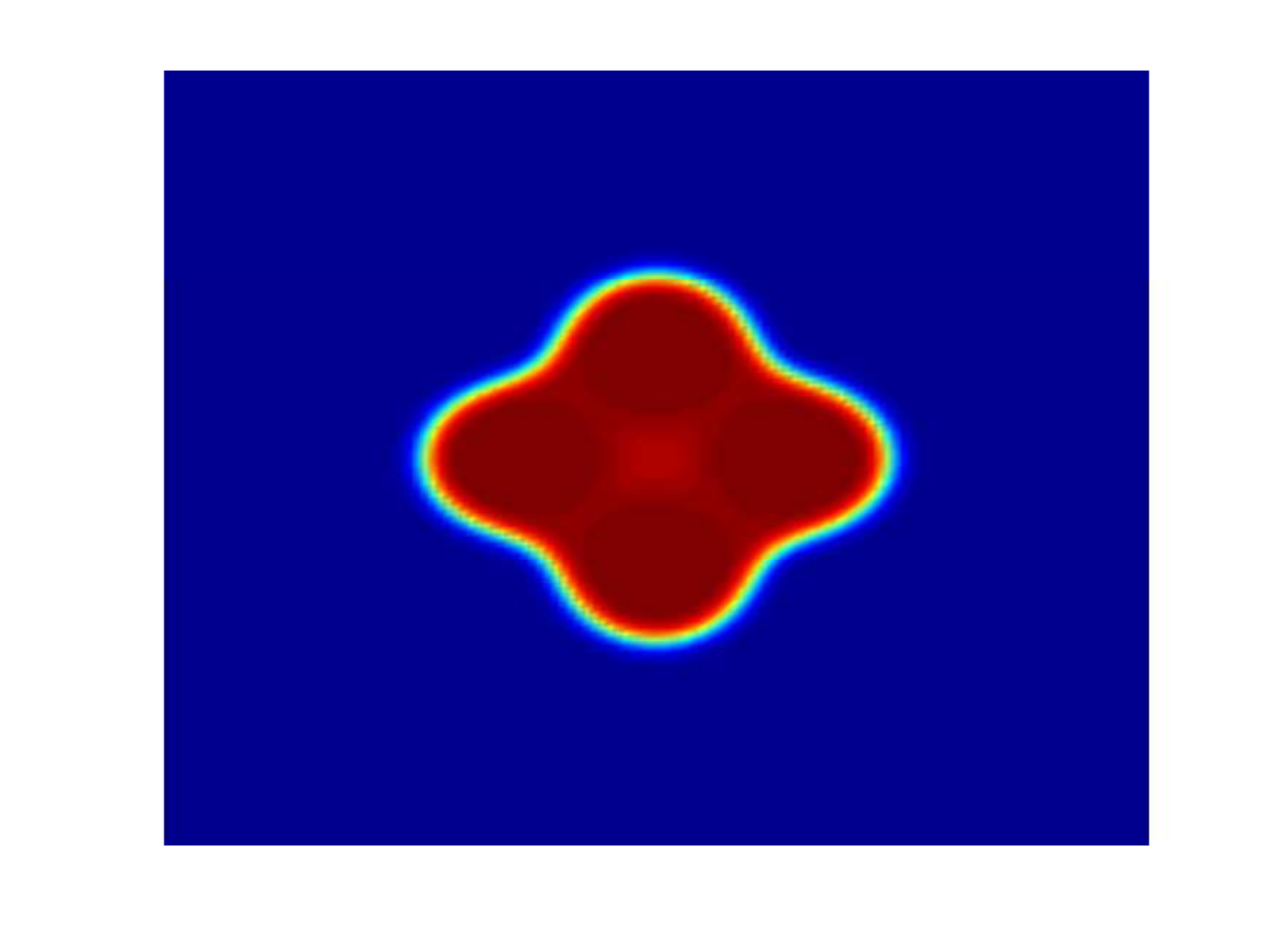}\\
\includegraphics[width=1.47in]{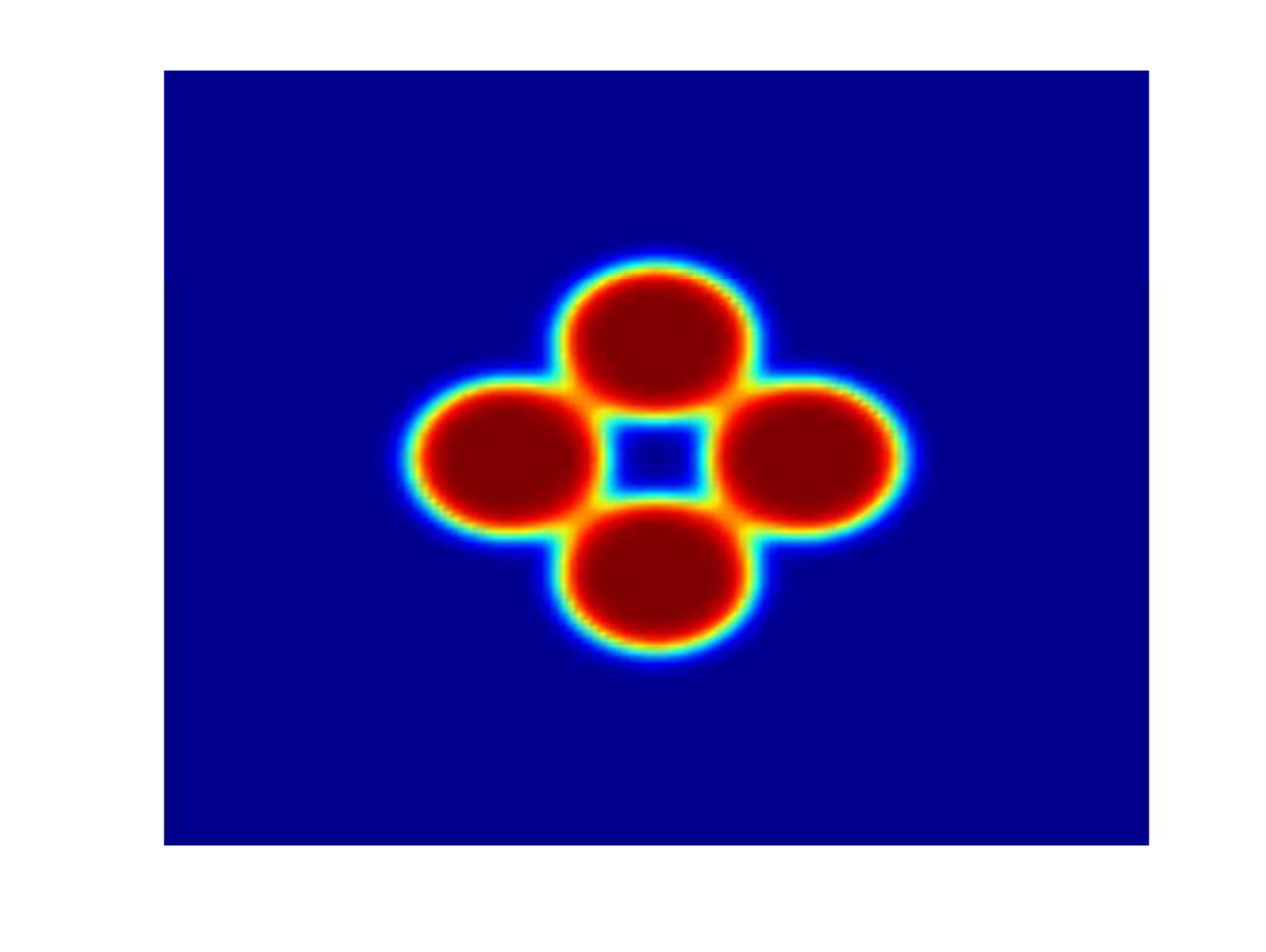}
\includegraphics[width=1.47in]{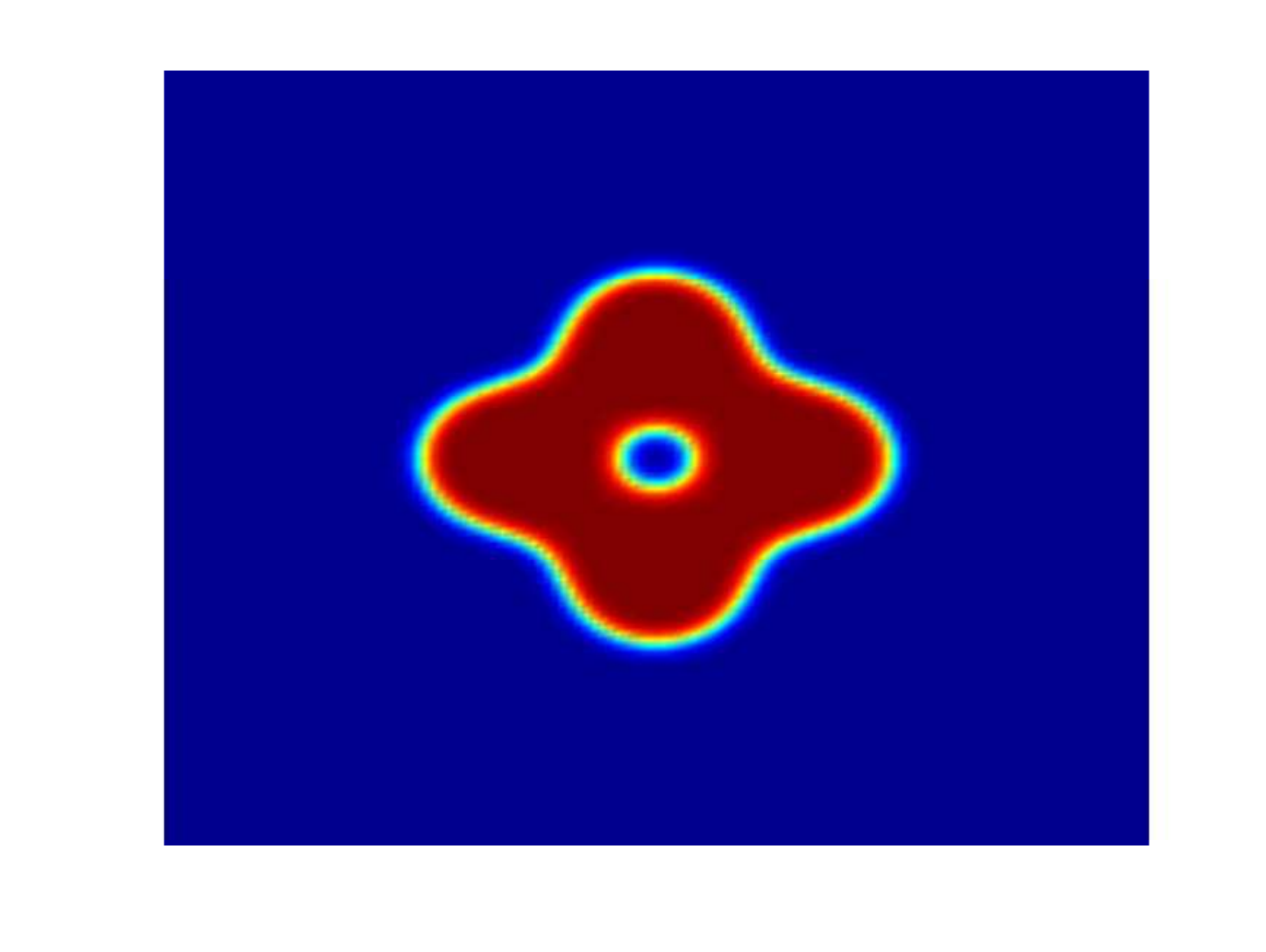}
\includegraphics[width=1.47in]{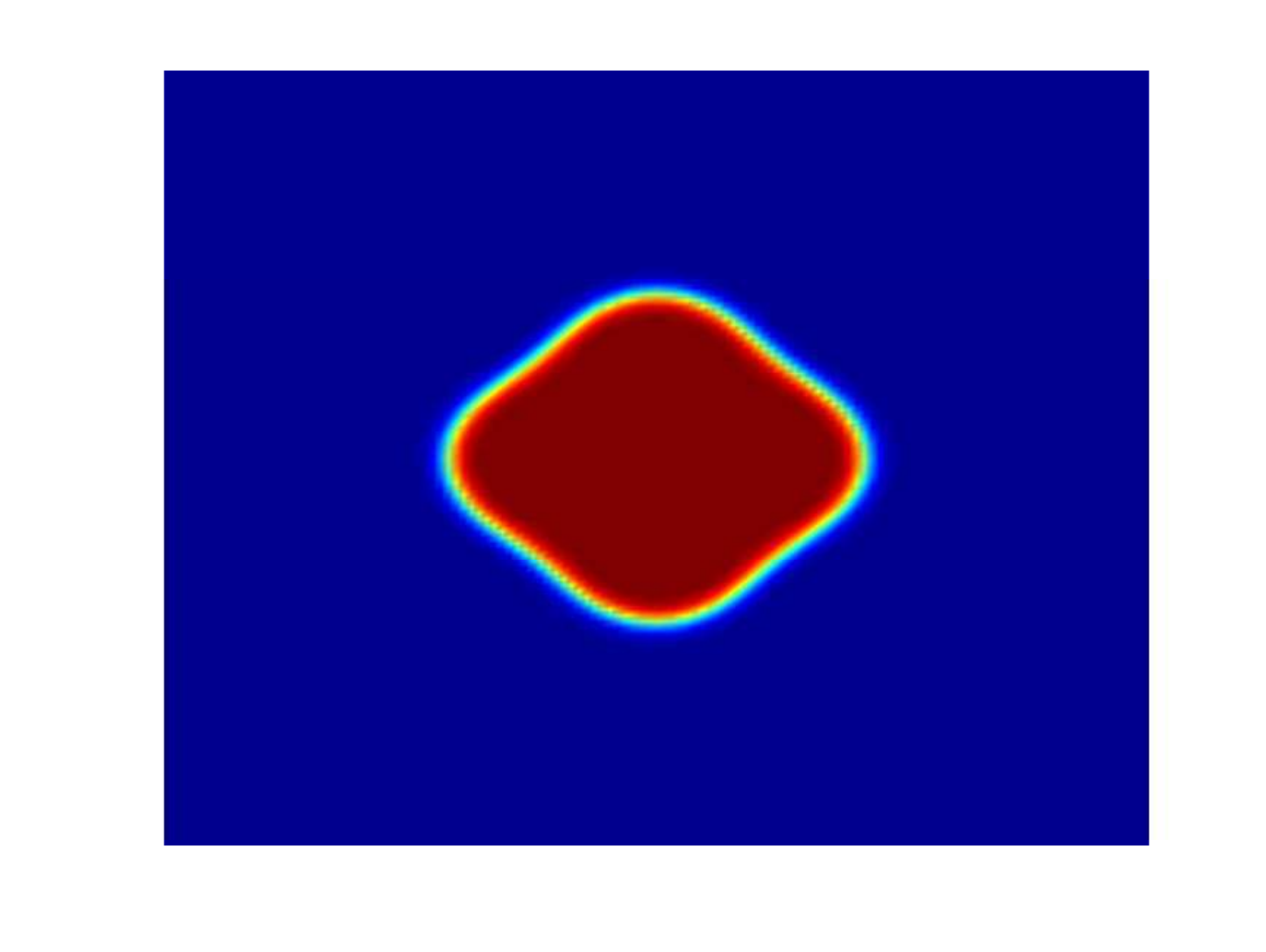}
\includegraphics[width=1.47in]{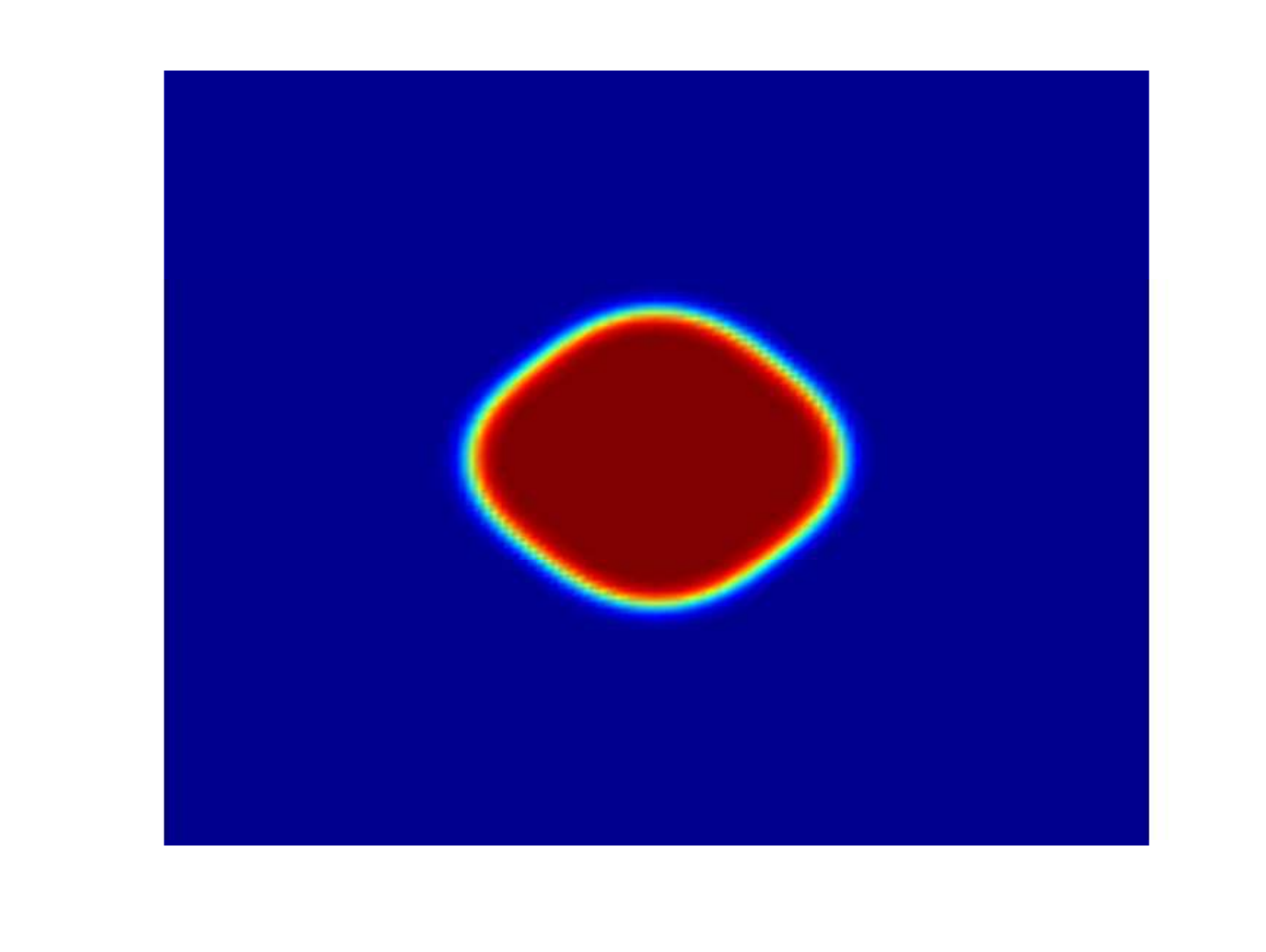}\\
\includegraphics[width=1.47in]{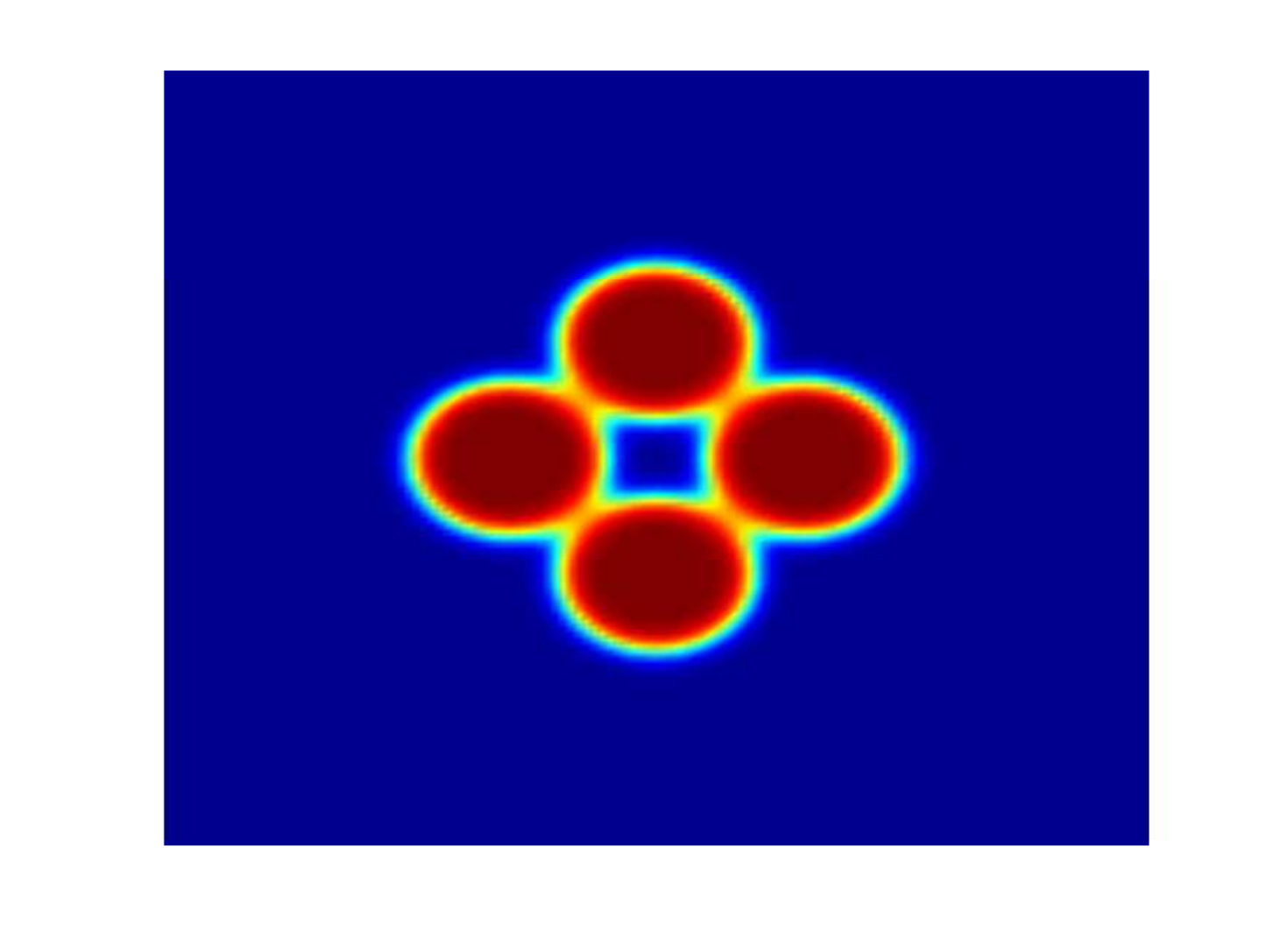}
\includegraphics[width=1.47in]{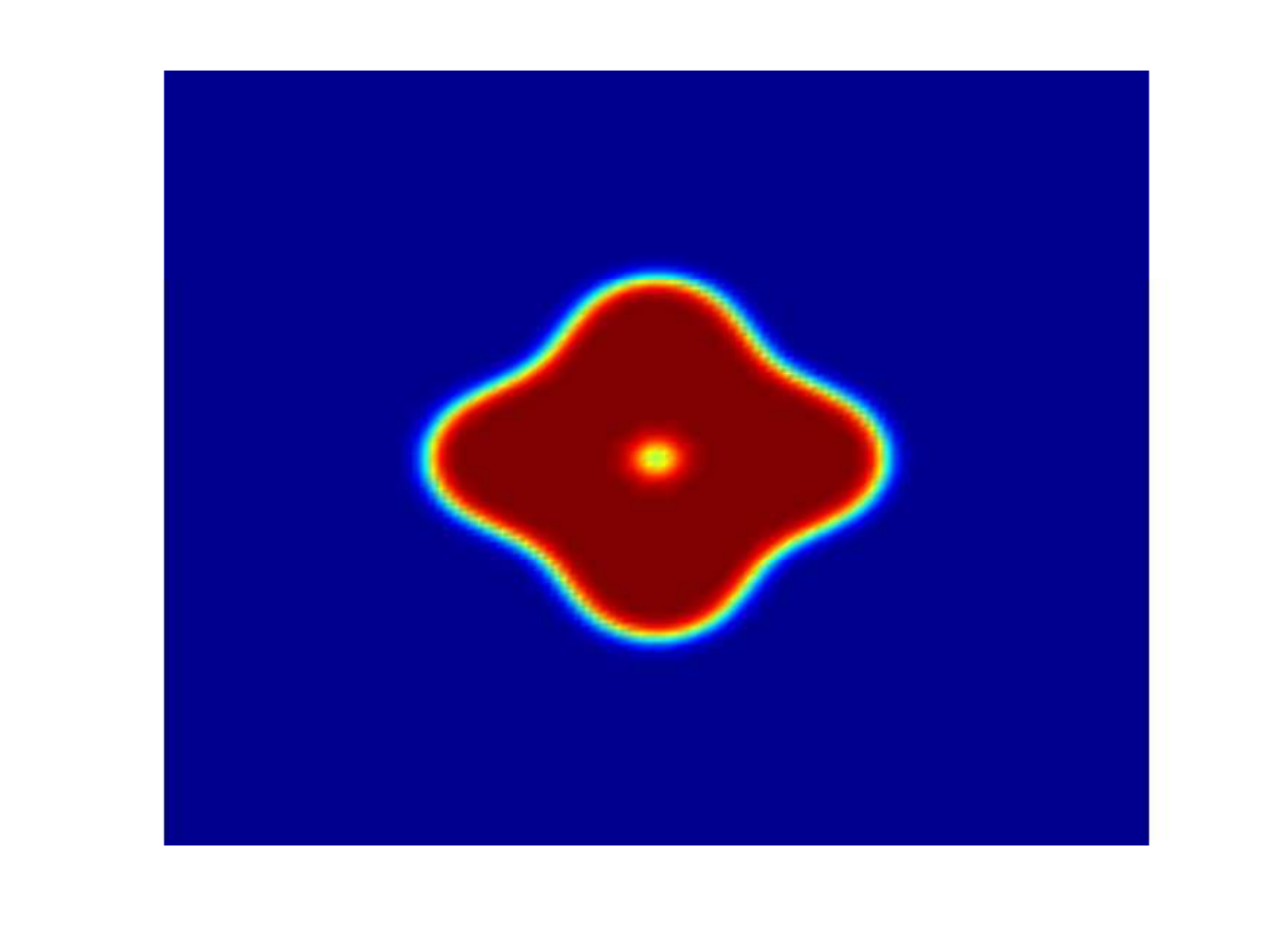}
\includegraphics[width=1.47in]{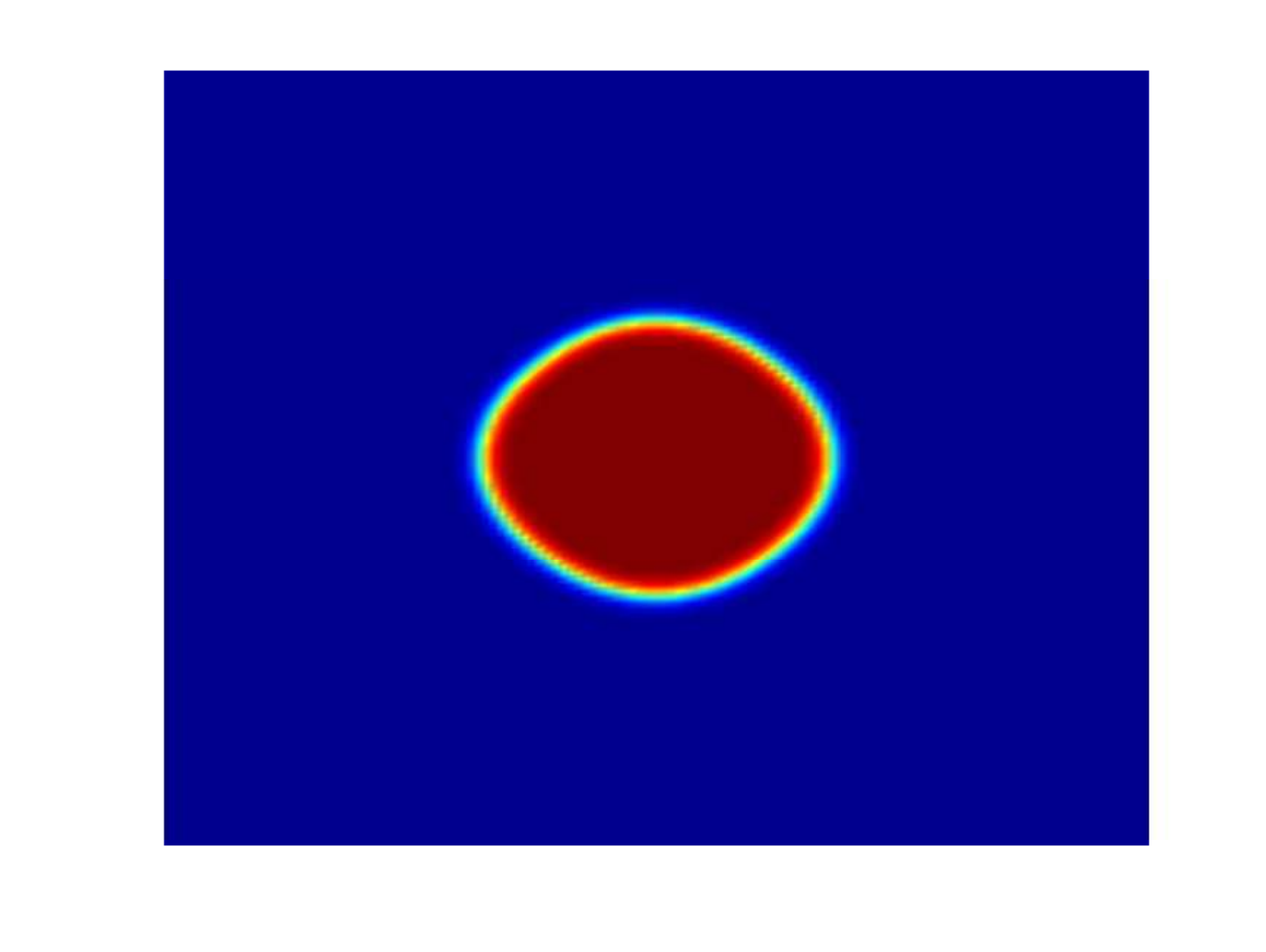}
\includegraphics[width=1.47in]{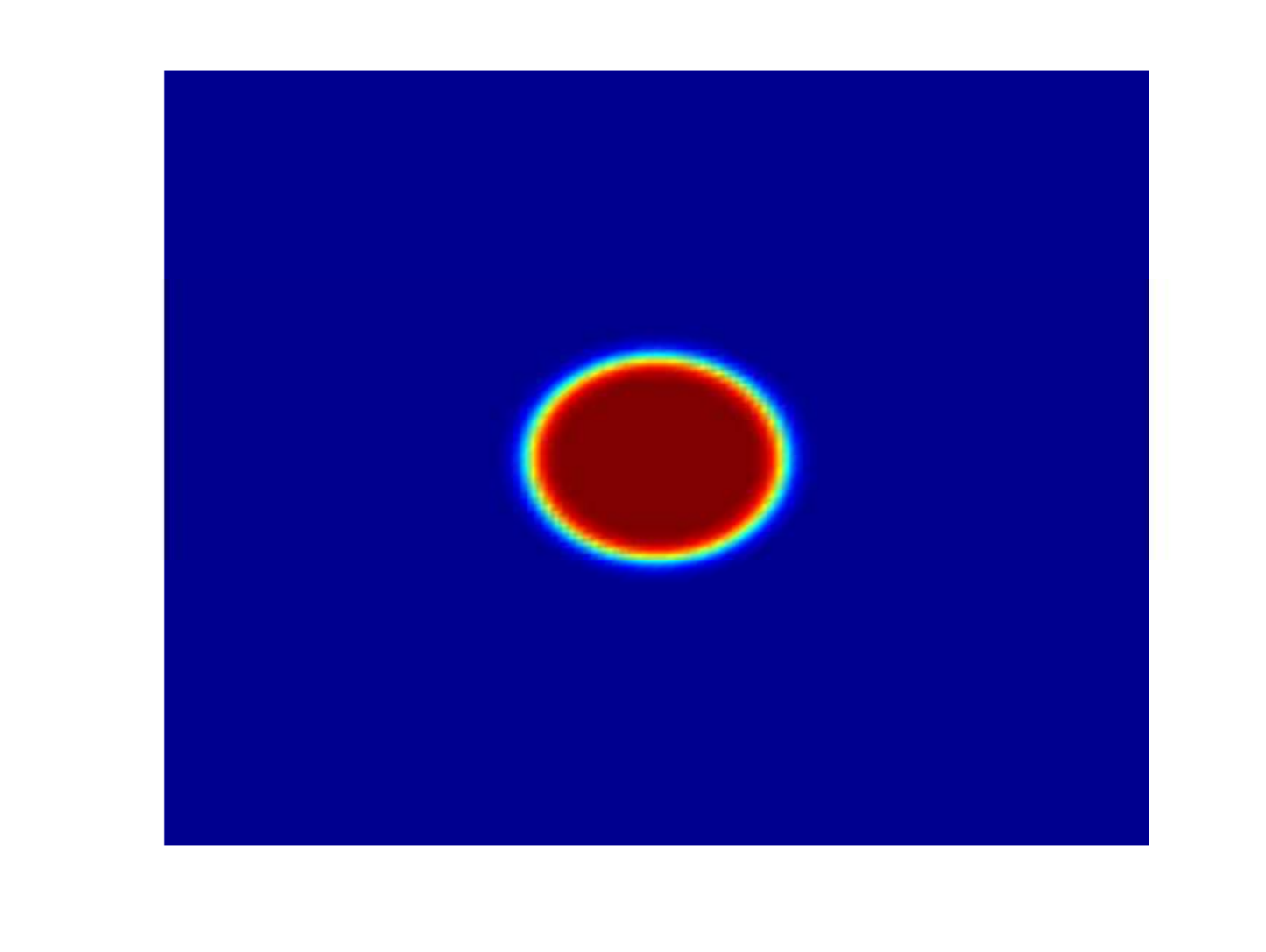}\\
\caption{Solution snapshots of Allen-Cahn equation \eqref{Problem-1} at
  $t=1, 30, 100, 200$ (from left to right)
  for three fractional order $\alpha=0.4,\,0.7$ and $0.9$ (from top to bottom), respectively.}
\label{Classical-AC-Drops}
\end{figure}
%%%%%%%%%%%%%%%%%%%%%%%%%%%%%%%%%%%%%%%%%%%%%%%%%%%%%%%%%%%%%%%%%%%%%%%%%%%%%%%%%%%%%

%%%%%%%%%%%%%%%%%%%%%%%%%%%%%%%%%%%%%%%%%%%%%%%%%%%%%%%%%%%%%%%%%%%%%%%%%%%%%%%%%%%%%%%%%%%%%
\begin{figure}[htb!]
\centering
\includegraphics[width=3.0in,height=2.0in]{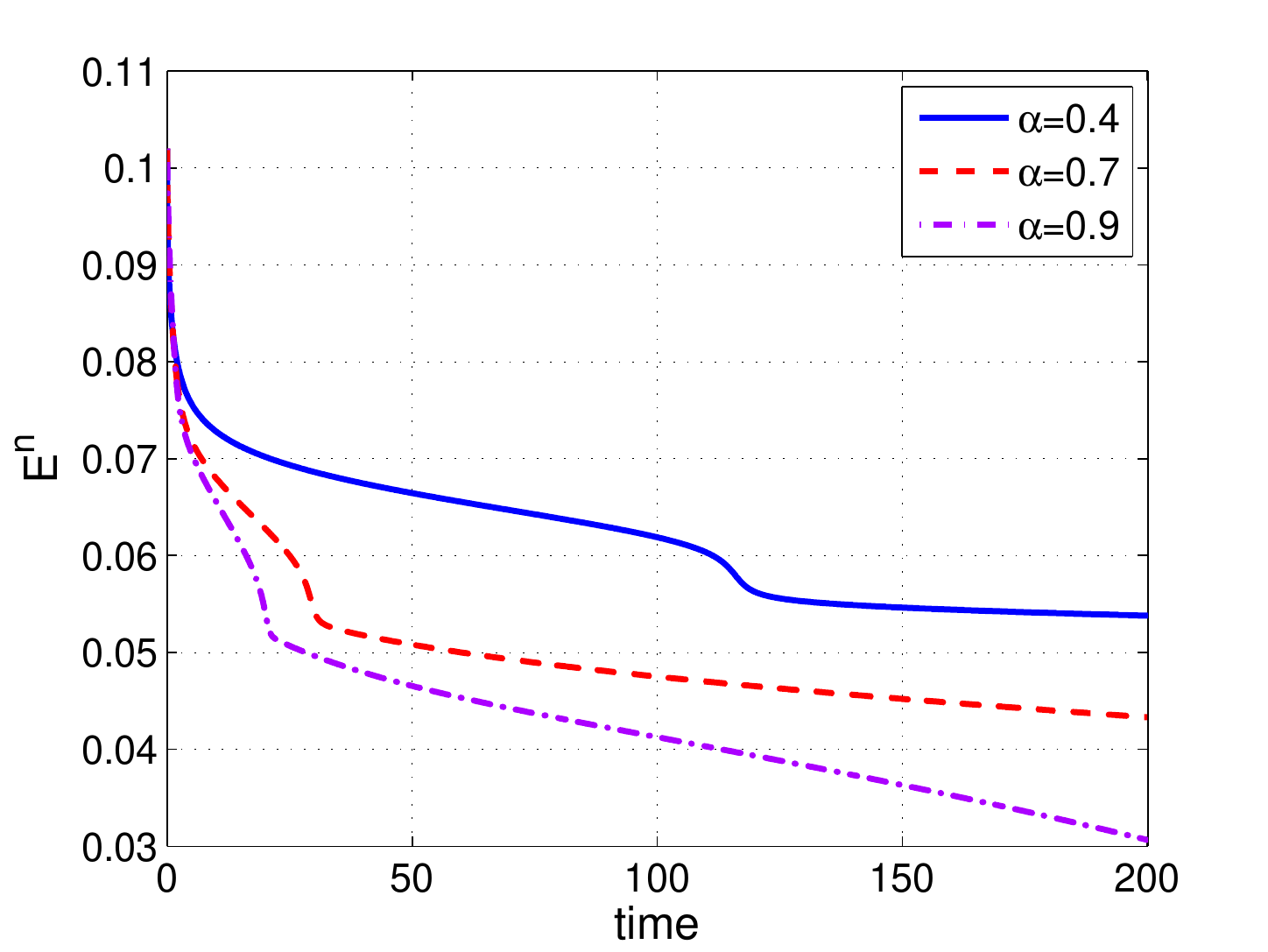}
\includegraphics[width=3.0in,height=2.0in]{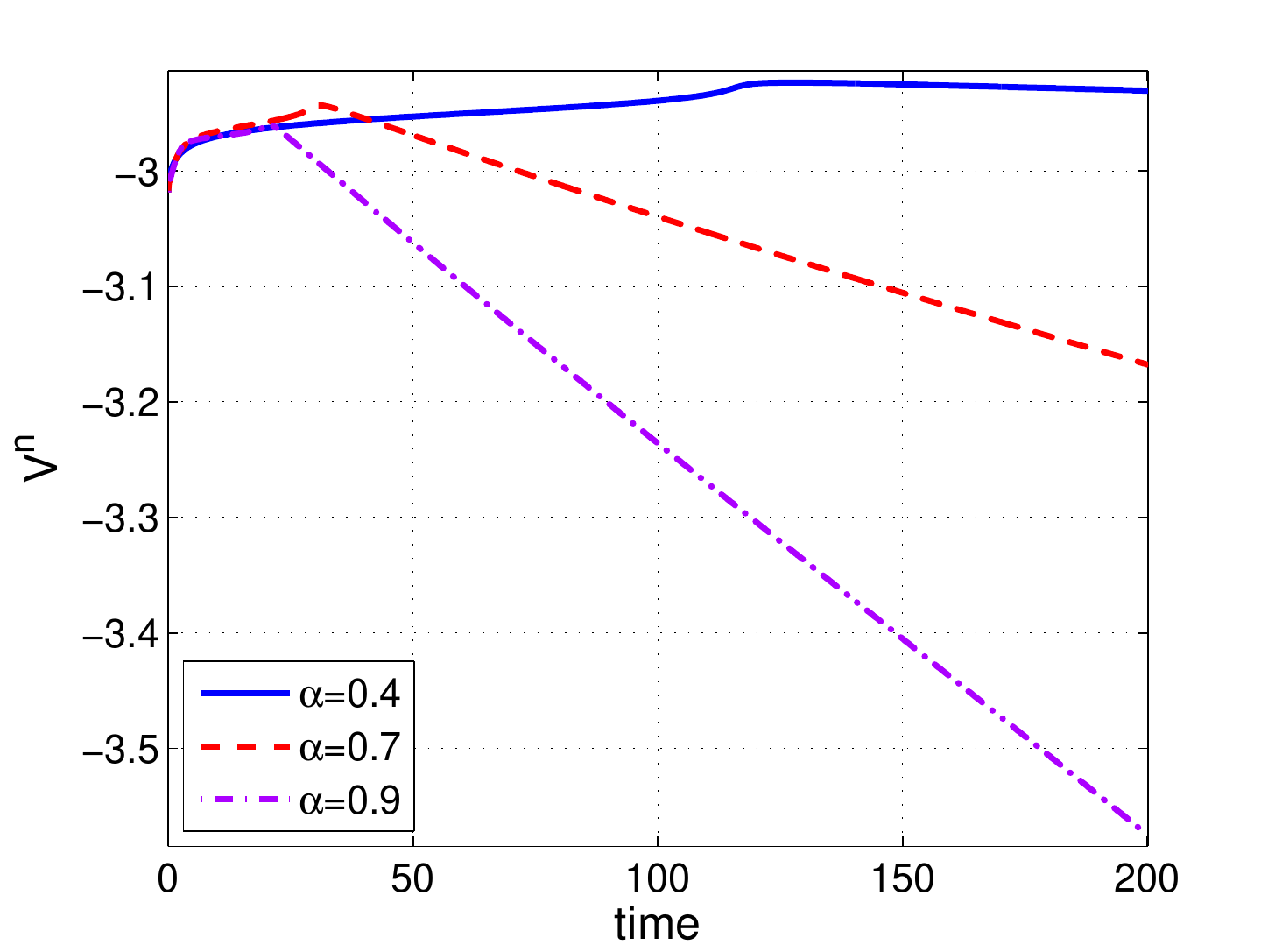}
\caption{Evolutions of energy and volume curves (from left to right) for
  the time-fractional Allen-Cahn equation \eqref{Problem-1} with fractional order $\alpha=0.4,\,0.7$ and $0.9$, respectively.}
\label{Comparison-Classical-AC-Energy-Mass}
\end{figure}
%%%%%%%%%%%%%%%%%%%%%%%%%%%%%%%%%%%%%%%%%%%%%%%%%%%%%%%%%%%%%%%%%%%%%%%%%%%%%%%%%%%%%

We examine different time-stepping approaches for simulating
the model \eqref{Problem-2} until the final time $T=30$ with a fractional order $\alpha=0.9$.
Always, put $T_{0}=0.01$, $N_{0}=30$ and $\gamma=3$ in the starting cell $[0, T_{0}]$.
We consider the \emph{Grade Step} approach using the uniform mesh with $N_1=2970$,
and the \emph{Adaptive Step} approach with parameters
$\kappa=10^{6}$, $\tau_{\min}=\tau_{N_{0}}=10^{-3}$
and $\tau_{\max}=10^{-1}$.
From Figure \ref{Comparison-Adaptive-Energy-Curves},
the discrete energy curves generated by using the adaptive time steps practically
coincide with those by \emph{Grade Step} approach, for both the CN-IEQ and CN-SAV methods.
As expected, the \emph{Adaptive Step} approach uses small time steps when
the energy dissipates fast, and generates large time steps otherwise.
In the remainder interval $(T_0,T]$, we put 2970 points on the uniform mesh,
while the total number of adaptive time steps are 667.
So the adaptive time-stepping strategy is computationally efficient.

%%%%%%%%%%%%%%%%%%%%%%%%%%%%%%%%%%%%%%%%%%%%%%%%%%%%%%%%%%%%%%%%%%%%%%%%%%%%%%%%%%%%%%%%%%
\begin{figure}[htb!]
\centering
\includegraphics[width=1.47in]{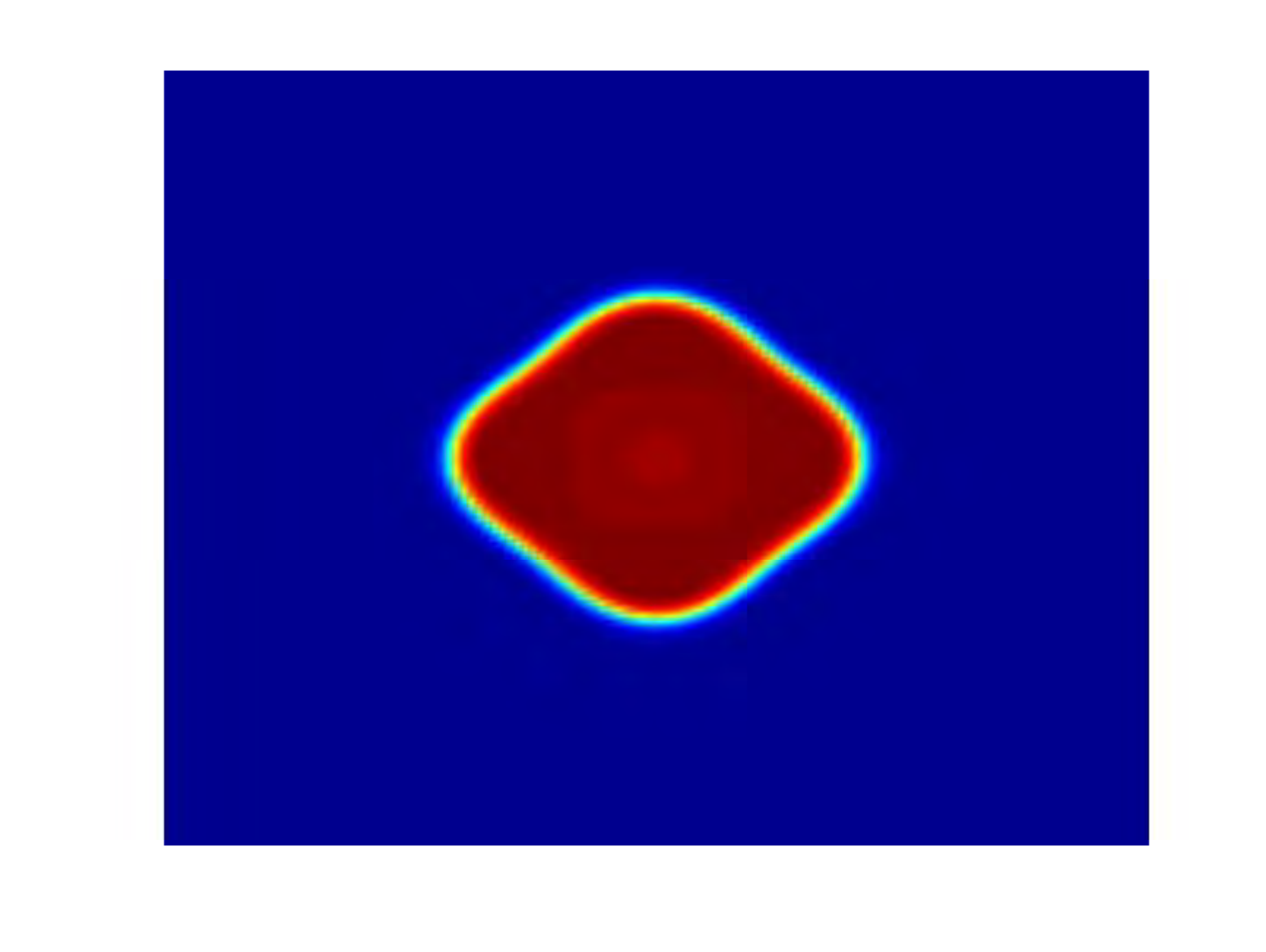}
\includegraphics[width=1.47in]{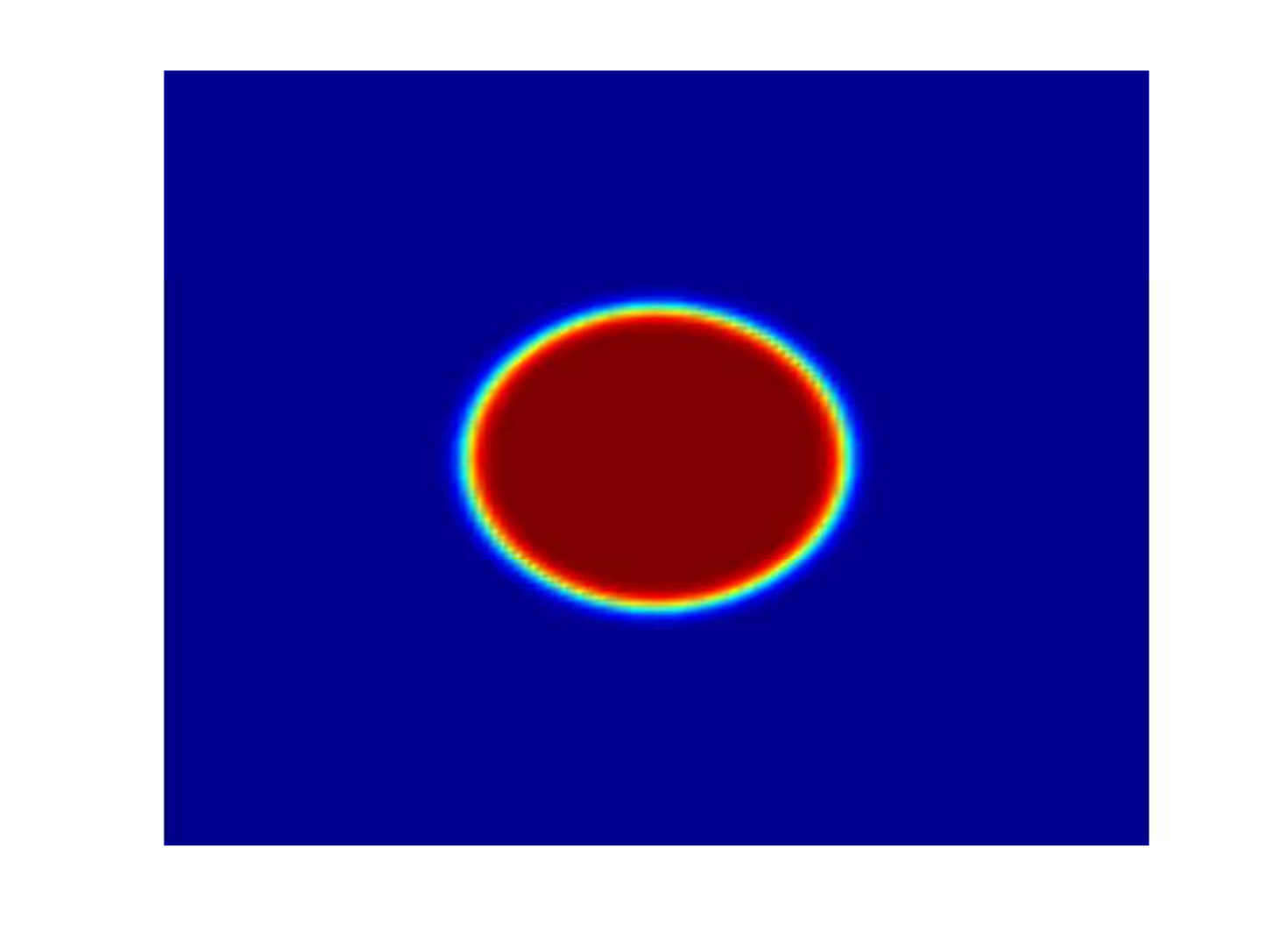}
\includegraphics[width=1.47in]{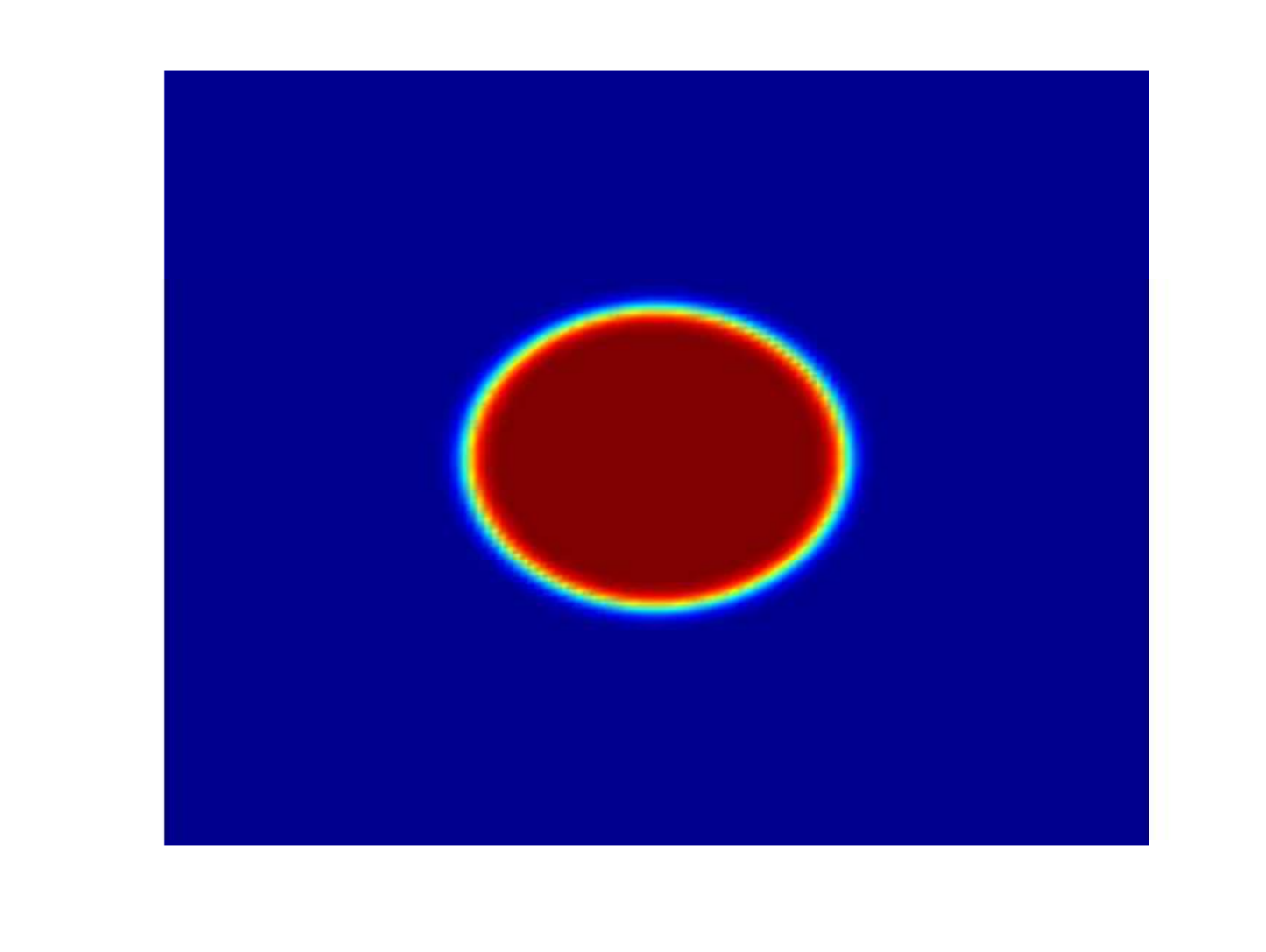}
\includegraphics[width=1.47in]{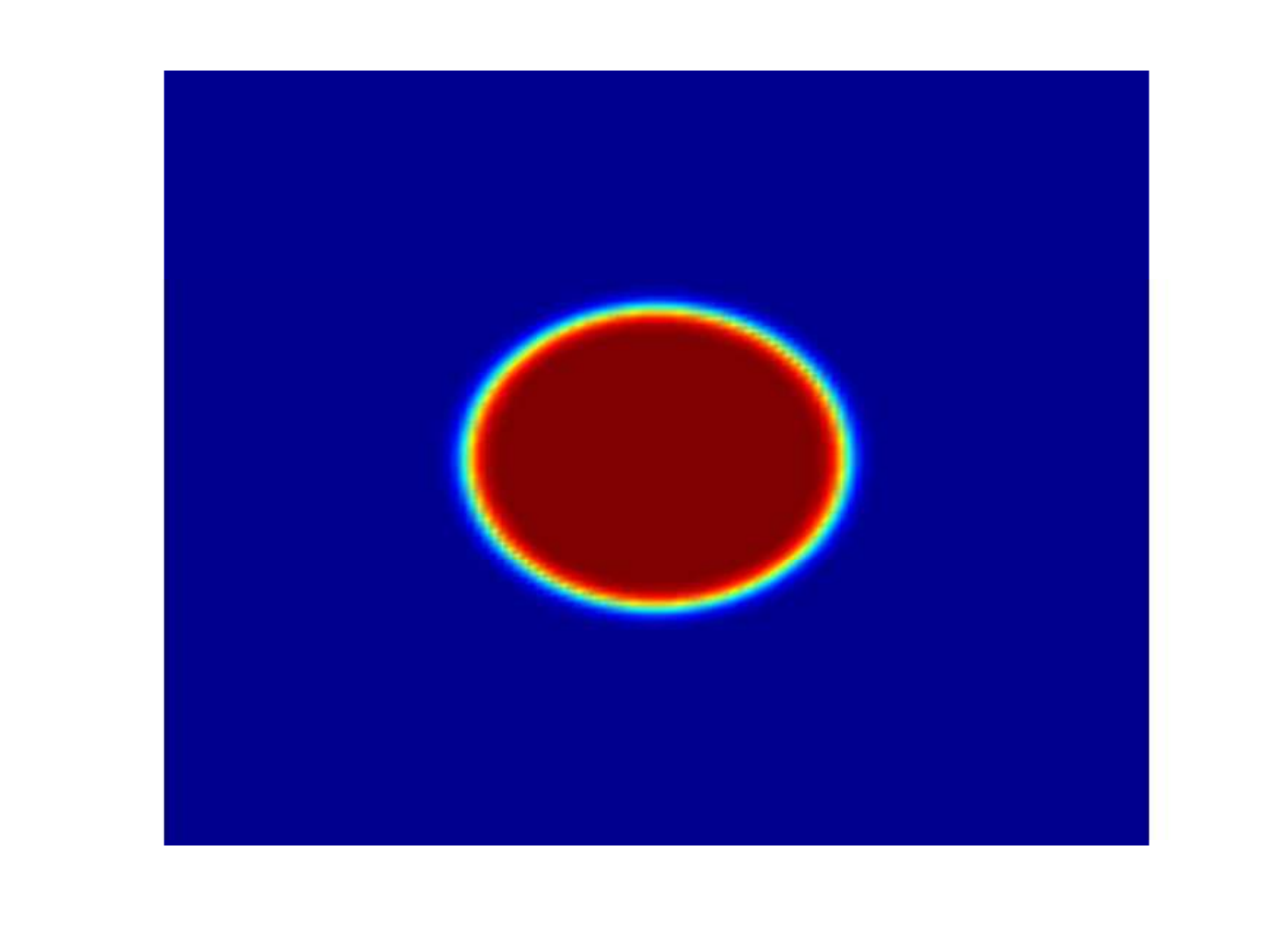}\\
\includegraphics[width=1.47in]{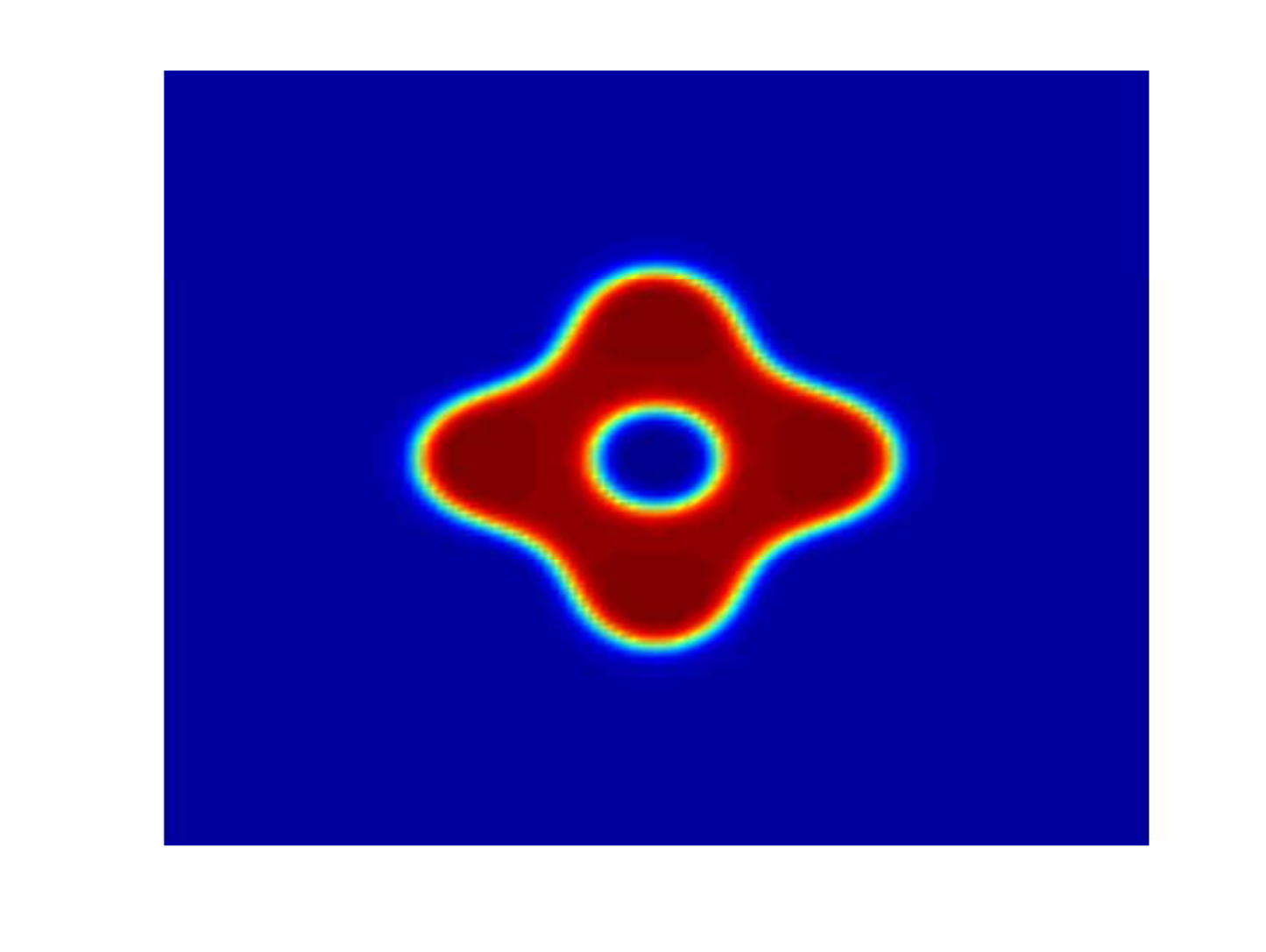}
\includegraphics[width=1.47in]{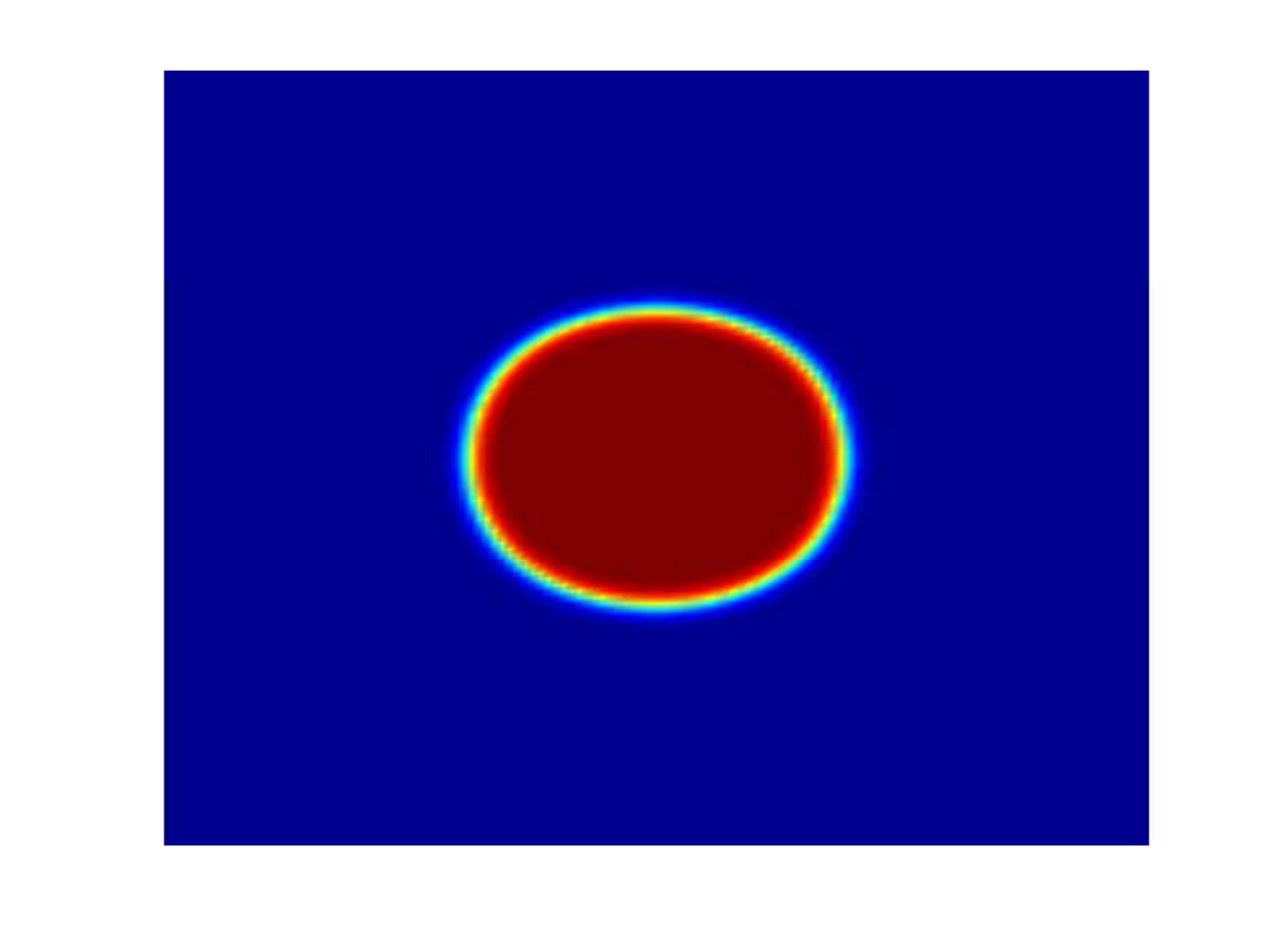}
\includegraphics[width=1.47in]{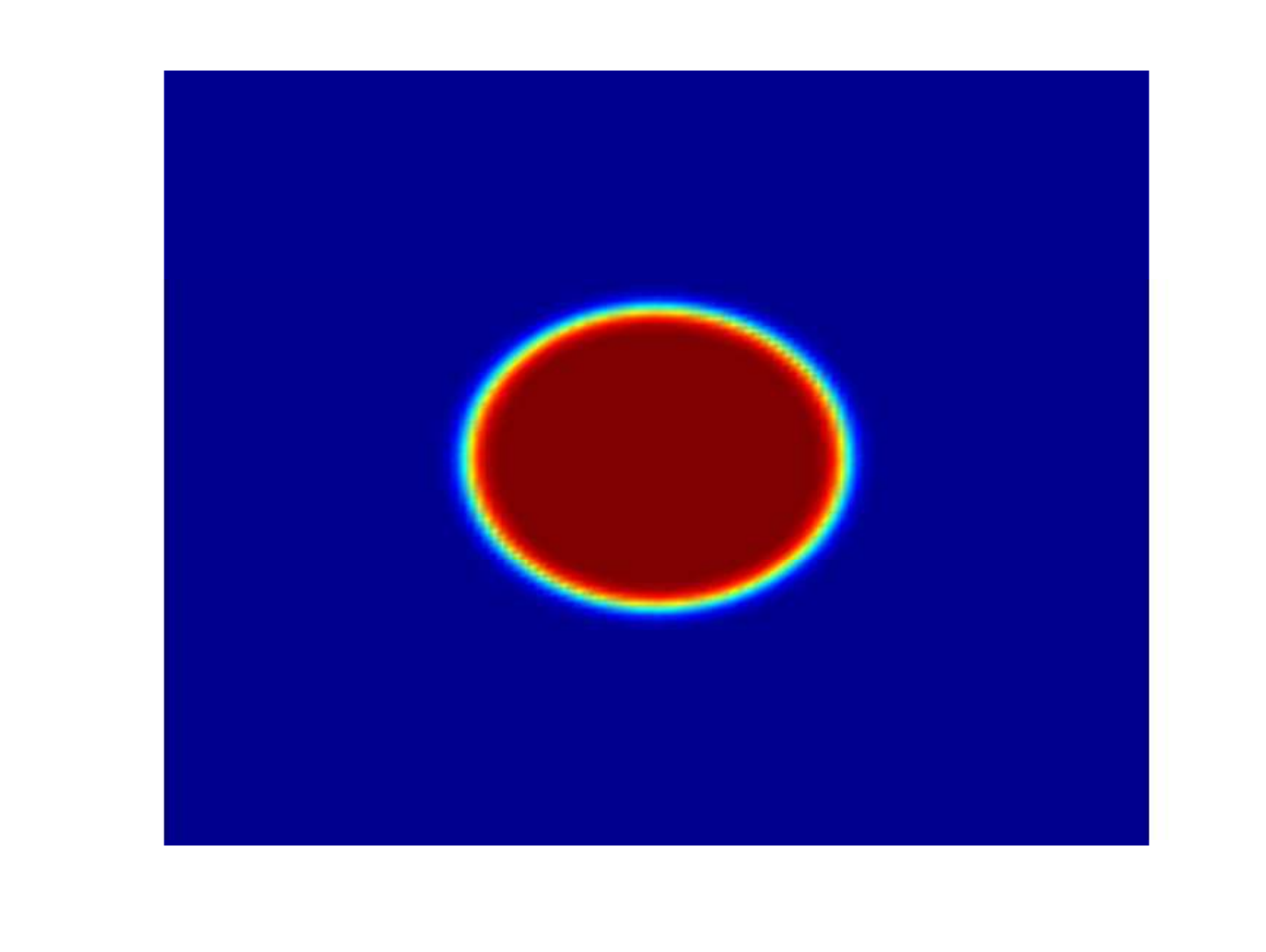}
\includegraphics[width=1.47in]{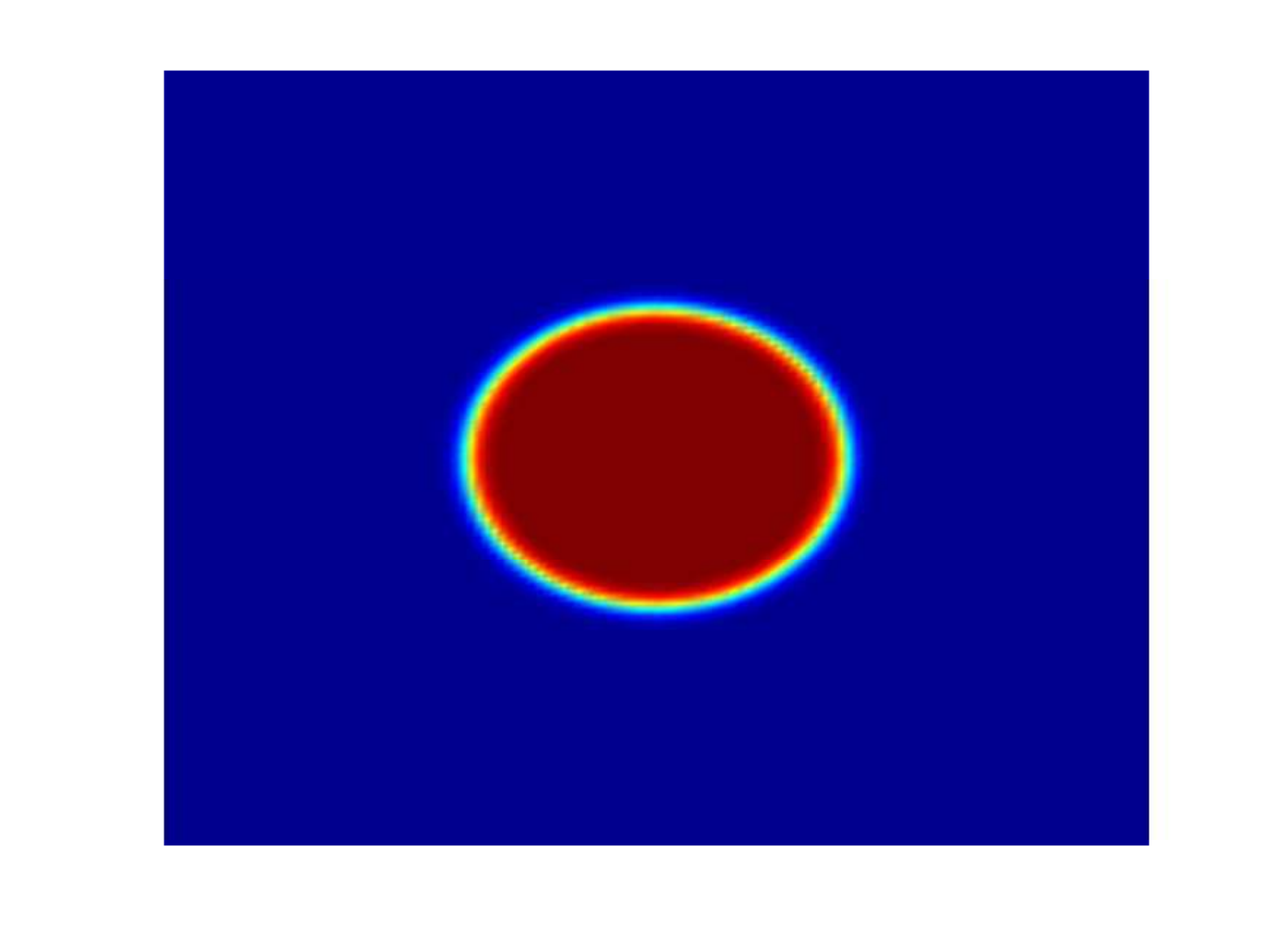}\\
\includegraphics[width=1.47in]{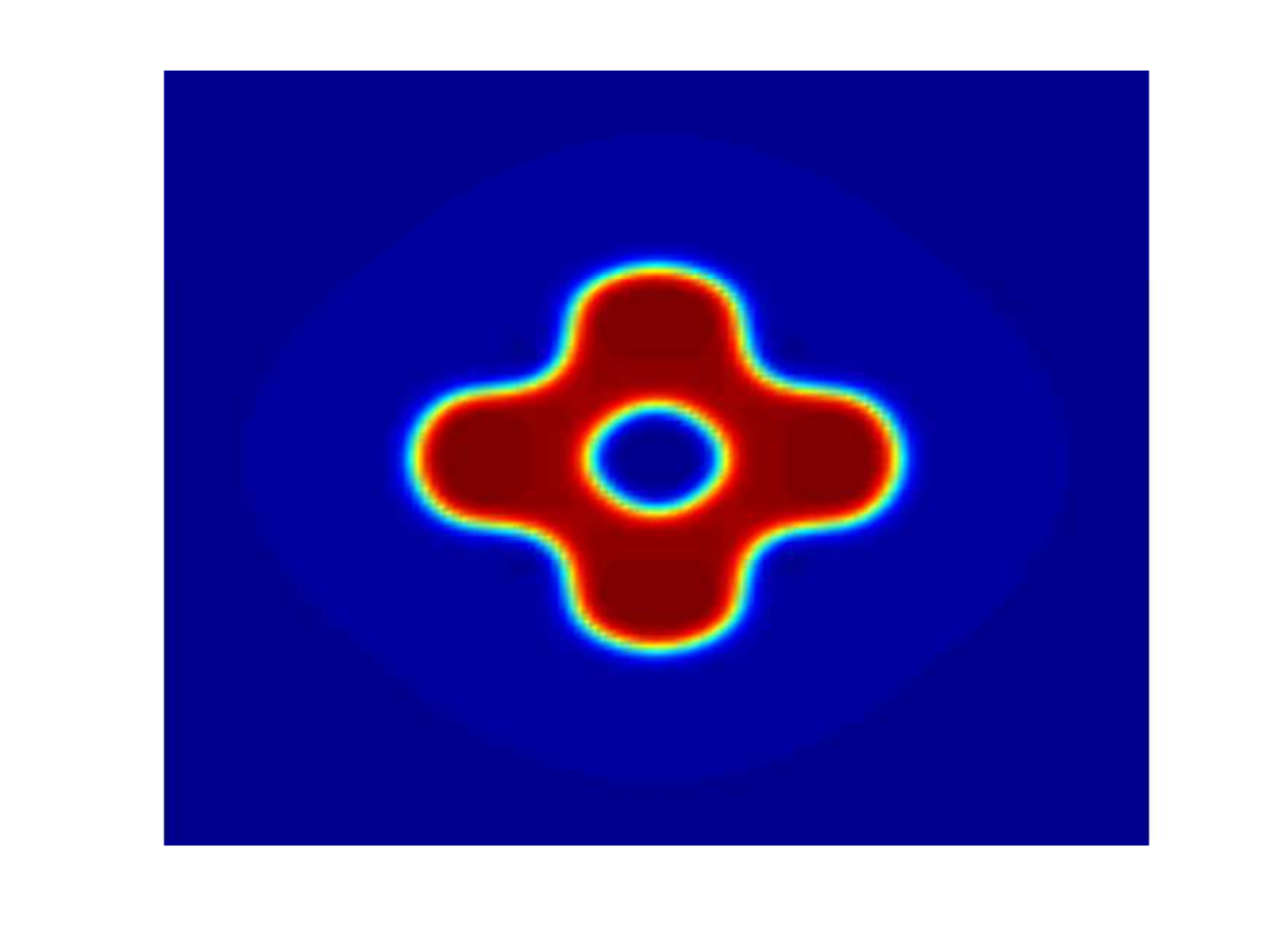}
\includegraphics[width=1.47in]{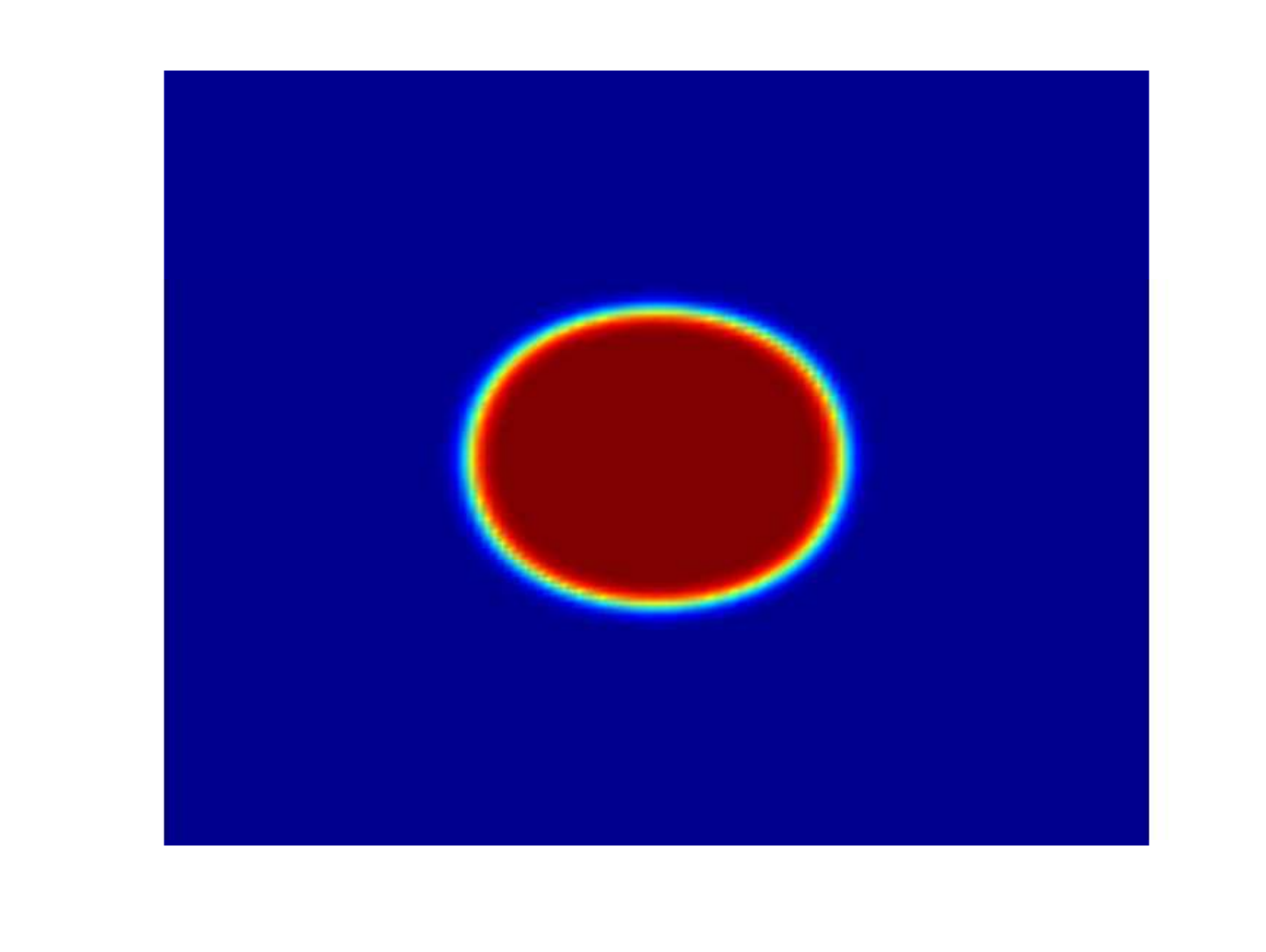}
\includegraphics[width=1.47in]{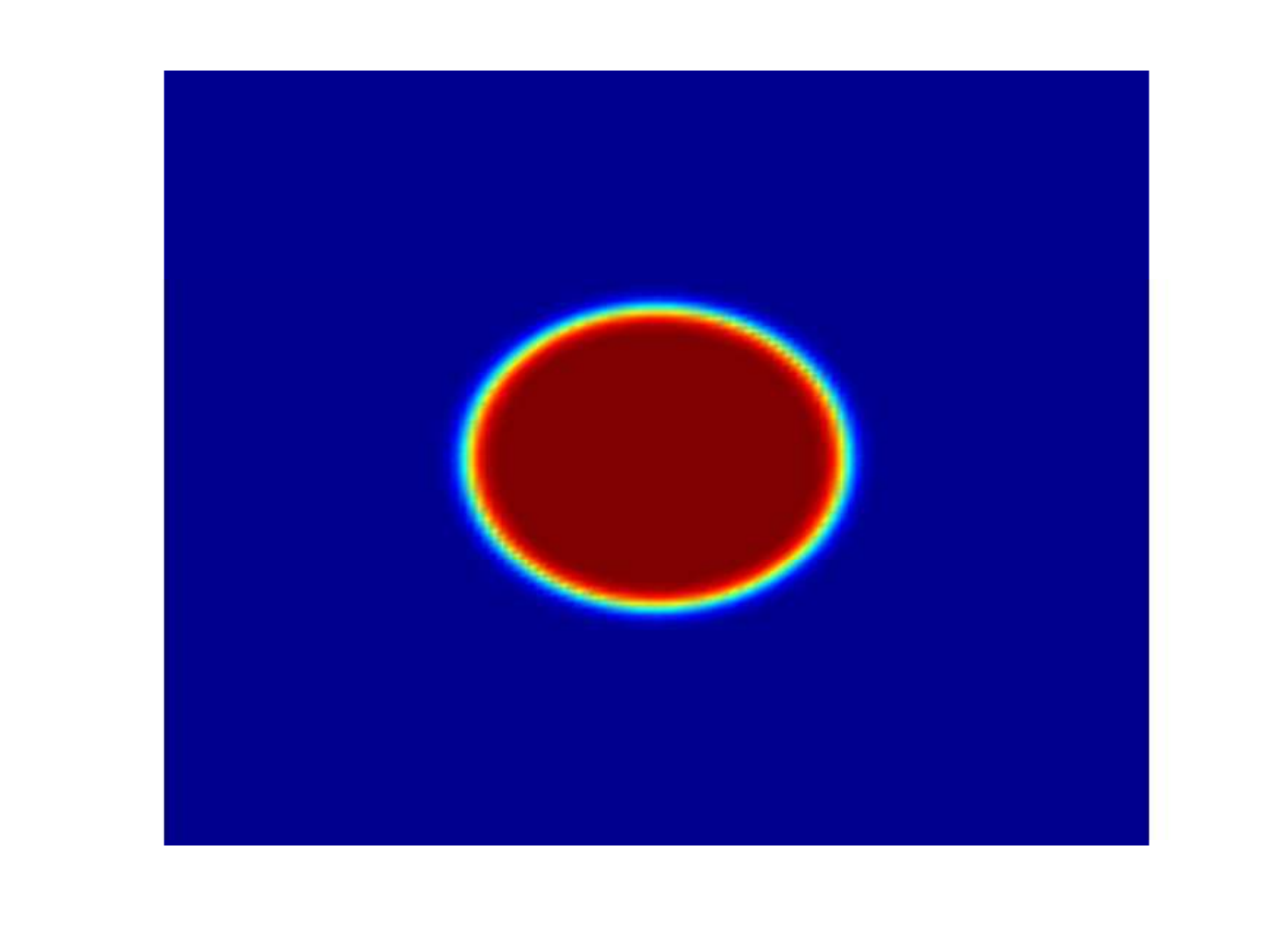}
\includegraphics[width=1.47in]{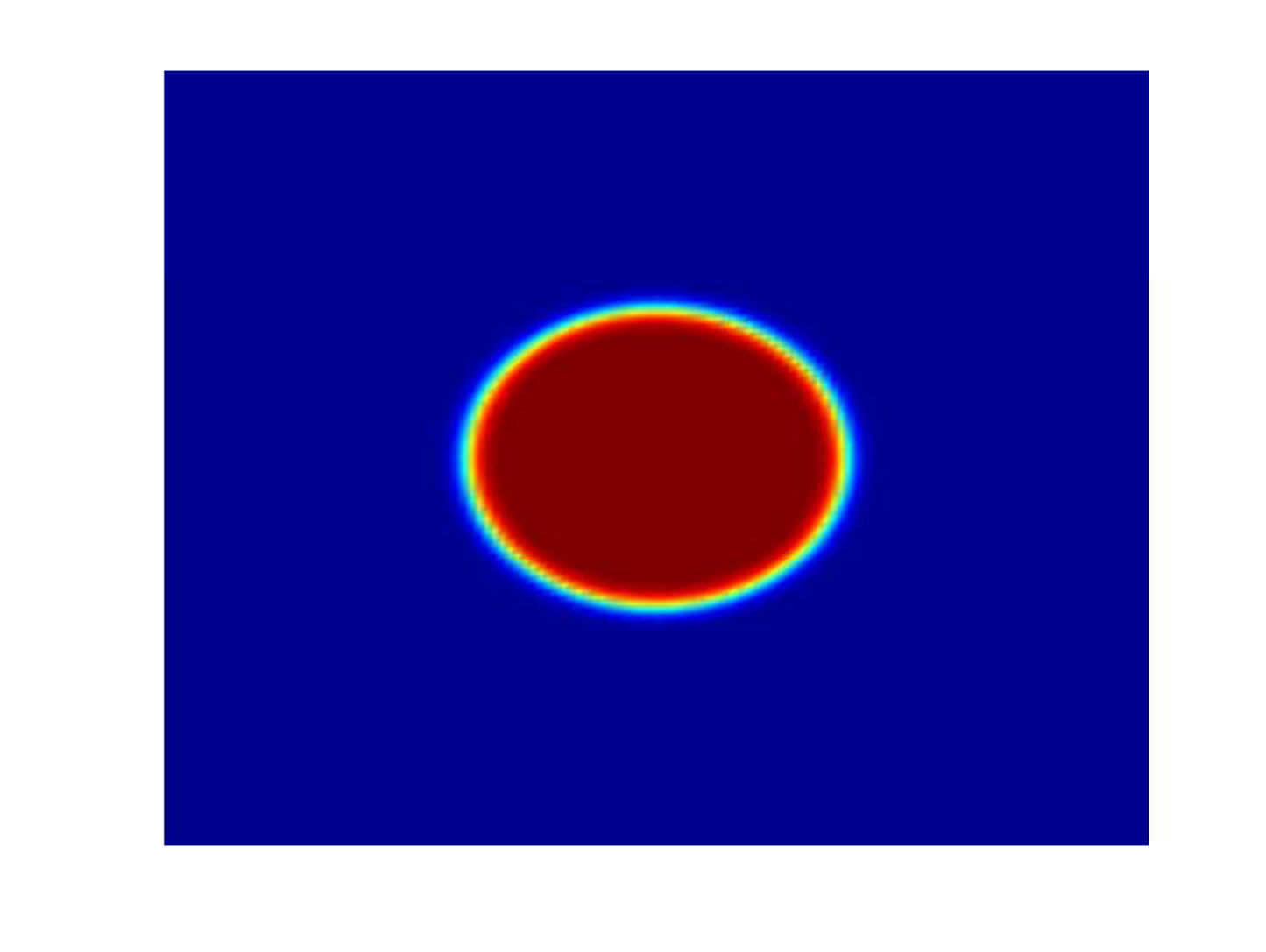}\\
\caption{Solution snapshots of the Cahn-Hilliard equation \eqref{Problem-3}
  at $t=0.01, 30, 100, 200$ (from left to right)
  for fractional order $\alpha=0.4,\,0.7$ and $0.9$ (from top to bottom), respectively.}
\label{CH-Drops}
\end{figure}
%%%%%%%%%%%%%%%%%%%%%%%%%%%%%%%%%%%%%%%%%%%%%%%%%%%%%%%%%%%%%%%%%%%%%%%%%%%%%%%%%%%%%

%%%%%%%%%%%%%%%%%%%%%%%%%%%%%%%%%%%%%%%%%%%%%%%%%%%%%%%%%%%%%%%%%%%%%%%%%%%%%%%%%%%%%%%%%%%%%
\begin{figure}[htb!]
\centering
\includegraphics[width=3.0in,height=2.0in]{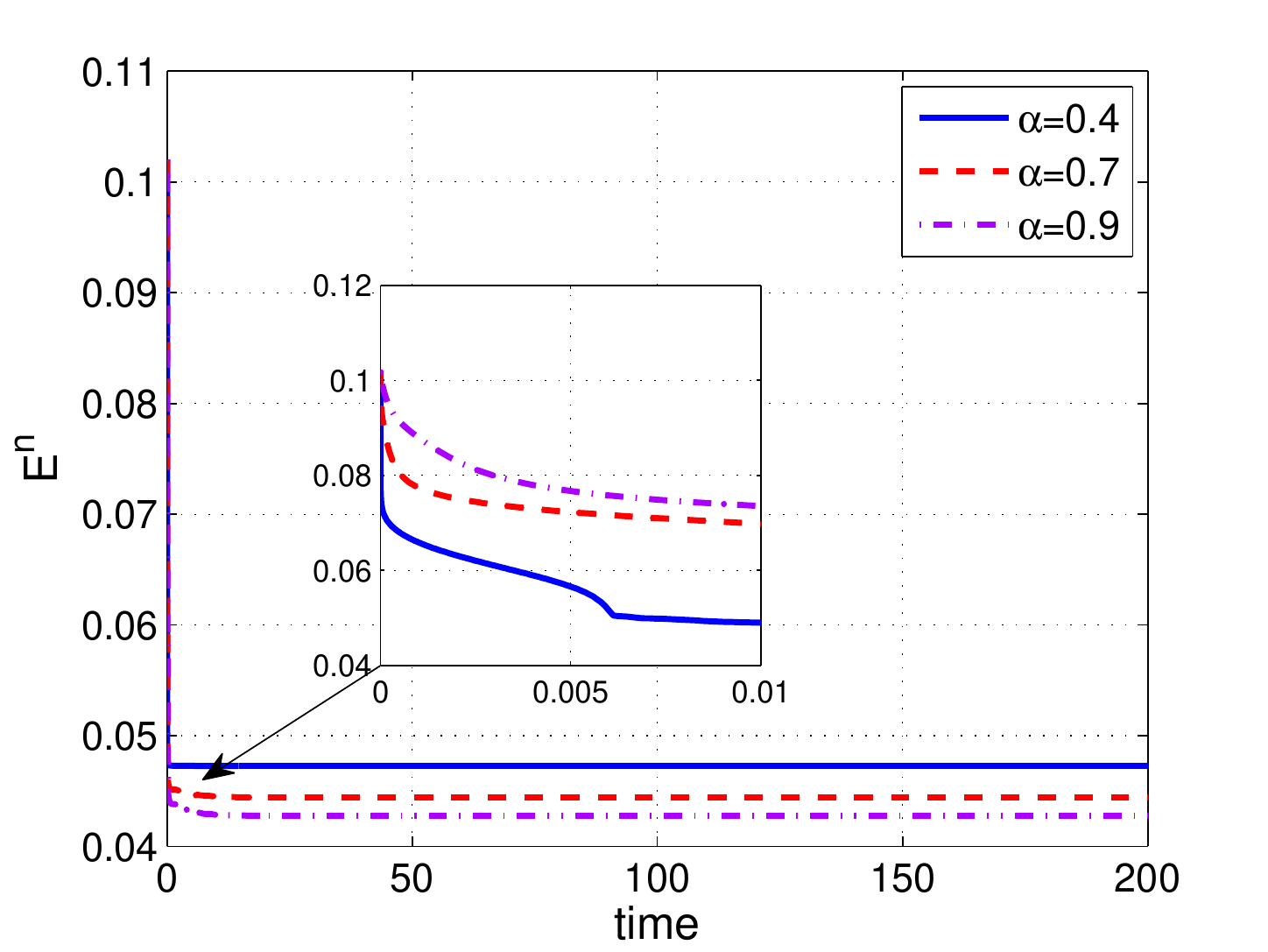}
\includegraphics[width=3.0in,height=2.0in]{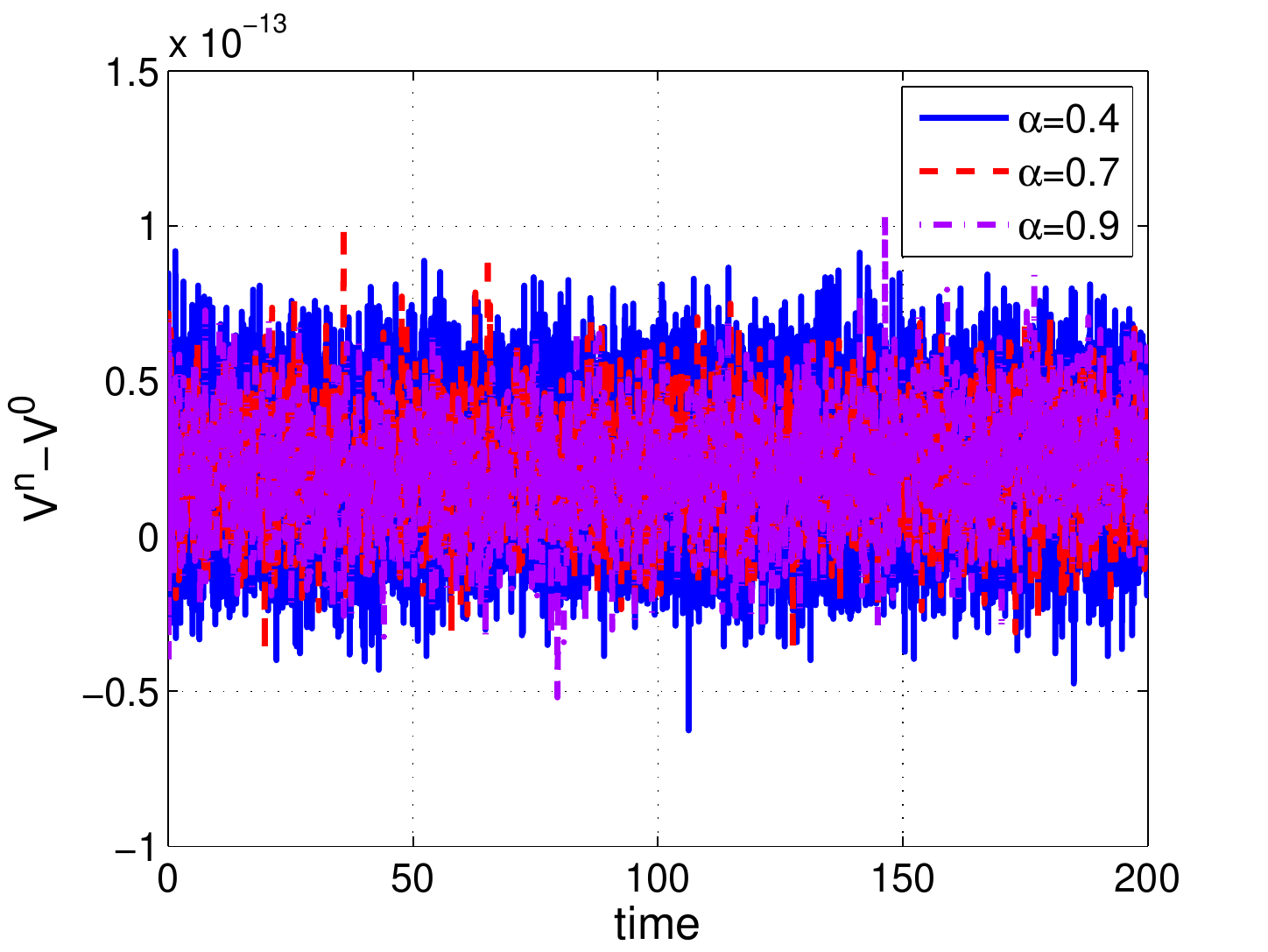}
\caption{Evolutions of energy and volume curves (from left to right)
for time-fractional Cahn-Hilliard equation \eqref{Problem-3}
 with fractional order $\alpha=0.4,\,0.7$ and $0.9$, respectively.}
\label{Comparison-CH-Energy-Mass}
\end{figure}
%%%%%%%%%%%%%%%%%%%%%%%%%%%%%%%%%%%%%%%%%%%%%%%%%%%%%%%%%%%%%%%%%%%%%%%%%%%%%%%%%%%%%

Now we simulate the merging of four drops
by the CN-SAV approach with the adaptive time-stepping strategy
using the mesh parameters $N_{0}=30$, $T_{0}=0.01$, $\kappa=10^{6}$, $\tau_{\min}=\tau_{N_{0}}=10^{-3}$
and $\tau_{\max}=10^{-1}$.
The solution snapshots are depicted in Figure \ref{Nonlocal-AC-Drops}
and the time evolution of the discrete energy and volume are plotted in Figure
\ref{Comparison-Nonlocal-AC-Energy-Mass}, respectively.
It is apparent that the initial separated four bubbles gradually coalesce into a single big
bubble and round up at the end of simulation corresponding to the minimization
of the interface area between two phases.
The coalescence speed of four bubbles are evidently affected by the fractional index $\alpha$,
that is, the larger the fractional order $\alpha$ is, the faster the coalescence.
Also, from Figure \ref{Comparison-Nonlocal-AC-Energy-Mass},
we see that the energy decreases in accord with the behavior of numerical solution,
and the volume is conserved just as predicted in Theorem \ref{Nonlocal-AC-CN-SAV-Volume}.
The solution of time-fractional Allen-Cahn equation \eqref{Problem-1} has a different behavior,
see Figures \ref{Classical-AC-Drops}-\ref{Comparison-Classical-AC-Energy-Mass},
the bubble shrinks and finally disappears because the equation \eqref{Problem-1} does not conserve the volume.
It is seen that the volume-preserving time-fractional Allen-Cahn equation \eqref{Problem-2}
may be a better choice for accurately simulating
the coalescence of bubbles than the non-volume-preserving version.

We also use the CN-SAV approach with the \emph{Adaptive Step}  approach
using $T_{0}=0.01$, $\kappa=10^{3}$ and $\tau_{\min}=\tau_{N_{0}}$
to simulate the time-fractional Cahn-Hilliard
equation \eqref{Problem-3}. The numerical results are given in Figures  \ref{CH-Drops}-\ref{Comparison-CH-Energy-Mass}.
For three different fractional order $\alpha=0.4,0.7,0.9$,
we take $N_{0}=300,100,30$, $\gamma=5,4,3$ and $\tau_{\max}=10^{-2},10^{-1},10^{-1}$, respectively.
Figure  \ref{CH-Drops} shows that the initial bubbles coalesce into one bubble quite rapidly
and the steady state is reached immediately. Furthermore, the energy falls off steeply
and decays faster for smaller fractional order$\alpha$,
see Figure \ref{Comparison-CH-Energy-Mass}, so that
 some extremely small time steps are required to capture this remarkable behavior.
Also, it is clear that the volume is conserved during the simulation.
Simple comparison from Figures  \ref{Comparison-Nonlocal-AC-Energy-Mass} and \ref{Comparison-CH-Energy-Mass}
shows that the energy of time-fractional Cahn-Hilliard
equation \eqref{Problem-3} dissipates much faster than that of time-fractional
Allen-Cahn equation \eqref{Problem-2} with the nonlocal volume constraint.
Correspondingly, the coalescence speed of initial bubbles of the former much faster
than that of the latter, see Figures \ref{Classical-AC-Drops} and \ref{CH-Drops}.

%%%%%%%%%%%%%%%%%%%%%%%%%%%%%%%%%%%%%%%%%%%%%%%%%%%%%%%%%%%%%%%%%%%%%%%%%%%%%%%%%%%%%
%%%%%%%%%%%%%%%%%%%%%%%%%%%%%%%%%%%%%%%%%%%%%%%%%%%%%%%%%%%%%%%%%%%%%%%%%%%%%%%%%%%%%

\subsection{Coarsening dynamics}

\begin{example}\label{Simulating-Coarsening-Dynamics}
We investigate the coarsening dynamics of the conservative time-fractional
phase field models with the model parameters
$\lambda=0.1$ and $\varepsilon=0.05$.
If not explicitly specified, we use $128 \times 128$ equal distanced meshes in space
to discretize the domain $\Omega=(0,2\pi)^{2}$.
Consider a randomly initial condition by assigning a random number varying from
$-0.001$ to $0.001$ at each grid points.
Taking the simulating parameters $\beta=4$ and $C_{0}=1$, we always apply
the present CN-SAV methods with the \emph{Adaptive Step}  approach using the following mesh parameters
$T_{0}=0.01$, $N_{0}=30$, $\kappa=10^{3}$, $\tau_{\min}=\tau_{N_{0}}=10^{-3}$ and $\tau_{\max}=10^{-1}$.
\end{example}

%%%%%%%%%%%%%%%%%%%%%%%%%%%%%%%%%%%%%%%%%%%%%%%%%%%%%%%%%%%%%%%%%%%%%%%%%%%%%%%%%%%%%%%%%%
\begin{figure}[htb!]
\centering
\includegraphics[width=1.47in]{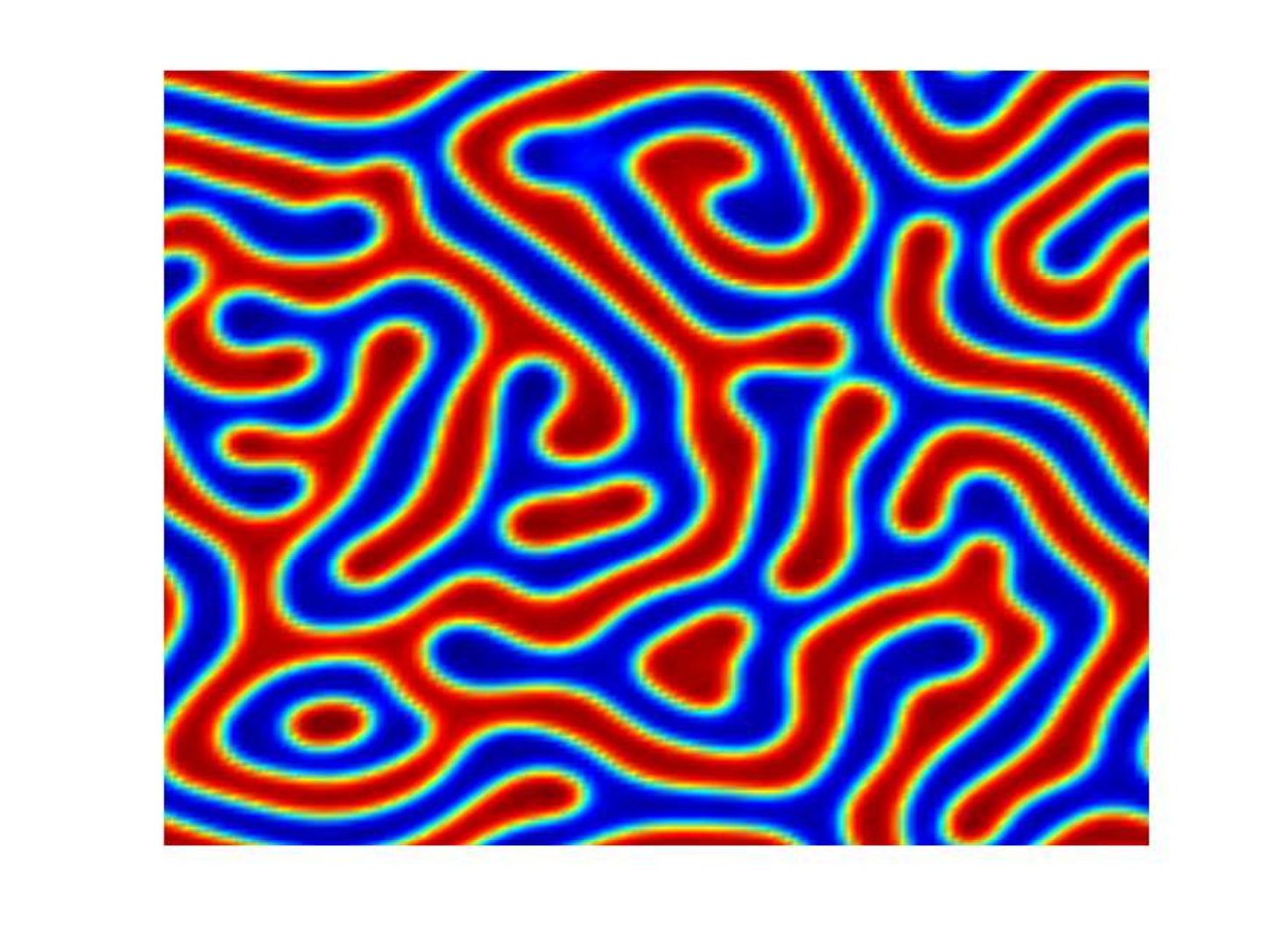}
\includegraphics[width=1.47in]{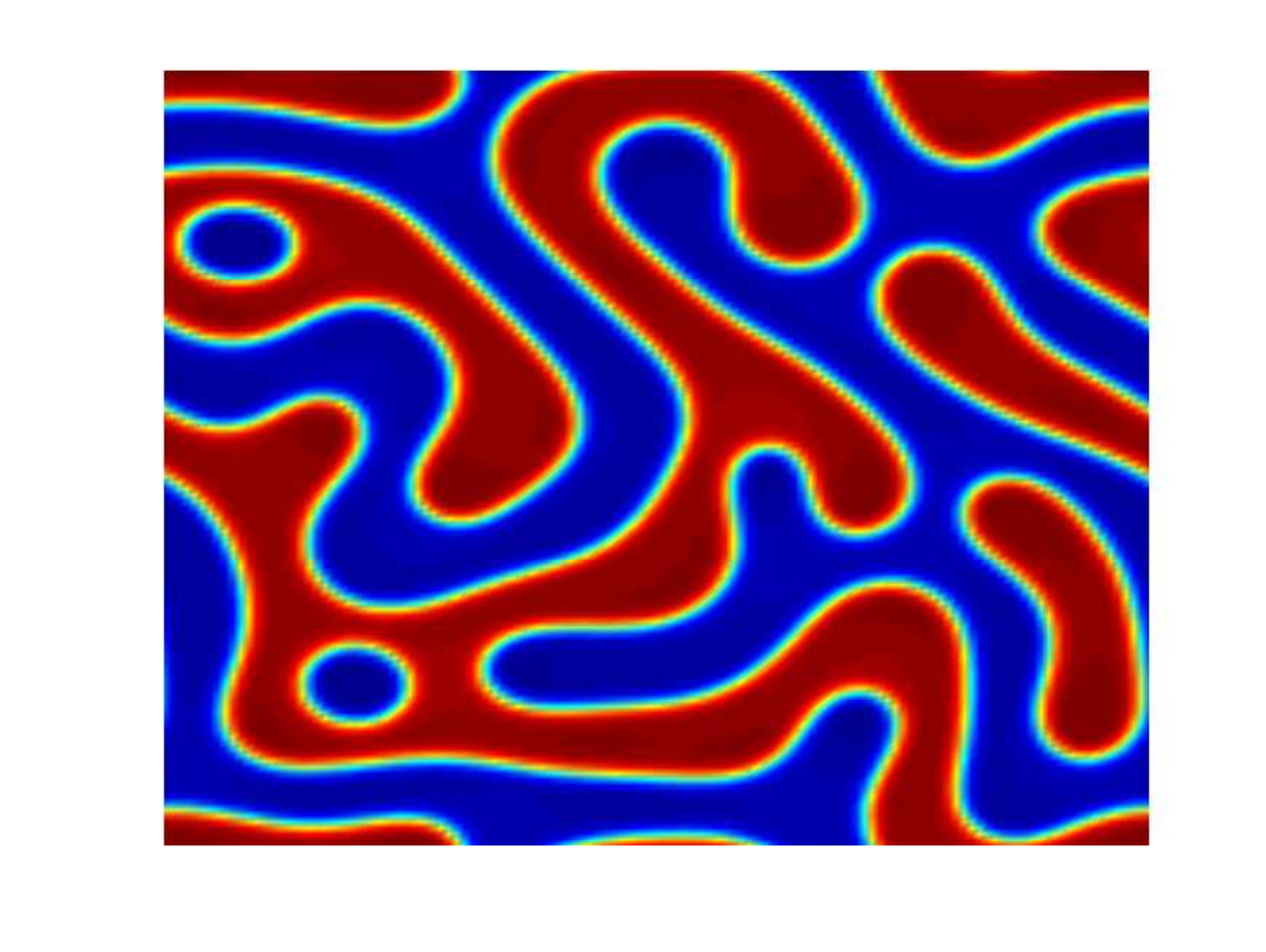}
\includegraphics[width=1.47in]{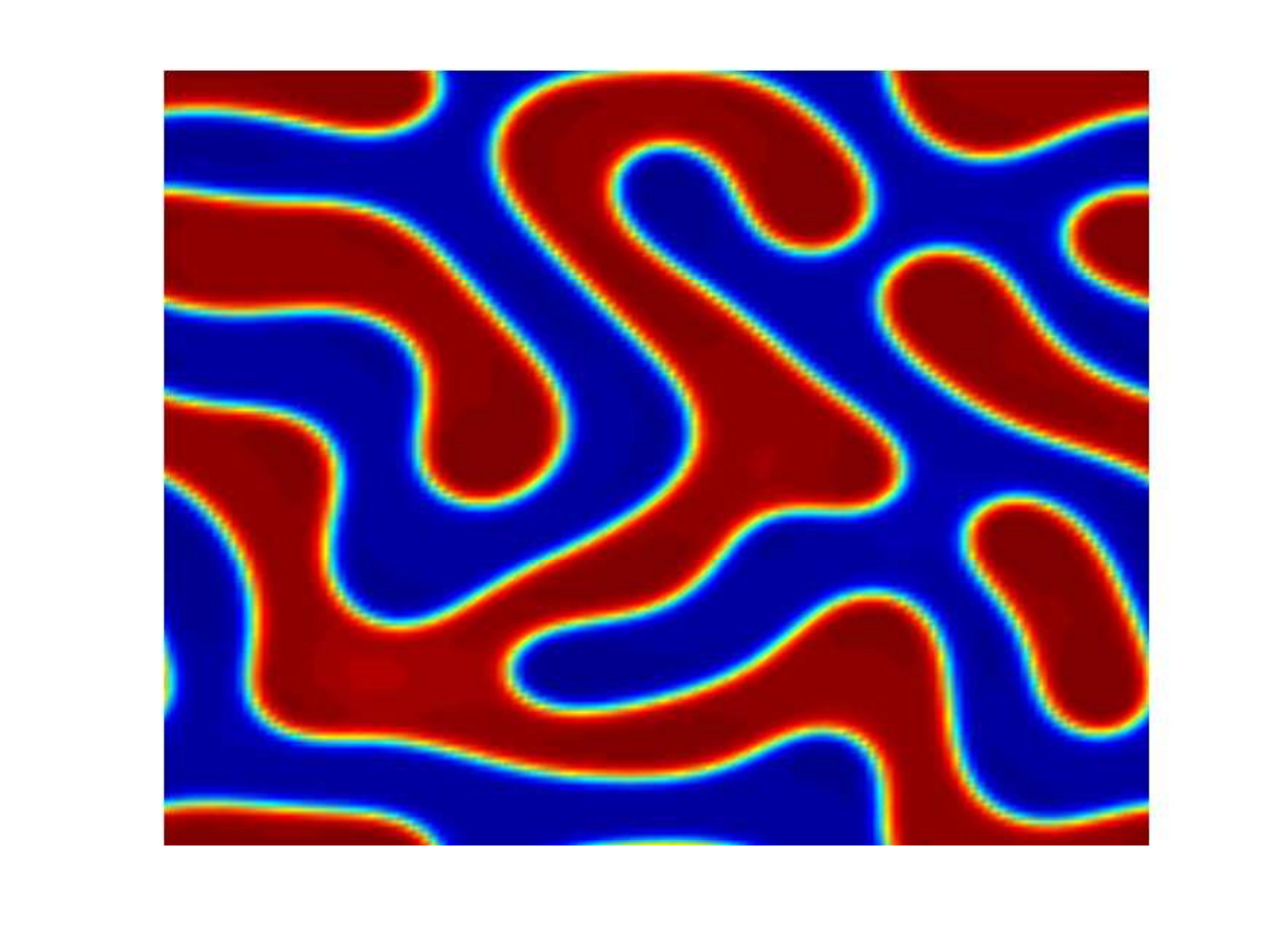}
\includegraphics[width=1.47in]{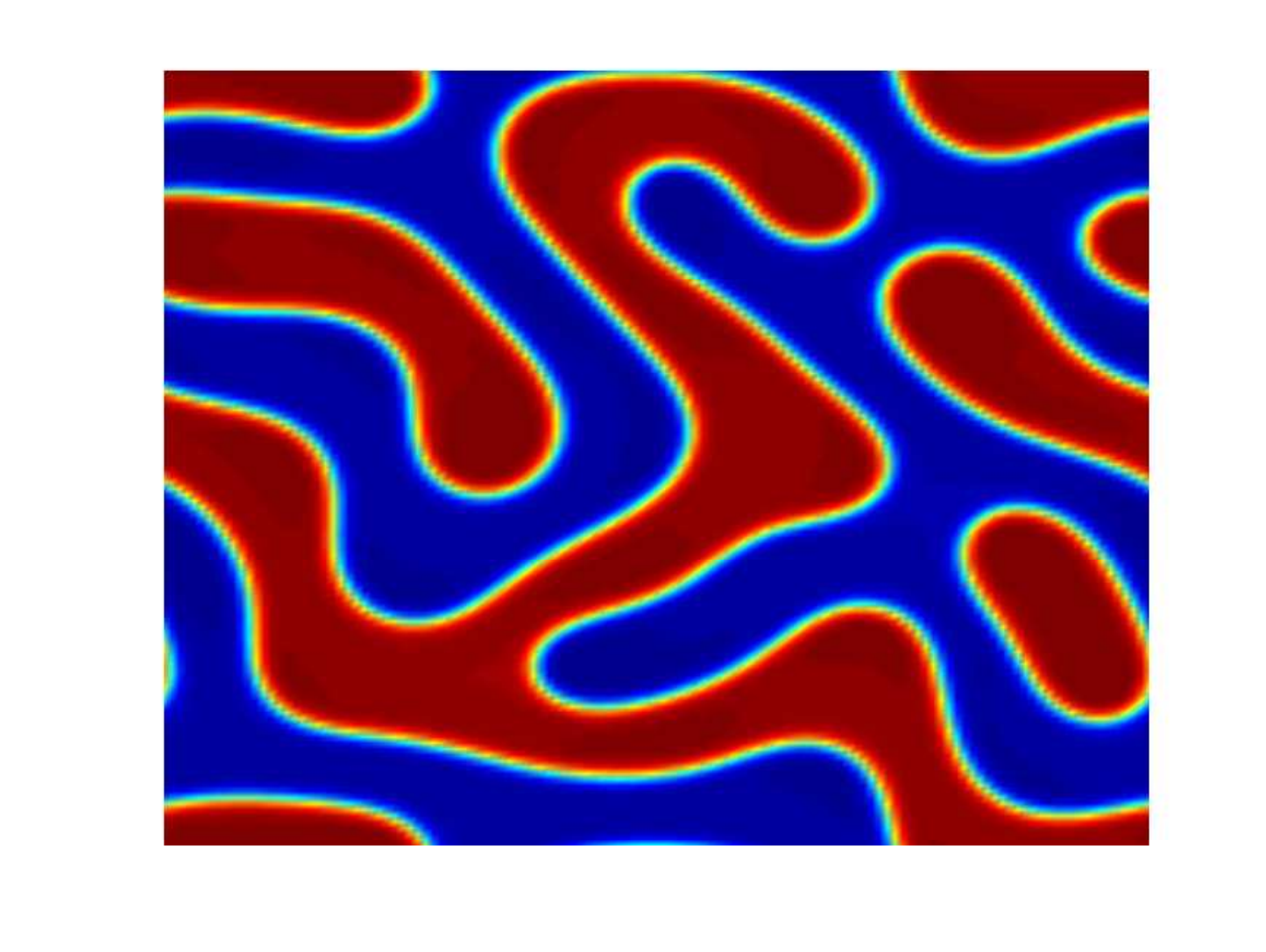}\\
\includegraphics[width=1.47in]{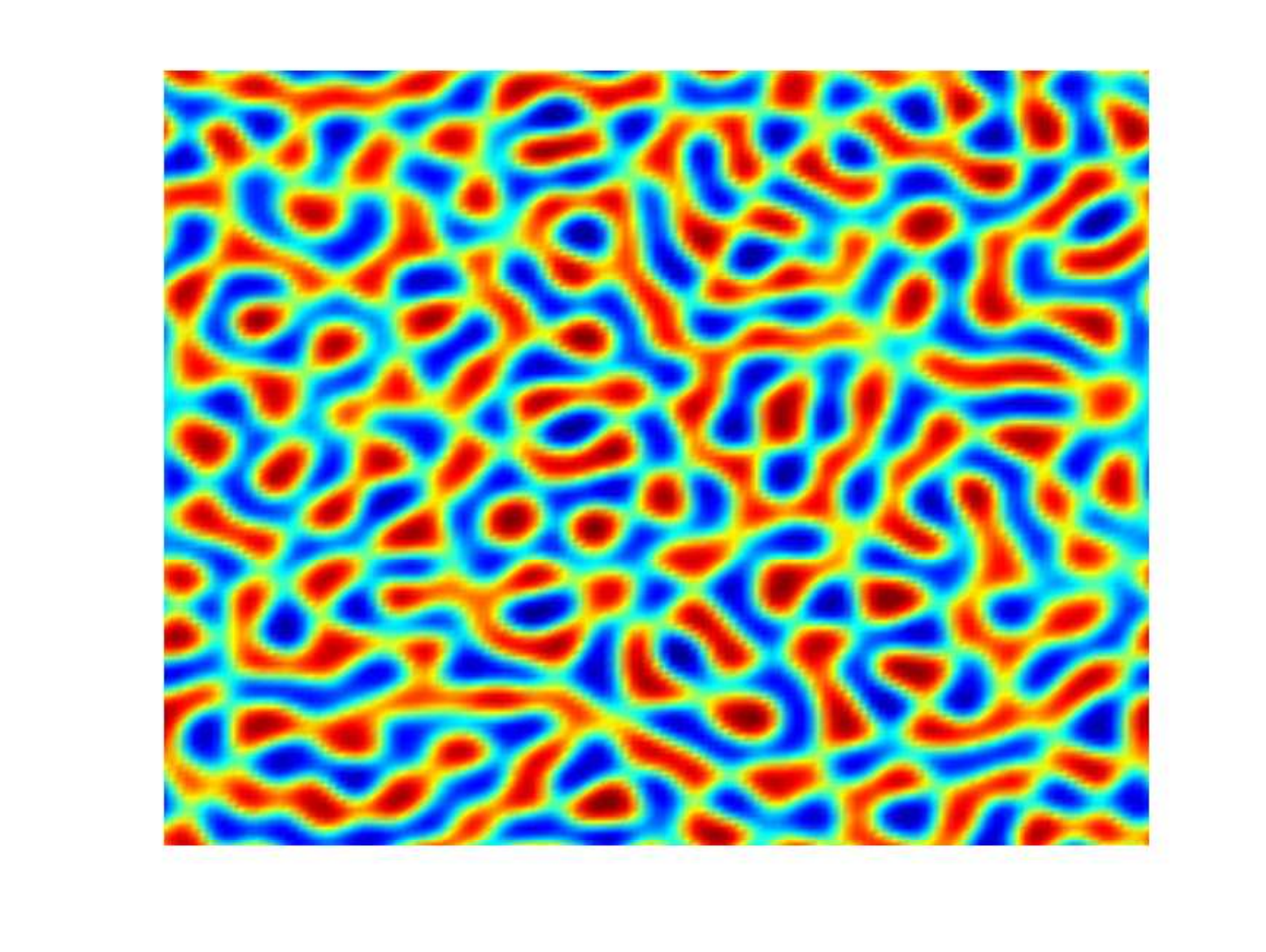}
\includegraphics[width=1.47in]{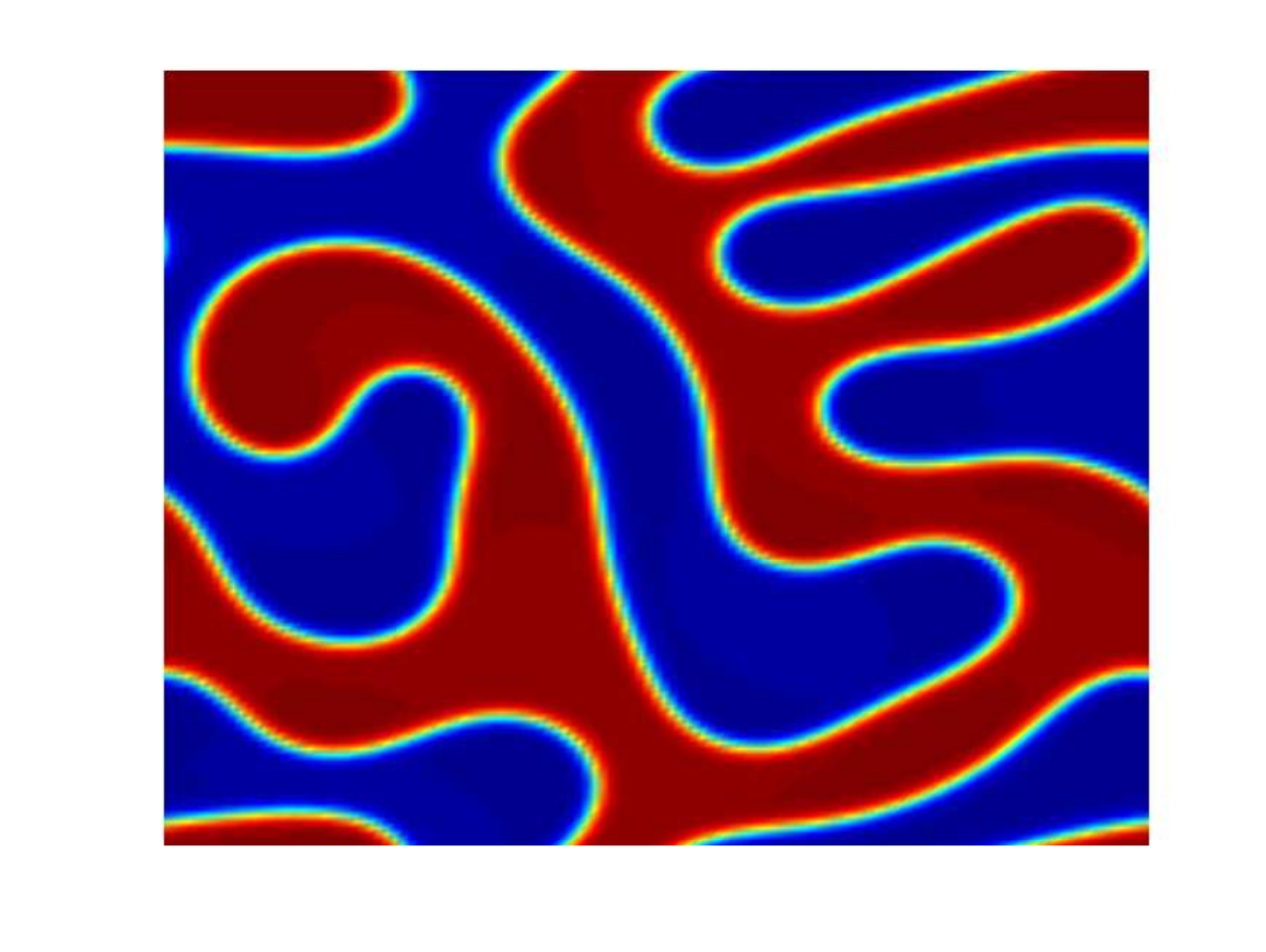}
\includegraphics[width=1.47in]{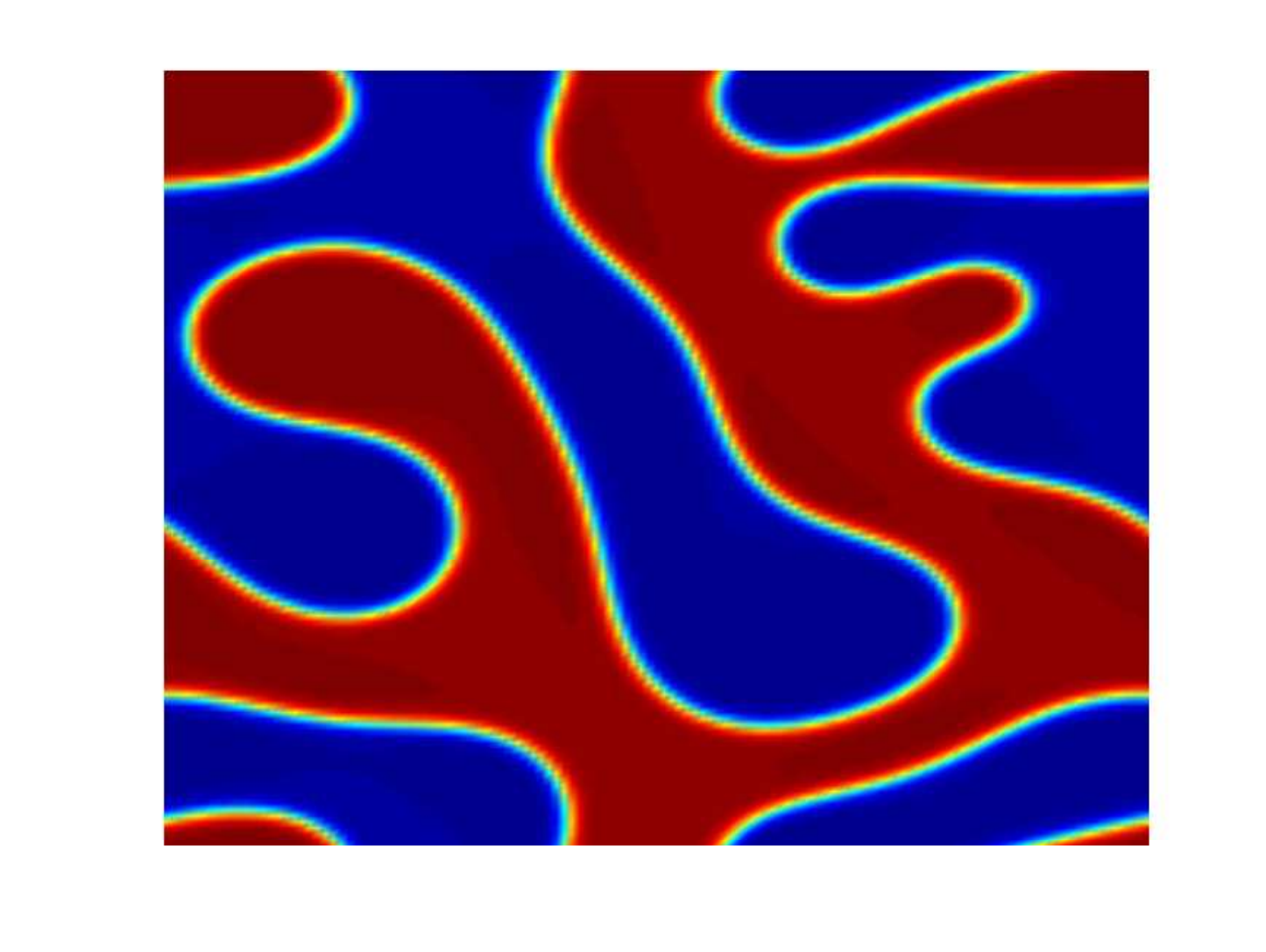}
\includegraphics[width=1.47in]{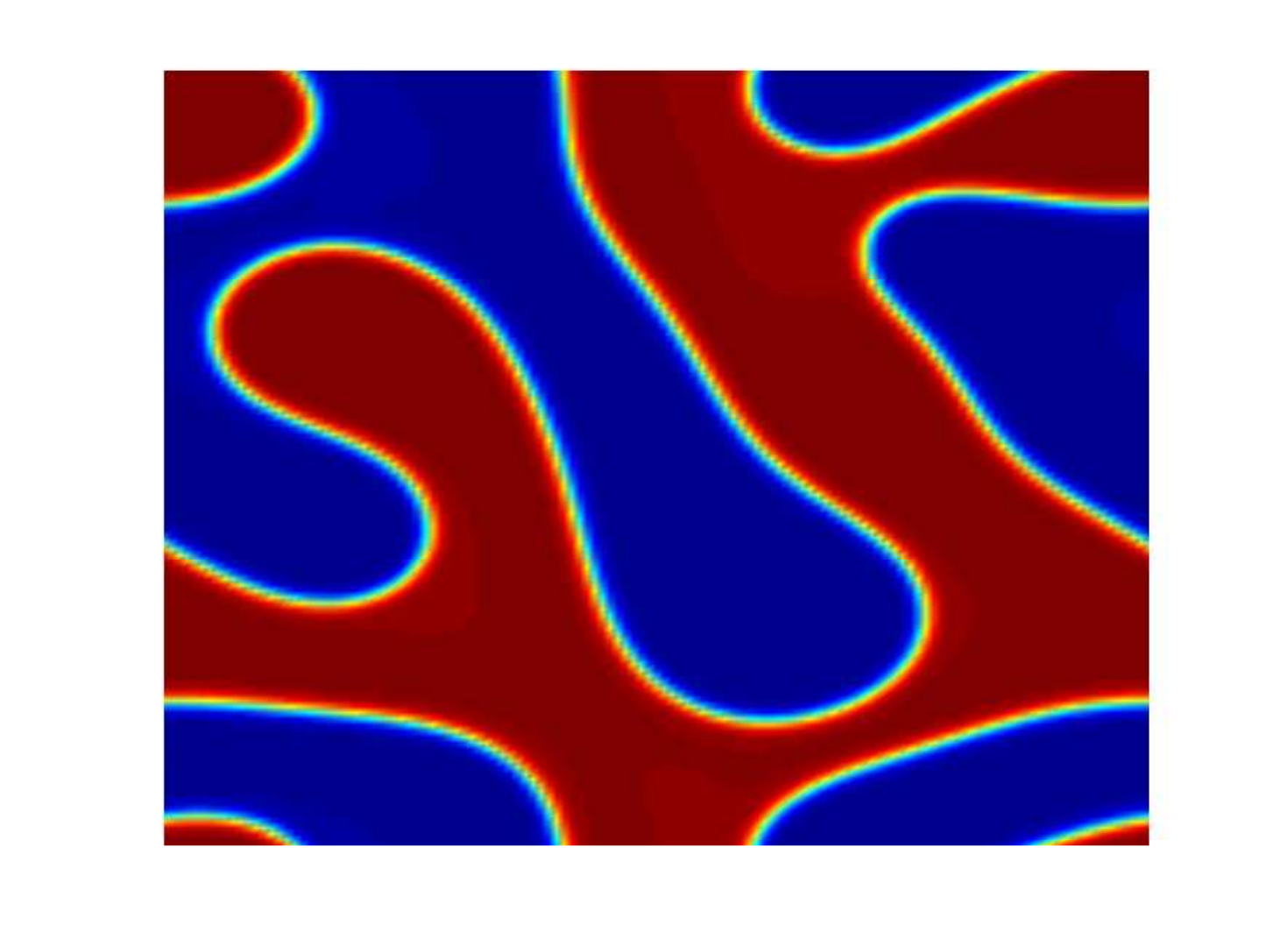}\\
\includegraphics[width=1.47in]{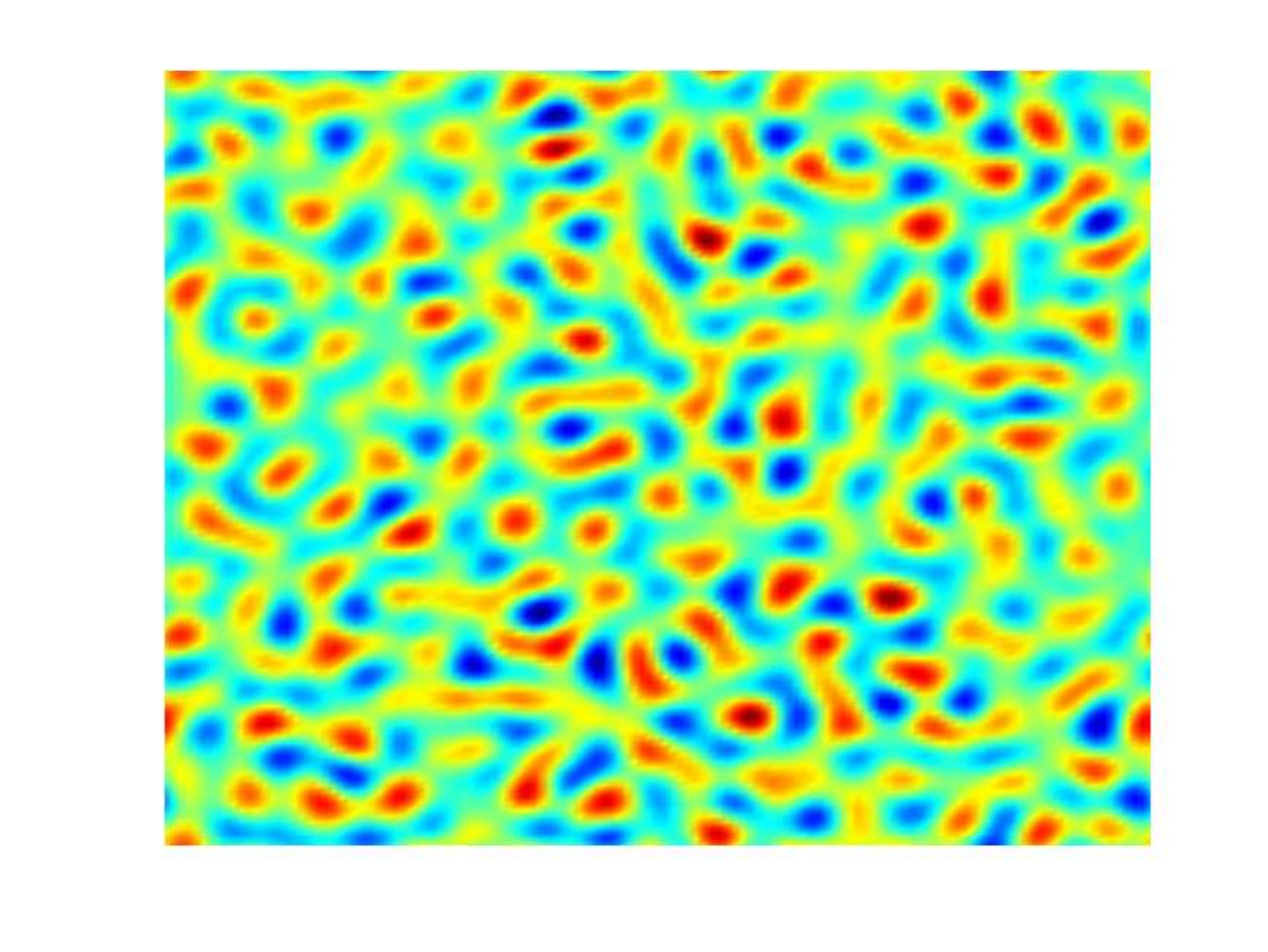}
\includegraphics[width=1.47in]{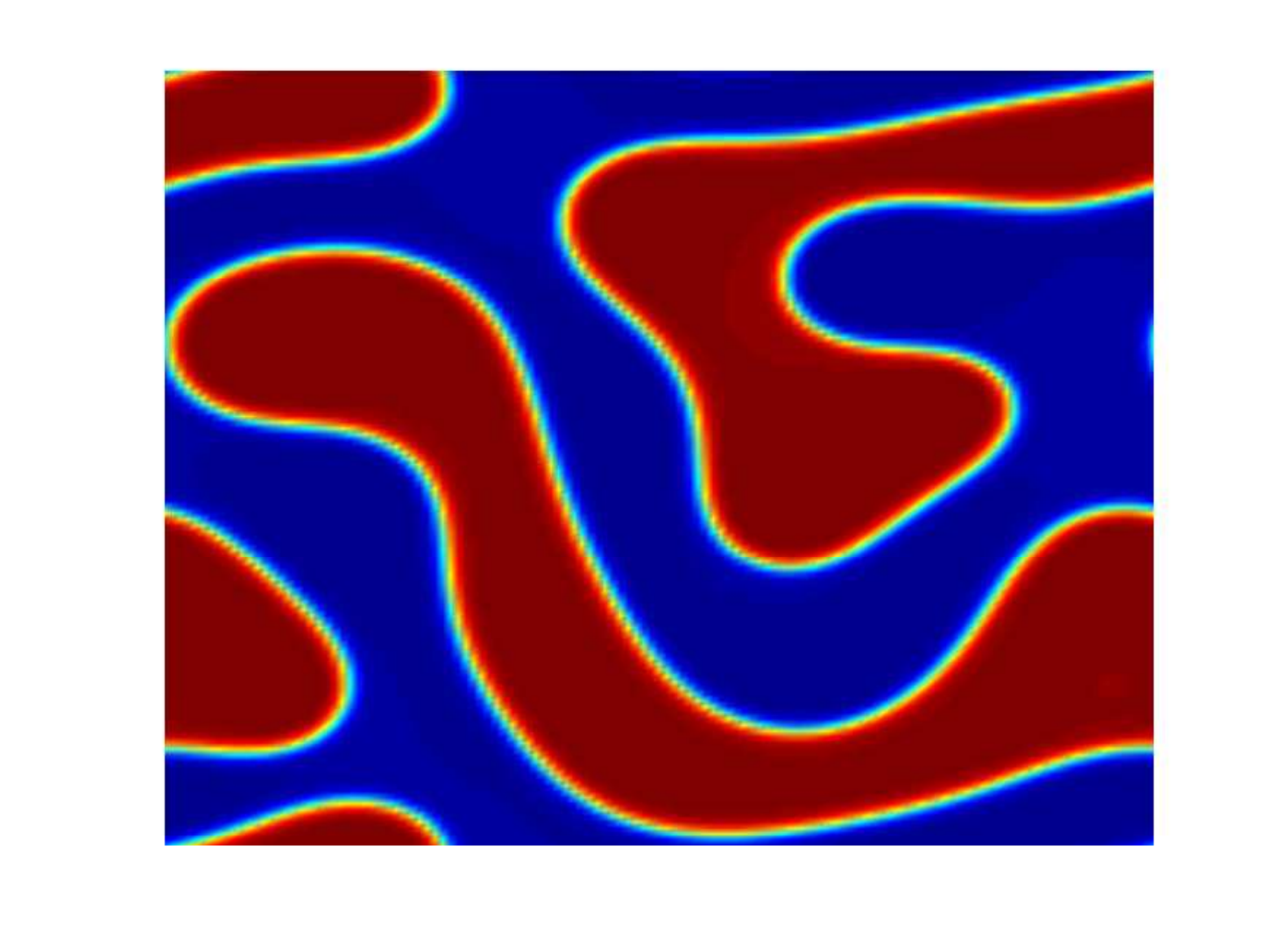}
\includegraphics[width=1.47in]{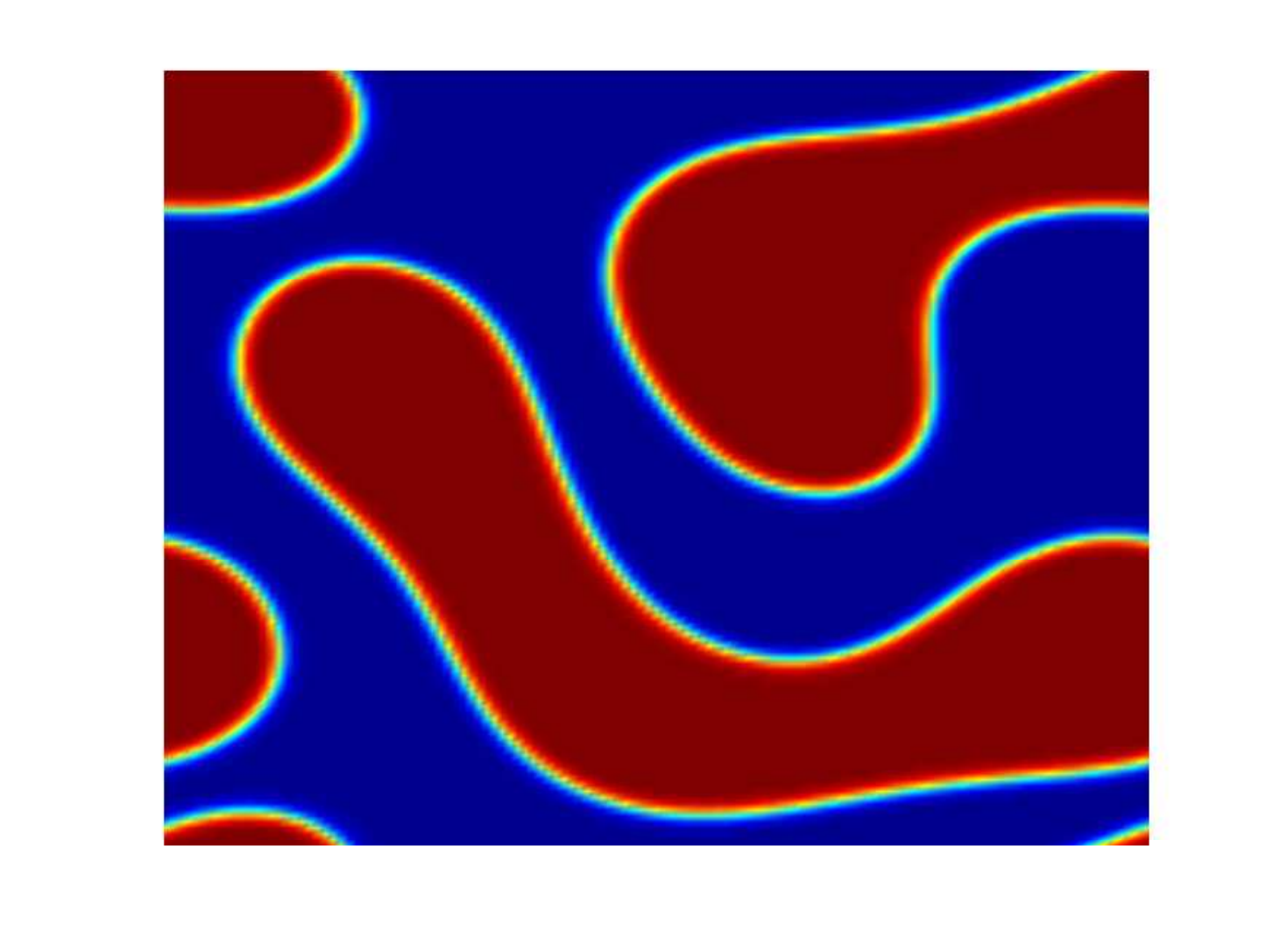}
\includegraphics[width=1.47in]{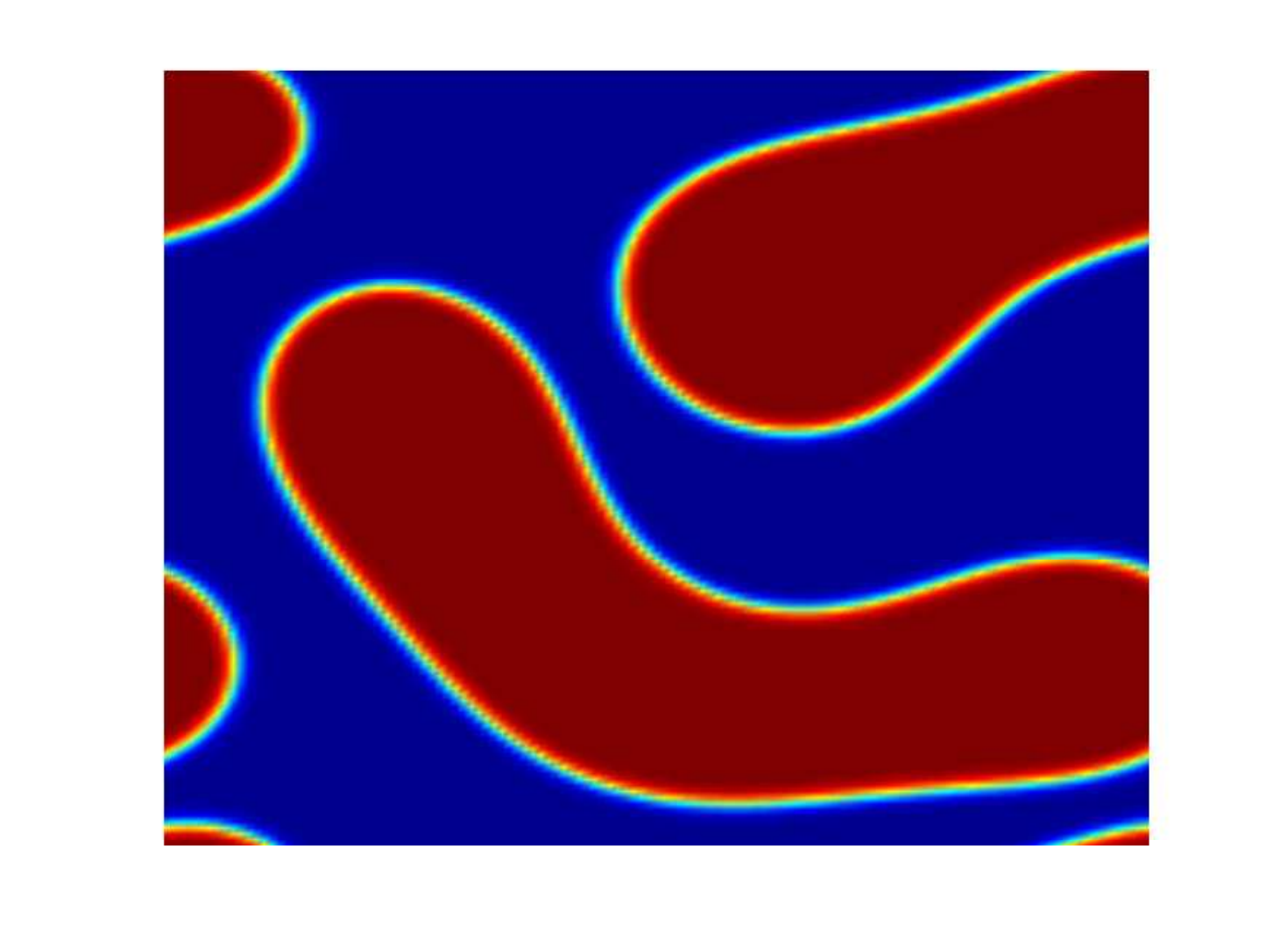}\\
\caption{Solution snapshots of coarsening dynamics
   of \eqref{Problem-3} at $t=1, 100, 300, 500$ (from left to right)
 with fractional orders $\alpha=0.4,\,0.7$ and $0.9$ (from top to bottom), respectively.}
\label{CH-Coarsening-Dynamic}
\end{figure}
%%%%%%%%%%%%%%%%%%%%%%%%%%%%%%%%%%%%%%%%%%%%%%%%%%%%%%%%%%%%%%%%%%%%%%%%%%%%%%%%%%%%%

%%%%%%%%%%%%%%%%%%%%%%%%%%%%%%%%%%%%%%%%%%%%%%%%%%%%%%%%%%%%%%%%%%%%%%%%%%%%%%%%%%%%%%%%%%%%%
\begin{figure}[htb!]
\centering
\includegraphics[width=2.0in]{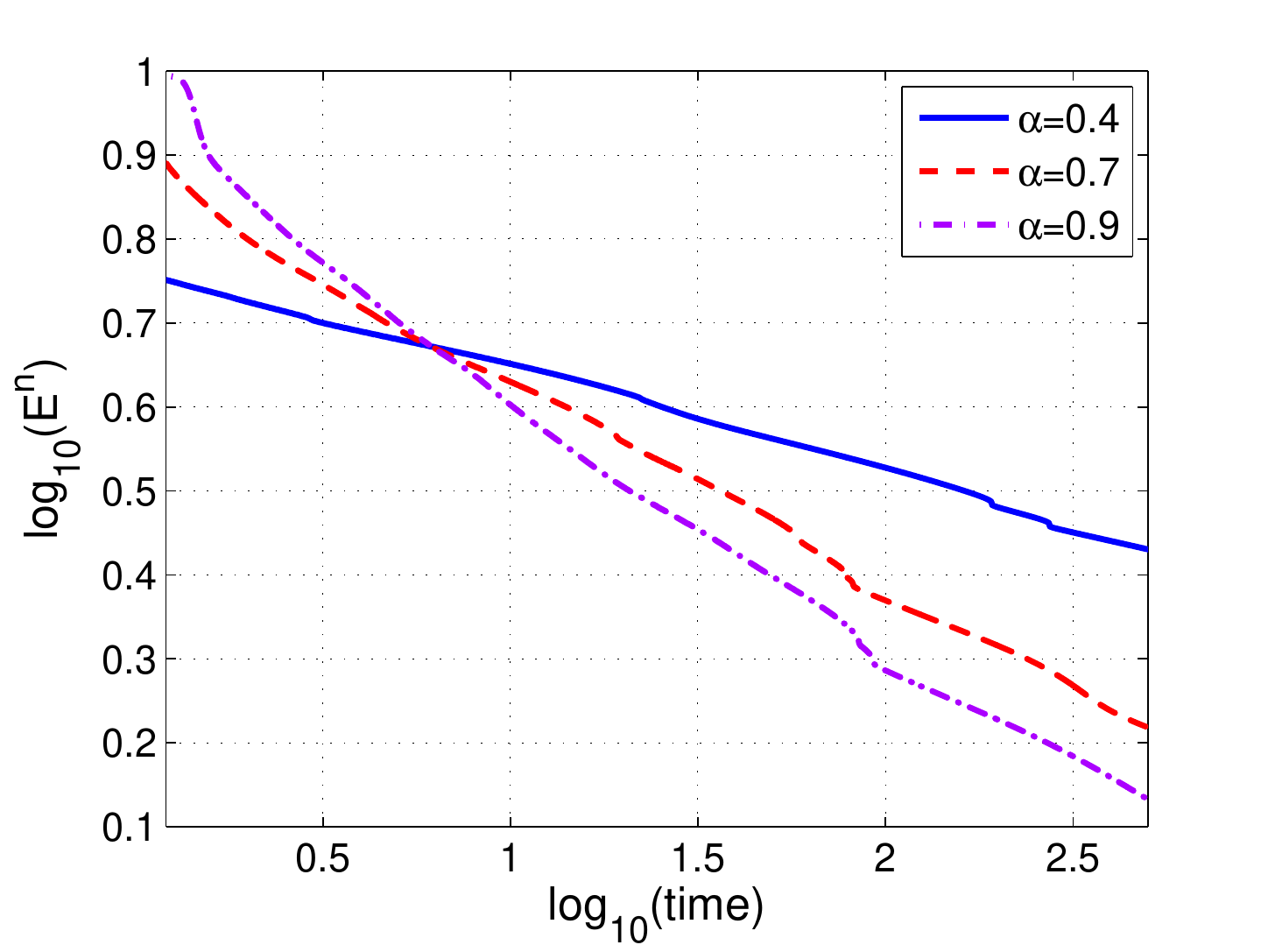}
\includegraphics[width=2.0in]{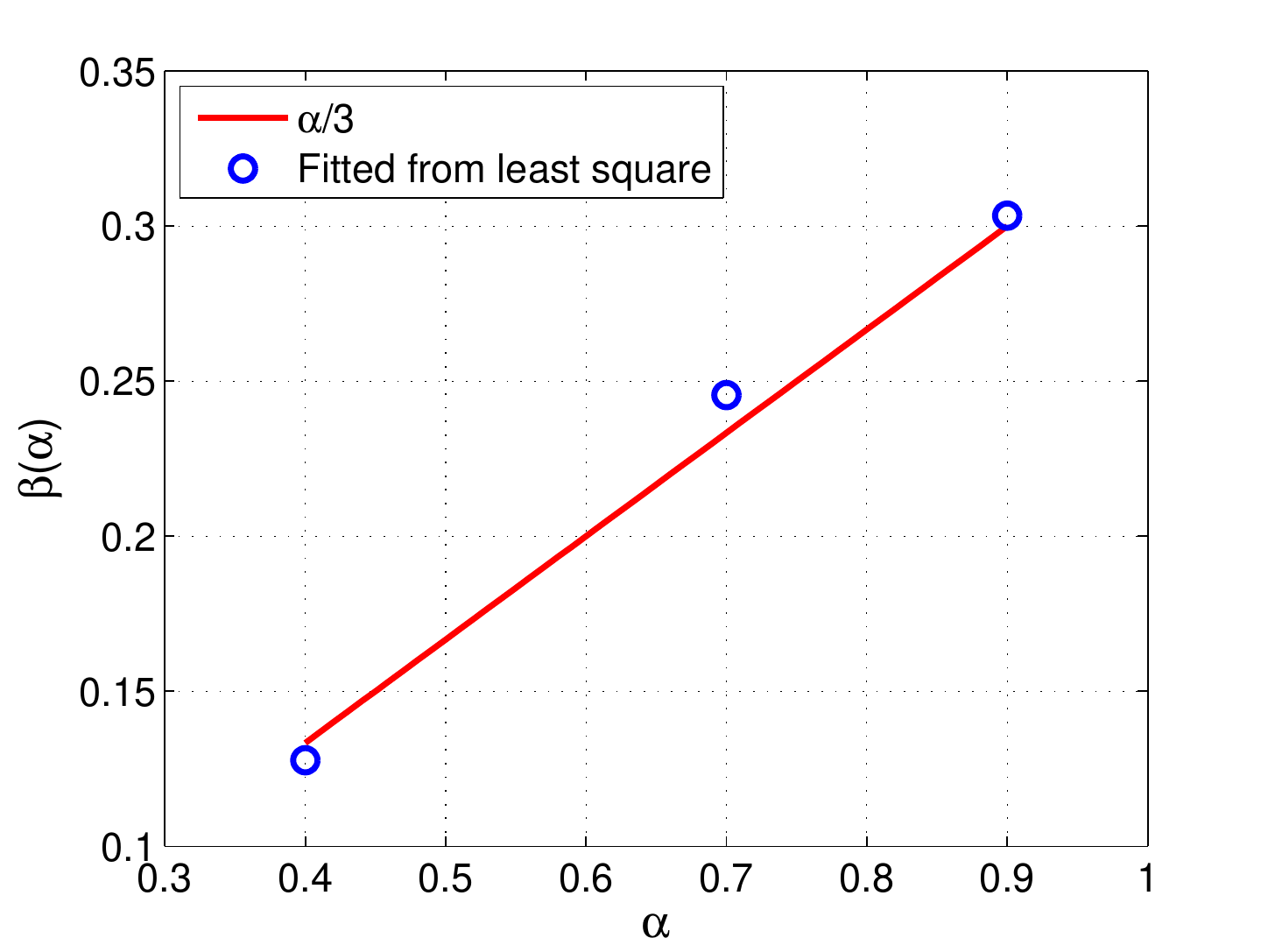}
\includegraphics[width=2.0in]{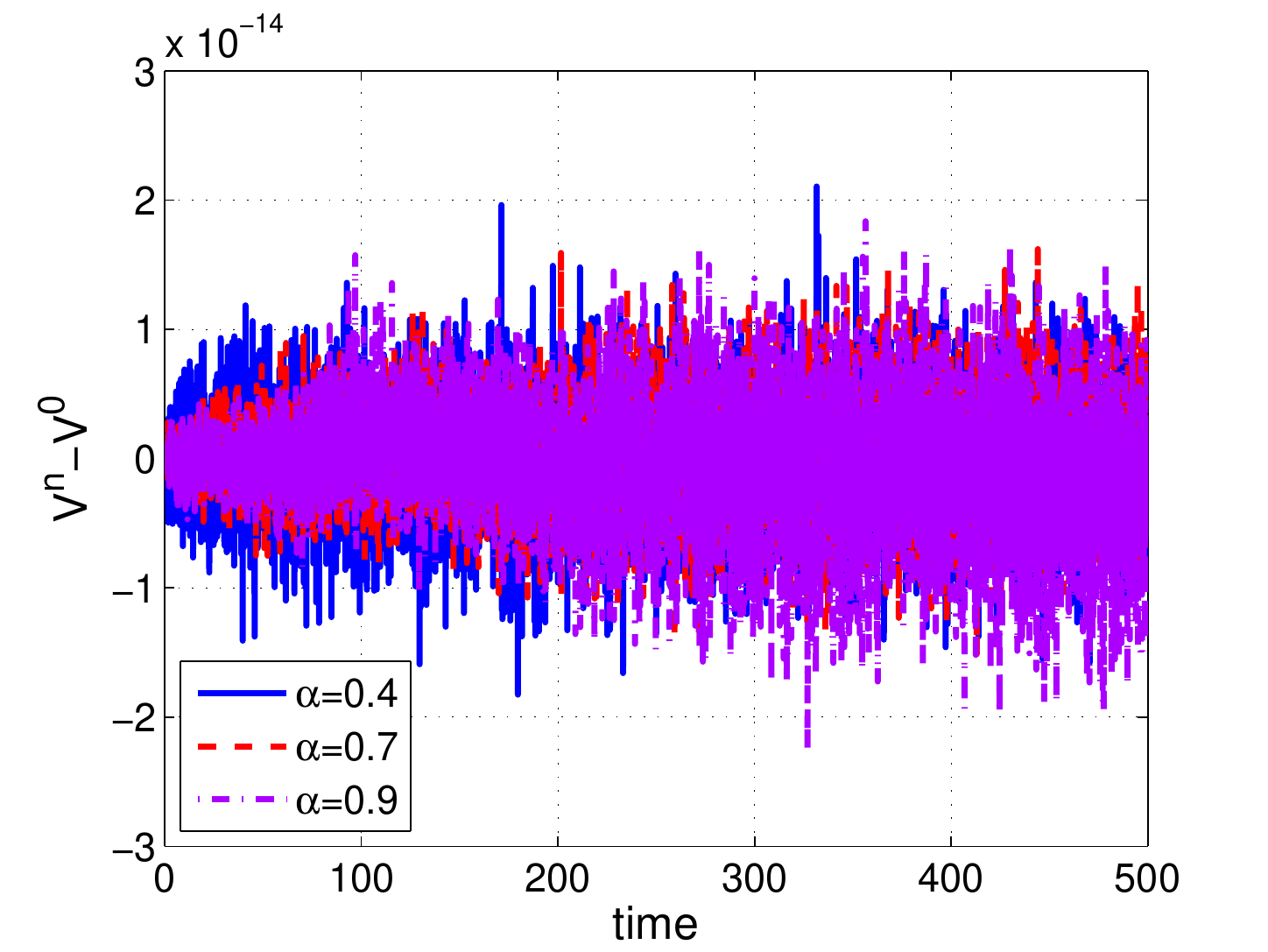}
\caption{Evolutions of energy, the least square fitted energy dissipation law scaling $\beta(\alpha)$ and volume conservation (from left to right)
of the  time-fractional Cahn-Hilliard equation for three fractional orders $\alpha=0.4,\,0.7$ and $0.9$, respectively.}
\label{CH-Dynamic-Energy-Mass}
\end{figure}
%%%%%%%%%%%%%%%%%%%%%%%%%%%%%%%%%%%%%%%%%%%%%%%%%%%%%%%%%%%%%%%%%%%%%%%%%%%%%%%%%%%%%

The snapshots of the coarsening dynamics of time-fractional Cahn-Hilliard equation \eqref{Problem-3}
with a variety of fractional order $\alpha$ at different time slots are depicted
in  Figure \ref{CH-Coarsening-Dynamic}. From the first column of Figure \ref{CH-Coarsening-Dynamic},
we find that the coarsening dynamics appear to be faster at the early time for smaller fractional order $\alpha$,
while it would be much slower as the time escapes.
In other words, the time-fractional Cahn-Hilliard model with larger fractional order $\alpha$
has faster evolution dynamics,
which is in good agreement with what we have observed in Example \ref{Simulating-Four-Drops}.
In Figure \ref{CH-Dynamic-Energy-Mass}, the energy dissipation law scaling $\beta(\alpha)$ is estimated by
doing the least square fit via the formula
$\log_{10}(E(\alpha,t))=\beta^{0}(\alpha)-\beta(\alpha)\log_{10}(t)$.
It is observed that the energy dissipates approximately as $O(t^{\frac{\alpha}{3}})$,
which is consistent with $O(t^{\frac{1}{3}})$ as $\alpha\rightarrow1$, as well-known.
As expected, the volume is also conserved during the coarsening process.

%%%%%%%%%%%%%%%%%%%%%%%%%%%%%%%%%%%%%%%%%%%%%%%%%%%%%%%%%%%%%%%%%%%%%%%%%%%%%%%%%%%%%%%%%%
\begin{figure}[htb!]
\centering
\includegraphics[width=1.47in]{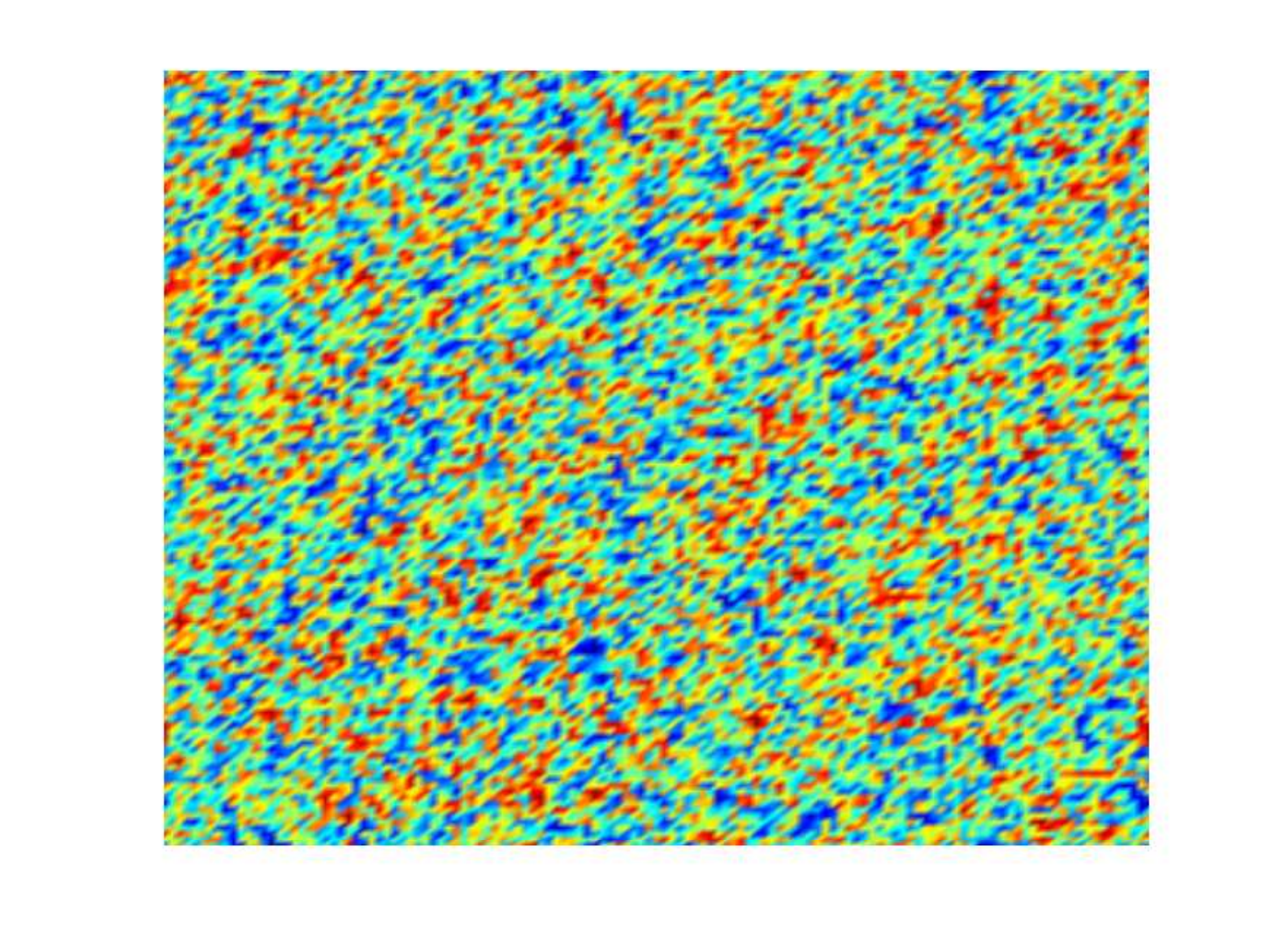}
\includegraphics[width=1.47in]{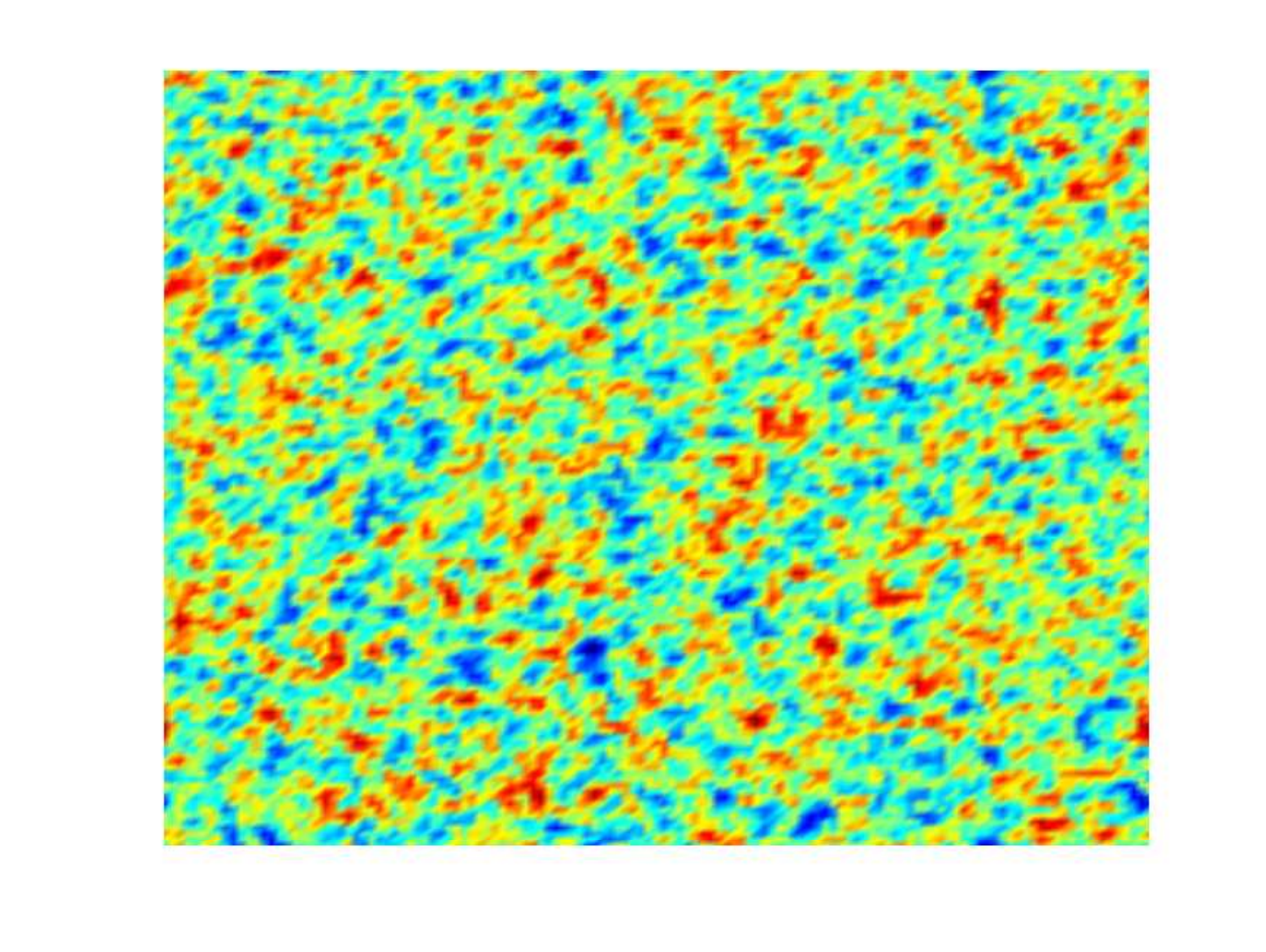}
\includegraphics[width=1.47in]{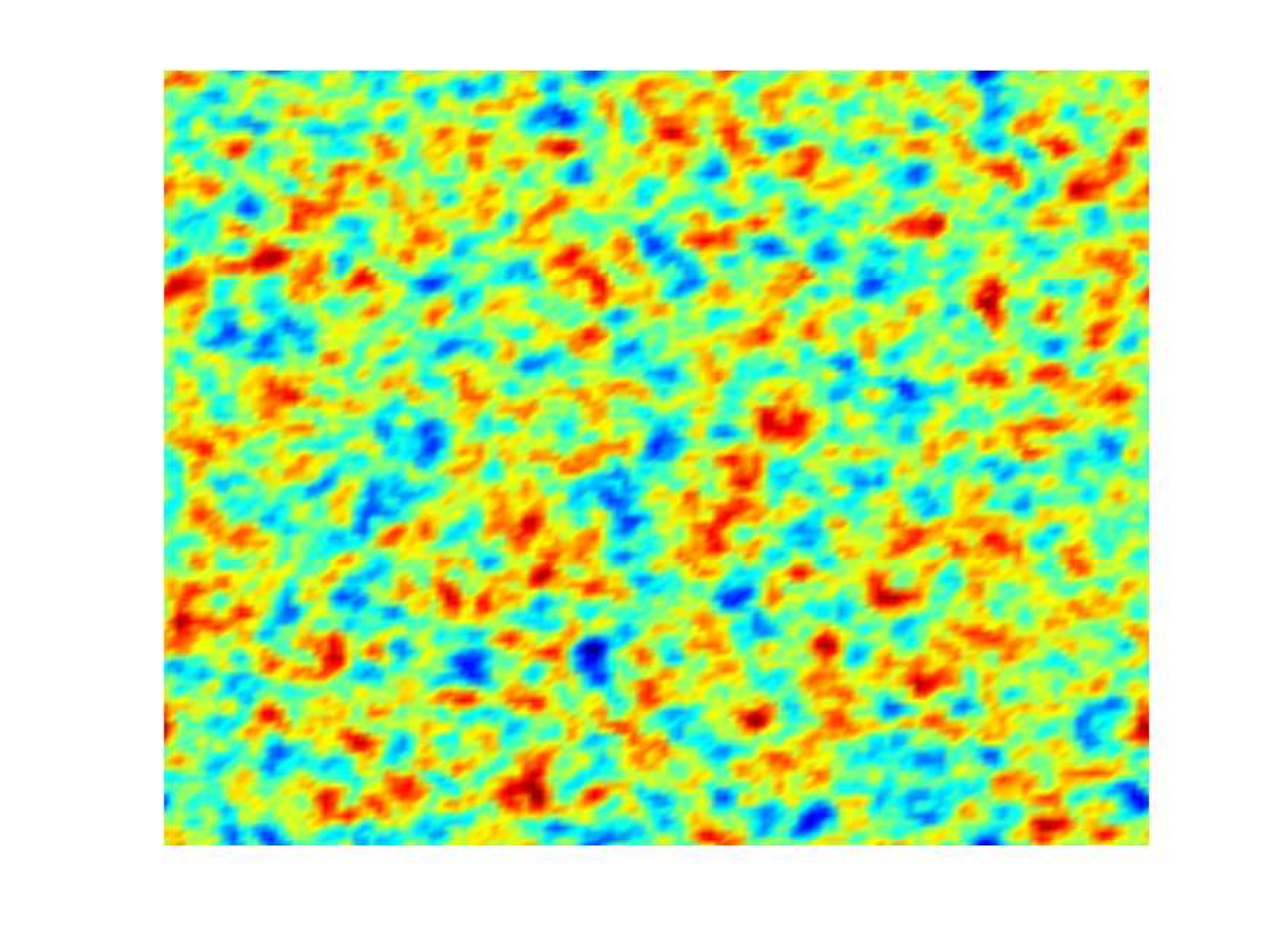}
\includegraphics[width=1.47in]{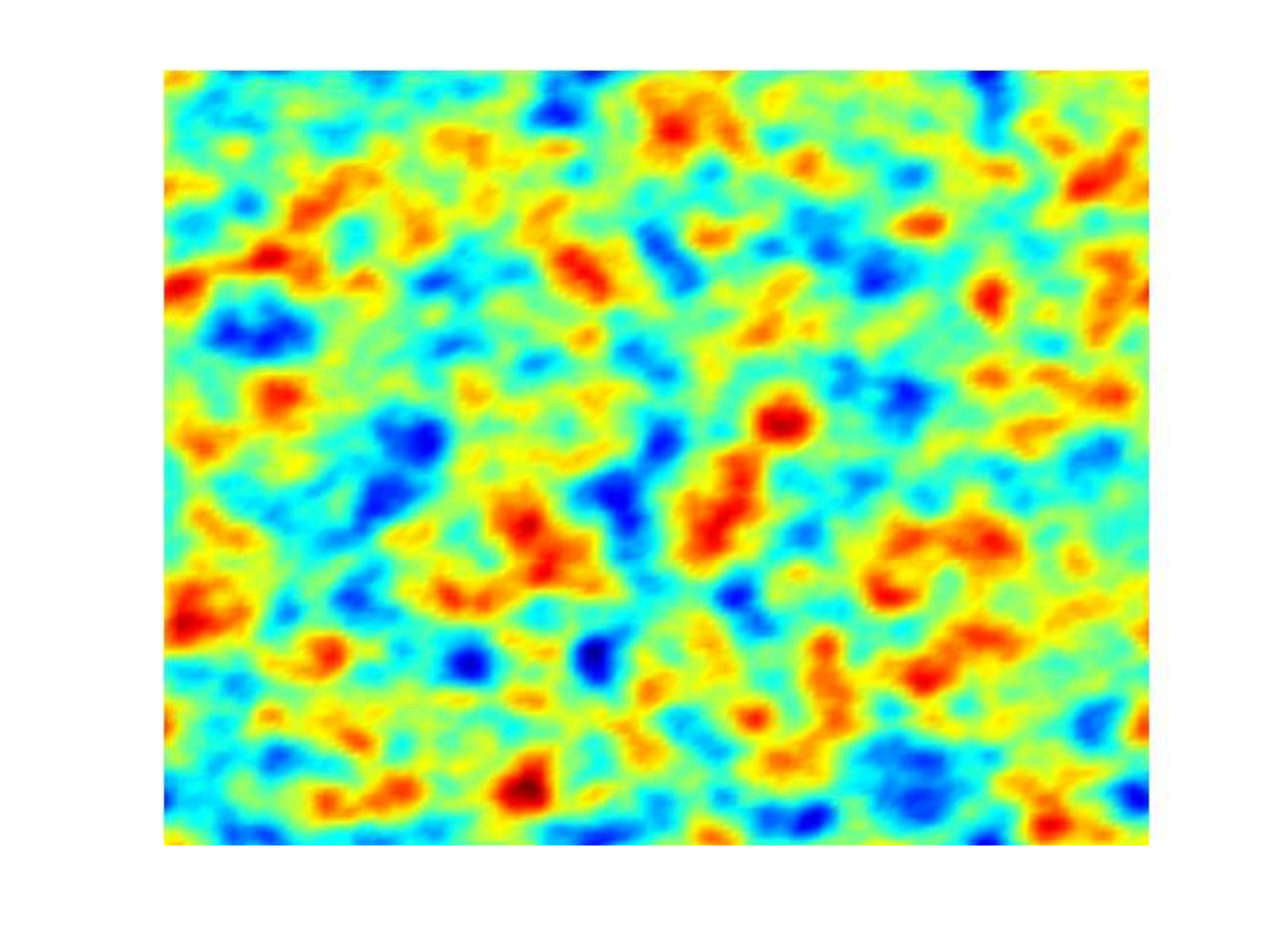}\\
\includegraphics[width=1.47in]{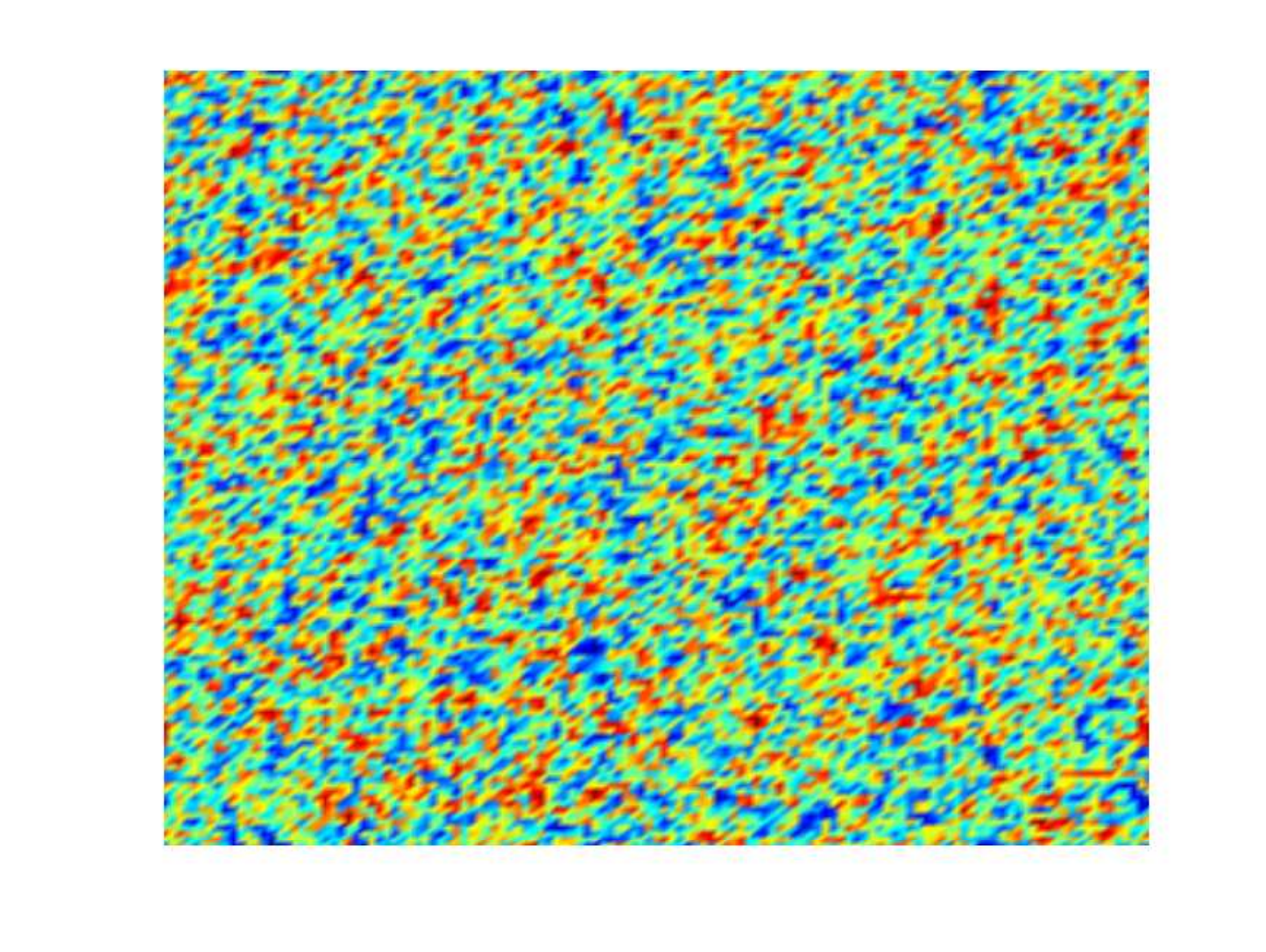}
\includegraphics[width=1.47in]{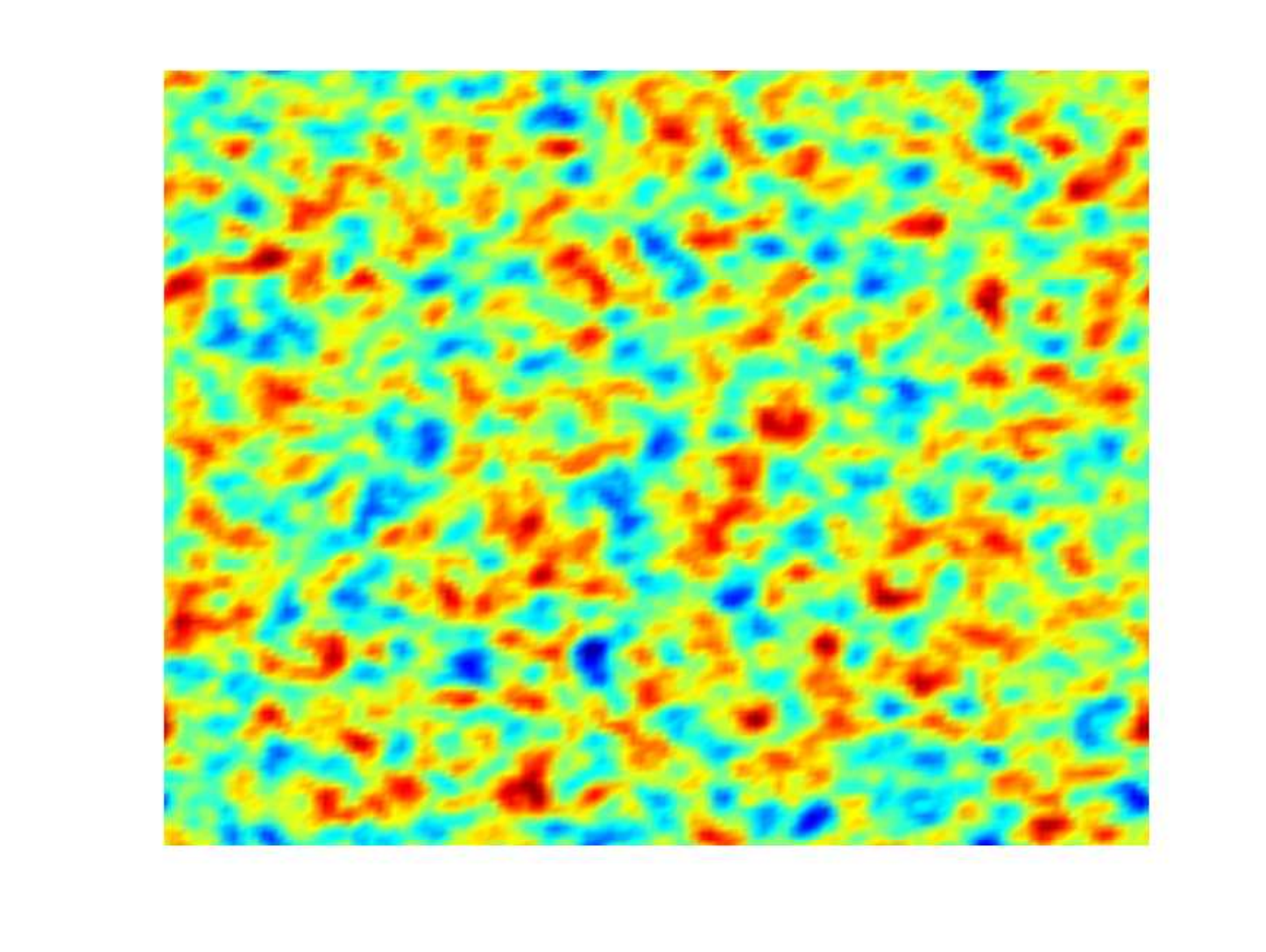}
\includegraphics[width=1.47in]{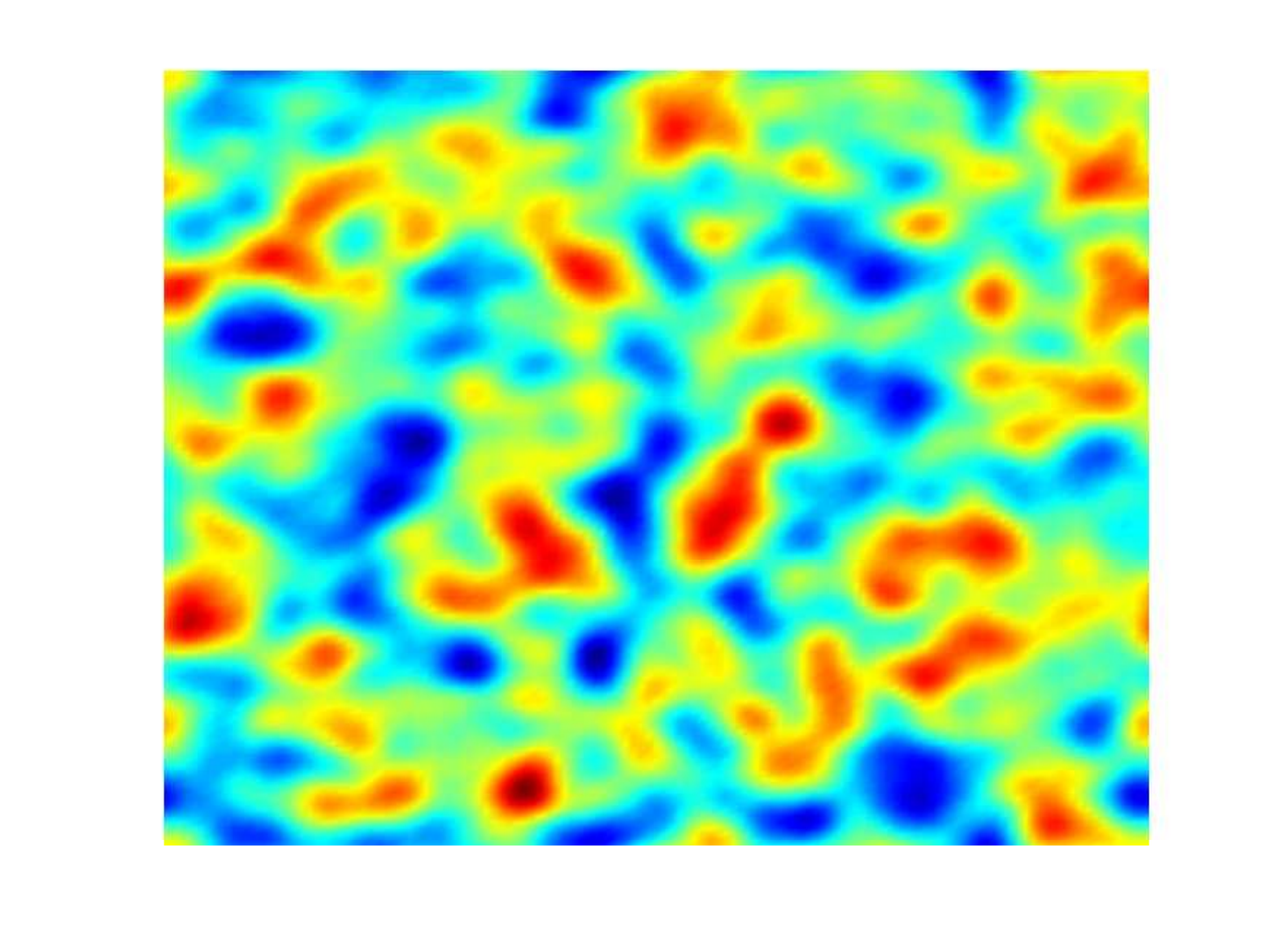}
\includegraphics[width=1.47in]{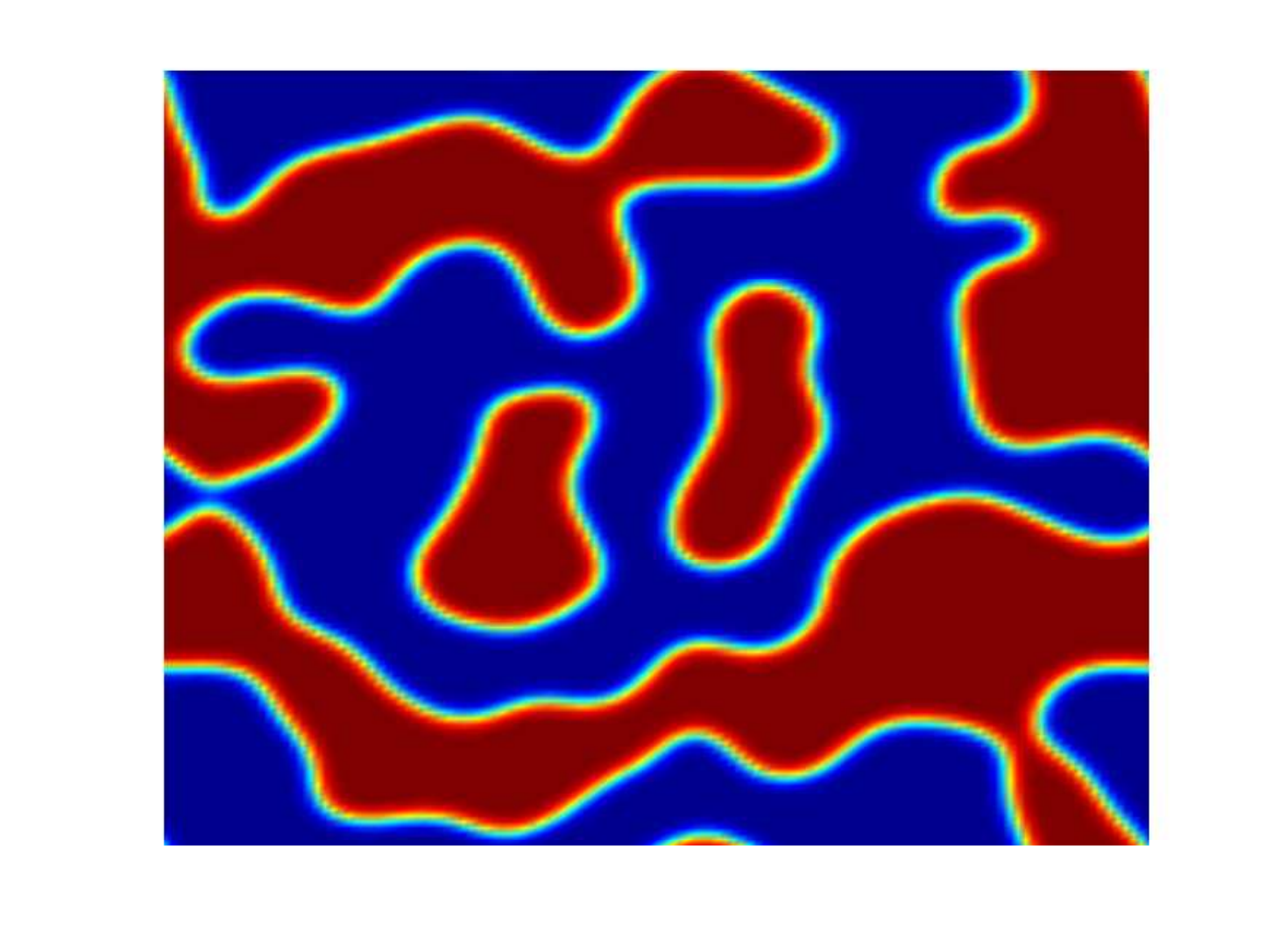}\\
\includegraphics[width=1.47in]{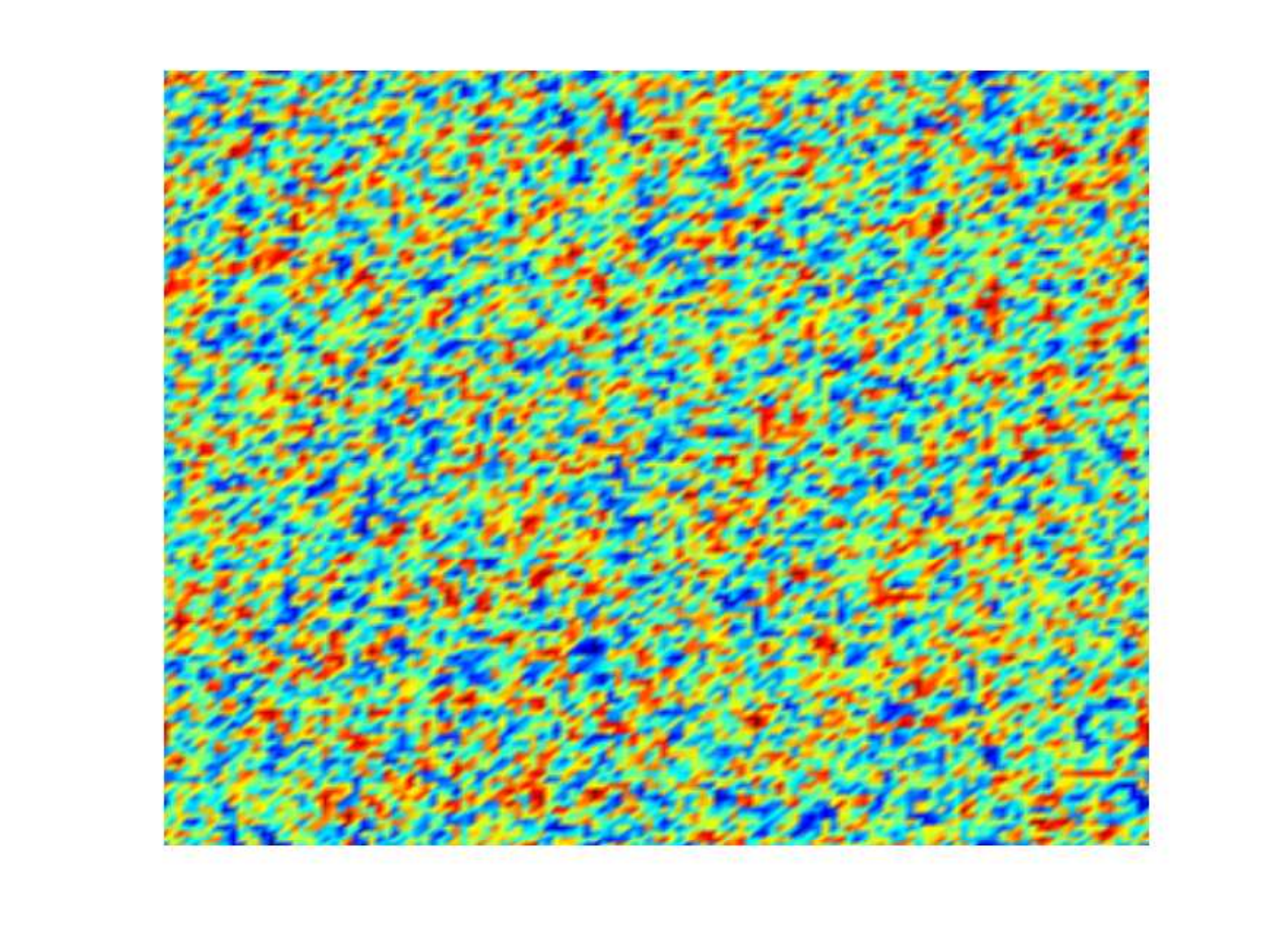}
\includegraphics[width=1.47in]{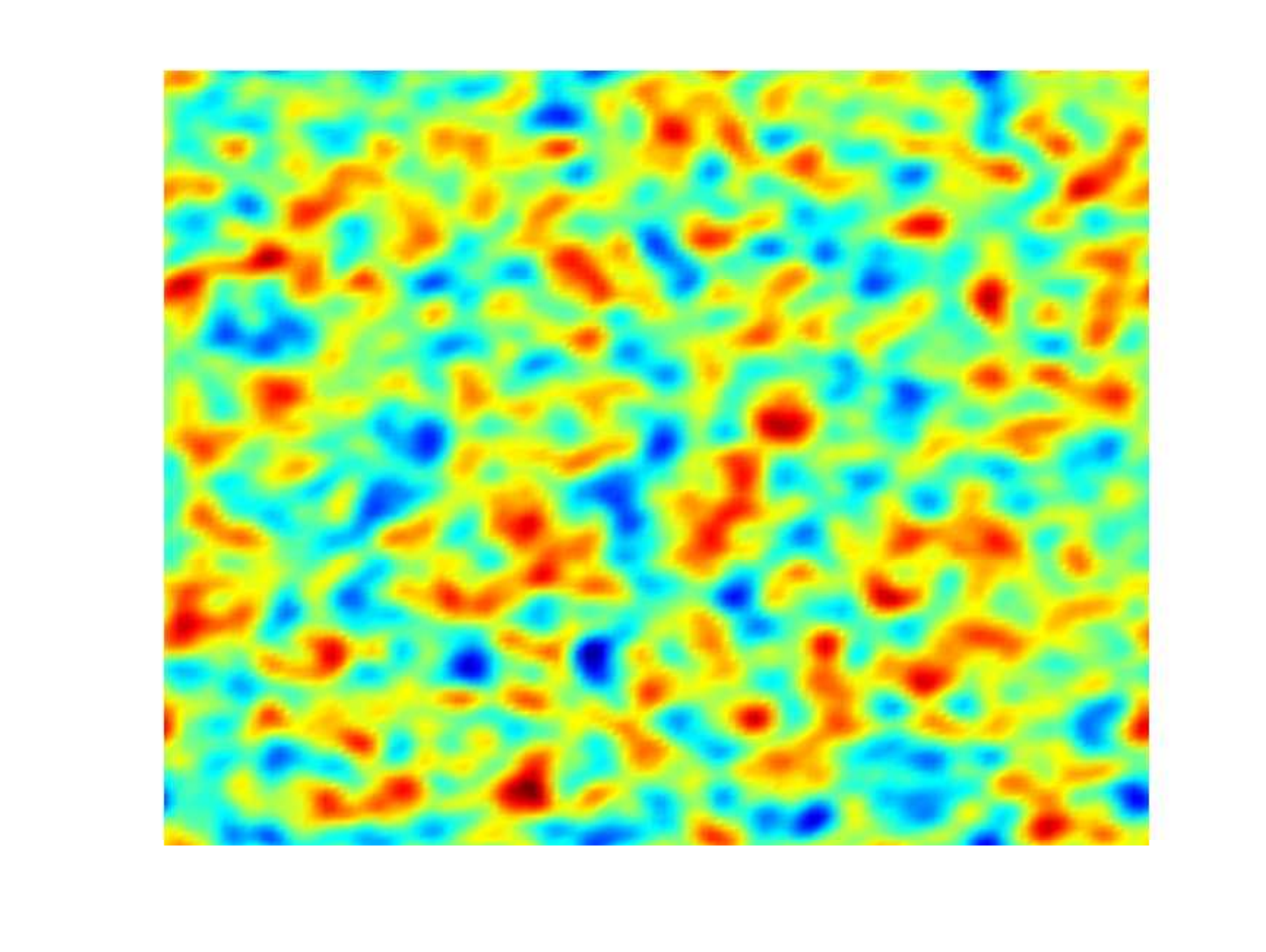}
\includegraphics[width=1.47in]{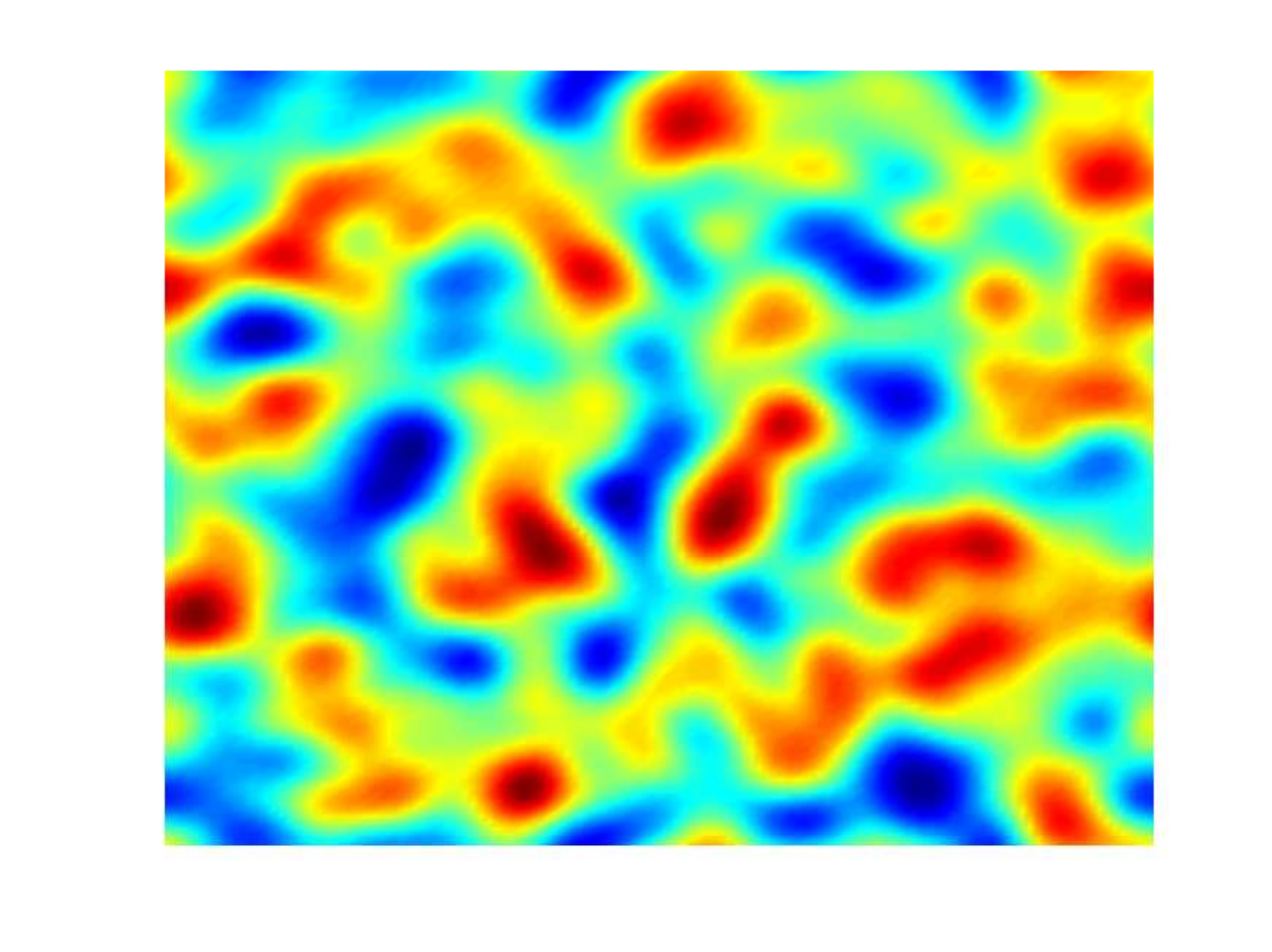}
\includegraphics[width=1.47in]{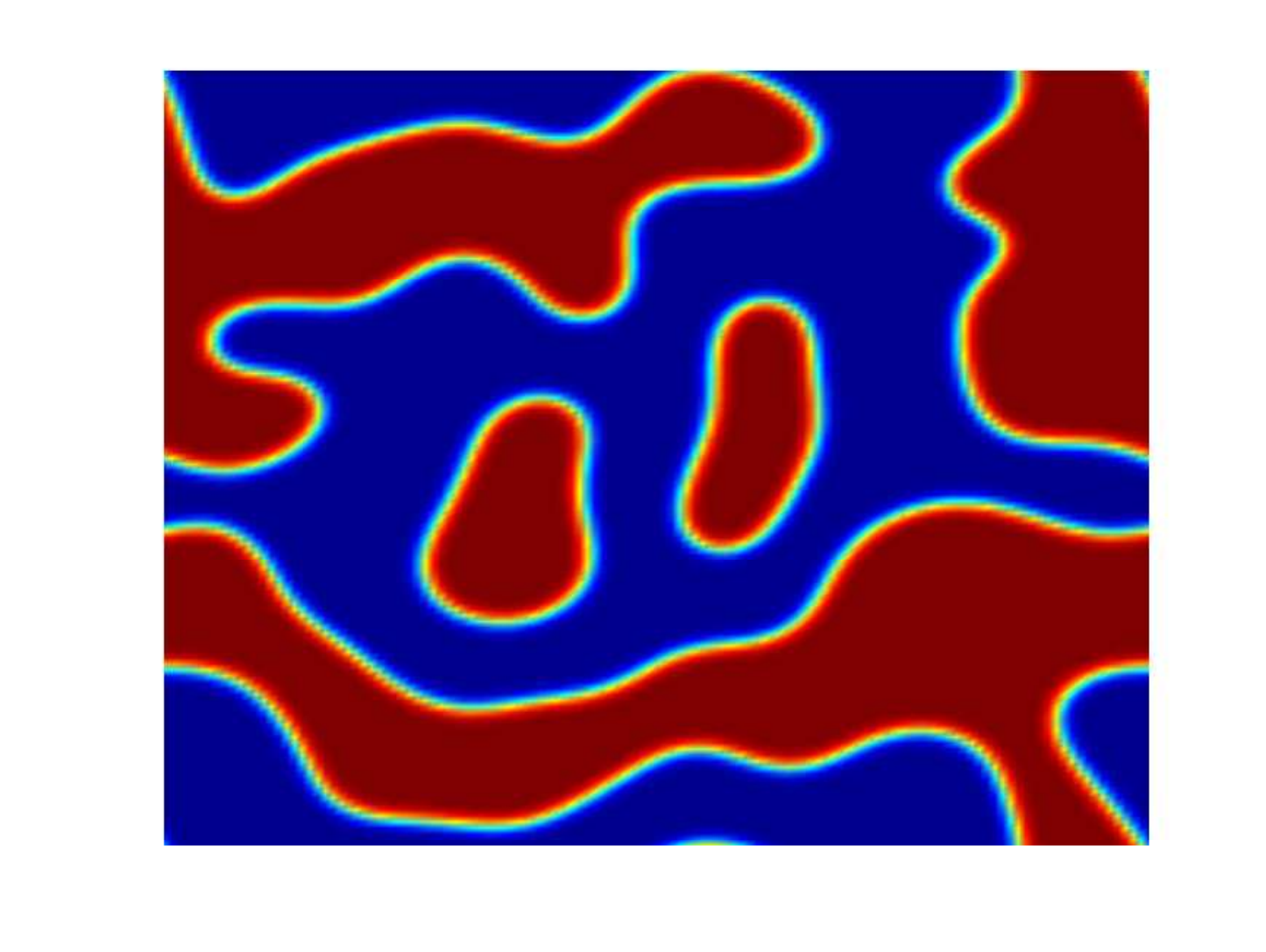}\\
\caption{Solution snapshots of coarsening dynamics of the
 conservative time-fractional Allen-Cahn equation at $t=1, 100, 300, 500$ (from left to right)
for three fractional orders $\alpha=0.4,\,0.7$ and $0.9$ (from top to bottom), respectively.}
\label{Nonlocal-AC-Coarsening-Dynamic}
\end{figure}
%%%%%%%%%%%%%%%%%%%%%%%%%%%%%%%%%%%%%%%%%%%%%%%%%%%%%%%%%%%%%%%%%%%%%%%%%%%%%%%%%%%%%

%%%%%%%%%%%%%%%%%%%%%%%%%%%%%%%%%%%%%%%%%%%%%%%%%%%%%%%%%%%%%%%%%%%%%%%%%%%%%%%%%%%%%%%%%%%%%
\begin{figure}[htb!]
\centering
\includegraphics[width=3.0in,height=2.0in]{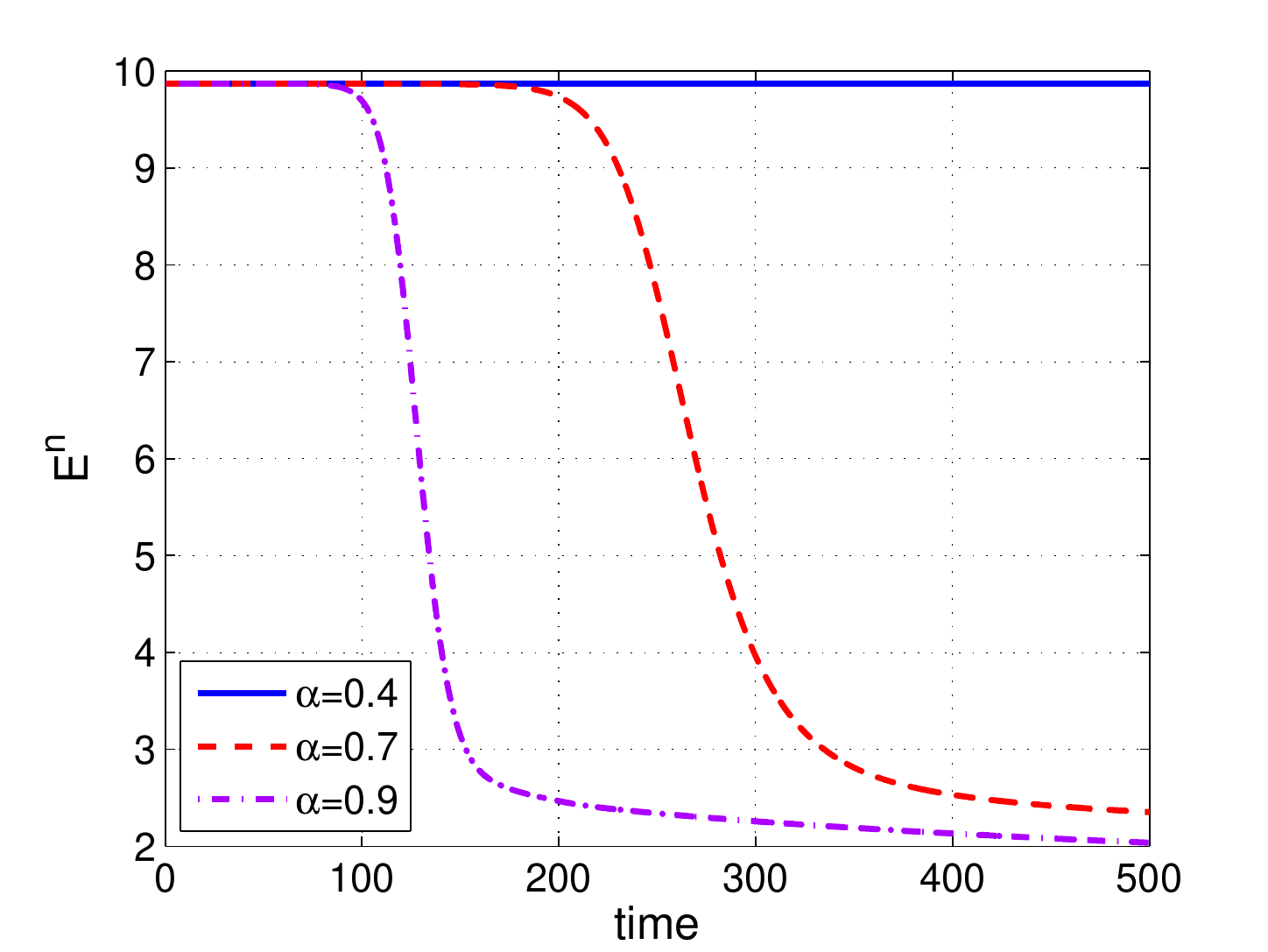}
\includegraphics[width=3.0in,height=2.0in]{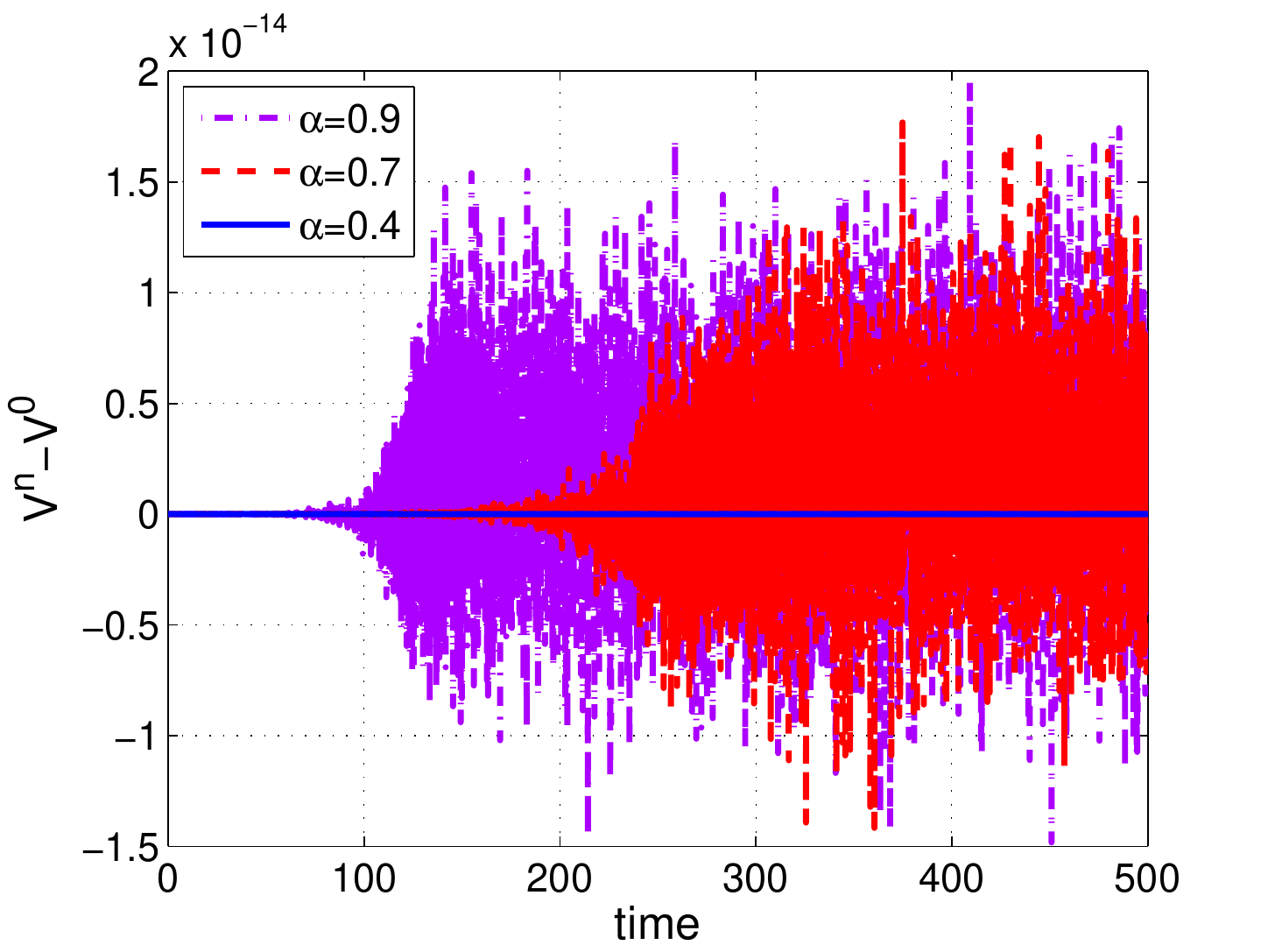}
\caption{Evolutions of energy and volume (from left to right) of
the conservative time-fractional Allen-Cahn equation for fractional orders $\alpha=0.4,\,0.7$ and $0.9$, respectively.}
\label{Nonlocal-AC-Dynamic-Energy-Mass}
\end{figure}
%%%%%%%%%%%%%%%%%%%%%%%%%%%%%%%%%%%%%%%%%%%%%%%%%%%%%%%%%%%%%%%%%%%%%%%%%%%%%%%%%%%%%
%%%%%%%%%%%%%%%%%%%%%%%%%%%%%%%%%%%%%%%%%%%%%%%%%%%%%%%%%%%%%%%%%%%%%%%%%%%%%%%%%%%%%

The coarsening snapshots of time-fractional Allen-Cahn equation with volume constraint are depicted
in Figure \ref{Nonlocal-AC-Coarsening-Dynamic}.
Compared with the numerical results in Figure \ref{CH-Coarsening-Dynamic},
these phase diagrams generated by the conservative time-fractional Allen-Cahn model
have no obvious difference with those produced by the time-fractional Cahn-Hilliard
model. Also, from Figure \ref{Nonlocal-AC-Dynamic-Energy-Mass},
we see that the new model \eqref{Problem-2} preserves the energy dissipation law
and the volume well. In summary,
the coarsening process of the time-fractional Allen-Cahn \eqref{Problem-2} with volume constraint
 is slower than that of the time-fractional Cahn-Hilliard model \eqref{Problem-3} because the energy dissipation rate of the former
is smaller; the mechanisms of coarsening dynamics of the two models are quite different
although both of them are volume-conserving.

\bibliographystyle{unsrt}%{siamplain}
\bibliography{Nonlocal-Allen-Cahn}
%%%%%%%%%%%%%%%%%%%%%%%%%%%%%%%%%%%%%%%%%%%%%%%%%%%%%%%%%%%%%%%%%%%%%%%%%%%%%%%%%%%%%%%%%%%
\end{document}